\documentclass[10pt]{iopart}
\usepackage{amssymb,mathrsfs}
\usepackage{amsbsy}
\usepackage{graphicx,rotating}
\usepackage[mathscr]{euscript}
\usepackage{color}
\usepackage{comment}
\usepackage{mathrsfs}
\usepackage{dsfont}
\usepackage{url}
\usepackage{tikz,pgfplots,pgfplotstable}
\usetikzlibrary{shapes}
\usepgfplotslibrary{external} 
\tikzexternaldisable

\expandafter\let\csname equation*\endcsname\relax
\expandafter\let\csname endequation*\endcsname\relax
\usepackage{amsmath}
\usepackage{mathtools}
\usepackage[pdfborder={0 0 0}, colorlinks=true, linkcolor=blue, urlcolor=blue]{hyperref}

\newcommand{\A}{\mathcal{A}}

\newcommand{\dd}{\mathbf{d}}

\newcommand{\Gm}{\mathbf{G}}
\newcommand{\Hess}{\mathcal{H}}

\newcommand{\hm}{\hat{m}}

\newcommand{\tm}{\tilde{m}}
\newcommand{\mm}{\mathbf{m}}
\newcommand{\nn}{\mathbf{n}}

\newcommand{\R}{\mathbb{R}}       %
 
\newcommand{\Rs}{\mathcal{R}}

\newcommand{\xx}{\mathbf{x}}

\newcommand{\pto}{\mathcal{F}}

\newcommand{\tv}{{\scriptscriptstyle\text{TV}}}
\newcommand{\tve}{{\scriptscriptstyle\text{TV},\varepsilon}}
\newcommand{\vtv}{{\scriptscriptstyle\text{VTV}}}

\newcommand{\cg}{{\scriptscriptstyle\text{cg}}}
\newcommand{\ncg}{{\scriptscriptstyle\text{ncg}}}
\newcommand{\nuc}{*}

\newcommand{\hRs}{\hat{\Rs}}

\def\addressices{Institute for Computational Engineering \& Sciences, The
  University of Texas at Austin, Austin, TX, USA}
\def\addressgeomech{Department of Geological Sciences and Department
  of Mechanical Engineering, The University of
  Texas at Austin, Austin, TX, USA}
\def\addressnyu{Courant Institute of Mathematical Sciences, New York University,
New York, NY, USA}

\DeclareMathAlphabet{\mathup}{OT1}{\familydefault}{m}{n}

\begin{document}

\title[Structural similarity and regularization for joint inverse PDE
  problems] {A comparative study of structural similarity and
  regularization for joint inverse problems governed by PDEs}
\author{Benjamin~Crestel$^1$, Georg~Stadler$^2$ and Omar~Ghattas$^{1,3}$}

\address{
$^1$\addressices\\
$^2$\addressnyu\\
$^3$\addressgeomech
}
\ead{\url{ben.crestel@utexas.edu}, 
\url{stadler@cims.nyu.edu} and
\url{omar@ices.utexas.edu}}
\vspace{10pt}

\begin{abstract}
Joint
inversion refers to the simultaneous inference of multiple parameter
fields from observations of systems governed by single or
multiple forward models. In many cases these parameter fields
reflect different attributes of a single medium and are thus spatially
correlated or structurally similar. By imposing prior information on
their spatial correlations via a joint regularization term, we seek to
improve the reconstruction of the parameter fields relative to
inversion for each field independently.
One of the main challenges is
to devise a joint regularization functional that conveys the spatial
correlations or structural similarity between the fields
while at the same time permitting scalable and efficient solvers for
the joint inverse problem.
We describe several joint regularizations that are motivated by these
goals: 
a cross-gradient and a normalized cross-gradient structural similarity
term, the vectorial total variation, and a joint regularization
based on the nuclear norm of the gradients. Based on numerical results
from three classes of inverse problems with piecewise-homogeneous
parameter fields, we conclude that the vectorial total variation functional is
preferable to the other methods considered.  Besides resulting in good
reconstructions in all experiments, it allows for scalable,
efficient solvers for joint inverse problems governed by PDE forward models.
\end{abstract}


\section{Introduction}
\label{sec:intro}

In a joint inverse problem one seeks to reconstruct multiple parameter
fields from observational data and  forward models that map the
parameter fields to the data. 
In many cases these parameter fields reflect different attributes of a
single medium and are thus spatially correlated or structurally
similar. By imposing prior information on their spatial correlations
via a joint regularization term, we seek to improve the reconstruction
of the parameter fields relative to inversion for each field
independently.

We formulate the joint inverse problem as an optimization problem with
a regularized data misfit objective, governed by a forward model
that represents a single or multiple physical phenomena. In the
following, we restrict ourselves to forward models that take the form
of partial differential equations (PDEs)
characterized by two unknown parameter fields, $m_1$ and $m_2$, which
we seek to reconstruct from observational data~$\dd$. 
The parameter-to-observable map $\pto(m_1,m_2)$
typically involves solution of the forward PDEs given the parameter
fields, followed by application of the observation operator, which
restricts the PDE solution to the space of observables. The
optimization problem is thus
\begin{equation} \label{eq:joint1} \min_{(m_1, m_2)} \left\{ \frac12 |
  \pto(m_1, m_2) - \dd |^2  + \Rs(m_1, m_2) \right\}.
\end{equation}
The role played by $\Rs$ in~\eqref{eq:joint1} is discussed in the next
paragraph.  Here, we address two specific settings
for~\eqref{eq:joint1}.  In the first, the forward model in
$\pto(m_1,m_2)$ describes a single physical phenomenon.  An example of
such a joint inverse problem is inversion for the primary and
secondary wave speeds in the Earth given measurements of the
acceleration at the surface. Obtaining high quality reconstructions
for both parameter fields is known to be difficult without
incorporating some form of prior knowledge that couples the two
fields~\cite{EpanomeritakisAkccelikGhattasEtAl08,
  ManukyanMaurerNuber16, LiLiangAbubakarEtAl13}.  We refer to
formulation~\eqref{eq:joint1} as a single physics joint inverse
problem.

In the second type of joint inverse problem, we consider observations
$\dd_1$ and $\dd_2$ stemming from two distinct physical phenomena
respectively, each depending on a single parameter field.  In this
case the forward models of the physical phenomena are uncoupled, and
coupling occurs only via the inverse problem.  The corresponding
parameter-to-observable maps are denoted by $\pto_1(m_1)$
and $\pto_2(m_2)$, %
resulting in
\begin{equation} \label{eq:joint2}
  \min_{(m_1, m_2)} \left\{ \frac12 |
  \pto_1(m_1) - \dd_1 |^2 + \frac12 | \pto_2(m_2) - \dd_2 |^2 + \Rs(m_1, m_2)
  \right\} .
\end{equation}
This formulation emerges from the general case above by defining
$\pto(m_1,m_2) = [\pto_1(m_1), \pto_2(m_2)]^T$ and $\dd = [\dd_1, \dd_2]^T$.
In the context of subsurface exploration, just a few of the different
physical phenomena that can be combined in~\eqref{eq:joint2} include
electromagnetic and seismic
waves~\cite{AbubakarGaoHabashyEtAl12,SemerciPanLiEtAl14}, radar and
seismic waves~\cite{FengRenLiuEtAl17}, DC resistivity and seismic
waves~\cite{GallardoMeju03}, and current resistivity and groundwater
flow~\cite{SteklovaHaber16}.
The {\em joint regularization term}~$\Rs(m_1,m_2)$
in~\eqref{eq:joint1} and \eqref{eq:joint2} acts to impose regularity
on $m_1$ and~$m_2$ individually to combat ill-posedness, but can also
express structural similarity or spatial correlations between the two
parameter fields. 
The remainder of this section introduces several different choices
for~$\Rs$.
To isolate regularization from structural similarity,
we decompose the joint regularization term~$\Rs(m_1,m_2)$ into
\[ \Rs(m_1,m_2) = \gamma_1 \Rs_1(m_1) + \gamma_2 \Rs_2(m_2) + \gamma
\hRs(m_1,m_2) , \] 
with $\gamma, \gamma_1, \gamma_2 > 0$.  The terms~$\Rs_1$ and~$\Rs_2$
are regularization terms for each parameter field; here we take them
to be total variation~(TV) regularizations, since our target media
are piecewise-homogeneous (i.e., blocky). The 
term~$\hRs(m_1,m_2)$ incorporates the structural similarity
between $m_1$ and $m_2$. 
We now discuss several choices for
$\hRs$.  In~\cite{GallardoMeju03}, the authors introduce the
{\em cross-gradient} term
\[ \hRs(m_1, m_2) = \frac12 \int_\Omega | \nabla m_1 \times \nabla m_2 |^2 \,
dx, \]
which seeks to align gradients of the two parameter fields at each
point in the medium, i.e., level sets that
have the same shape. 
This seems to be the most popular choice in
geophysics~\cite{AbubakarGaoHabashyEtAl12, SemerciPanLiEtAl14,
  FengRenLiuEtAl17, GallardoMeju03, SteklovaHaber16}, and is discussed
in section~\ref{sec:crossgradient}.  Instead of the gradients of the
parameter fields, one can use normalized gradients. This results in
the normalized cross-gradient term
\[ \hRs_\ncg(m_1,m_2) = \frac12 \int_\Omega \left| \frac{\nabla m_1}{|\nabla m_1|}
\times \frac{\nabla m_2}{|\nabla m_2|} \right|^2 \, dx. \]
The {\em normalized cross-gradient} was first used in the
context of image registration~\cite{HaberModersitzki06b}, and is
discussed in section~\ref{sec:normcg}.  Alternatively, when an
empirical relationship between both parameters is known, one could use
it in place of the structural similarity
term~$\hRs$~\cite{AbubakarGaoHabashyEtAl12,HaberHoltzman13}; this
approach, however, can be problematic in practice as these
relationships are typically uncertain (thus introducing bias) and the
resulting optimization problems can be difficult to
solve~\cite{HaberHoltzman13,GallardoMeju11}.

Alternatively, a single joint regularization term can 
impose regularity on both parameter fields while also expressing a
preference for structural similarity. In particular, we consider the
{\em vectorial total variation}~(VTV) functional,
\[ \Rs(m_1, m_2) = \gamma \int_\Omega \sqrt{|\nabla m_1|^2 + |\nabla m_2|^2} \,
dx , \]
with $\gamma > 0$.  The VTV functional was introduced in the
context of multi-channel imaging~\cite{BlomgrenChan98,BressonChan08},
and later used in PDE-constrained joint inverse
problems~\cite{HaberHoltzman13}; it is discussed in
section~\ref{sec:vtv}.  A second term we consider is the {\em nuclear
  norm}, which was used in \cite{Holt14,KnollHollerKoestersEtAl17} to promote
gradient alignment of a vector-valued image. 
Building on this idea, in section~\ref{sec:nn} we introduce
a nuclear norm-based joint regularization term for 
PDE-constrained joint inverse problems.

The objective of this article is to construct and assess joint
regularization terms that are (1) efficient for inverse problems
governed by PDEs with infinite-dimensional parameter fields (and are
thus large-scale after discretization) and (2) perform well in
reconstructing sharp interfaces in the truth parameter fields.
Indeed, targeting large-scale inverse problems entails several unique
challenges that limit choices of the joint regularization term.
Nonlinear inverse problems such as~\eqref{eq:joint1}
and~\eqref{eq:joint2} must be solved iteratively, which requires
gradient- (and Hessian-) based optimization methods to limit the
number of optimization iterations, along with adjoint methods to limit
the number of PDE model solutions that must be carried out at each
iteration.
Moreover, the adjoint method efficiently provides only directional
second derivatives rather than full Hessians, the construction of
which would require as many PDE solves as there are parameters (or
observations).
For these reasons, unless otherwise specified, we employ an inexact
Hessian-free Newton--conjugate gradient method with backtracking
line-search~\cite{EpanomeritakisAkccelikGhattasEtAl08,
  AkcelikBirosGhattas02, NocedalWright06}. That is, we compute the
Newton search direction using the preconditioned conjugate gradient
method, with early termination to guarantee a descent direction and to
avoid over-solving~\cite{DemboEisenstatSteihaug82}.  The efficient
solution of the Newton system depends crucially on the choice of
preconditioner; we detail our choices for each joint regularization
functional in sections~\ref{sec:cg}, \ref{sec:vtv}, and \ref{sec:nn}.
An overview of the numerical methods we employ to solve large-scale
inverse problems governed by PDEs can be found in~\ref{sec:numopt}.

Besides practicality and efficiency, our comparison of joint inversion
methods focuses on the quality of the reconstructions. Truth
parameter fields in geophysical exploration  and medical imaging
problems often present sharp contrasts within parameter fields.
We focus on joint
regularization terms that can best preserve sharp edges in the
reconstructed images.  Motivated by these criteria and a literature
review, we identified the four candidates discussed above, namely (1)
the cross-gradient and (2) its normalized variant, both paired with
individual TV regularizations, (3) the VTV joint regularization, and
(4) a nuclear norm-based joint regularization.

\subsection{Contributions}

The main contributions of this article are as follows:
(1)~We review three joint regularization terms commonly found in the
literature
(cross-gradient paired with TV, normalized cross-gradient paired
with TV, and VTV joint regularization), and discuss their practical use for
large-scale joint inverse problems governed by PDEs.
We derive their first and second derivatives, and use them to study
properties of the different joint regularization terms.
(2)~We adapt a nuclear norm joint regularization term to the context of joint inverse problems
governed by PDEs. We discuss some of the resulting computational challenges,
and propose a solver to address them.
(3)~We carry out a detailed comparison of all four joint regularization
terms over a broad range of applications,
and discuss their practical performance to reconstruct
parameter fields with sharp interfaces.

\subsection{Paper overview}

In the next three sections, we introduce the four joint regularization
terms. The cross-gradient and normalized cross-gradient are discussed
in sections~\ref{sec:crossgradient} and \ref{sec:normcg}, the
vectorial total variation in section~\ref{sec:vtv}, and the nuclear
norm joint regularization in section~\ref{sec:nn}.
Section~\ref{sec:numerics} summarizes our numerical experiments.  In
section~\ref{sec:poissonpoisson}, we report on several multiple
physics joint inverse
problems of the form \eqref{eq:joint2}, 
in which the two parameters fields arise as coefficients in
two independent Poisson equations, respectively.  We use this example
to illustrate some key features of each joint regularization term.
In section~\ref{sec:acoustic}, we consider a single physics joint
inverse problem of the form~\eqref{eq:joint1} for the acoustic wave
equation, in which we invert for the bulk modulus and the
density. Finally, in section~\ref{sec:poissonacoustic}, we study a
multiple physics joint inverse problem with two different forward
models, one an elliptic PDE and the other an acoustic wave equation.
Section~\ref{sec:ccl} provides concluding remarks.

\section{Cross-gradient terms}
\label{sec:cg}

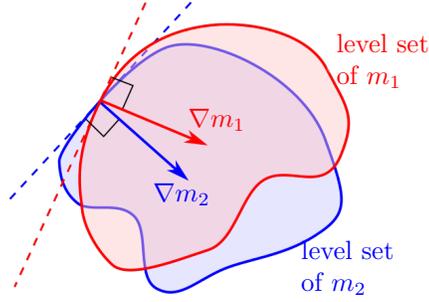
\begin{figure}[bt]
\centering
\tikzsetnextfilename{levelsets}
\begin{tikzpicture}
\node[anchor=south west, inner sep=0] (image1) at (0,0)
{\includegraphics[height=0.3\textwidth,trim=50 20 60 70, clip=true]
{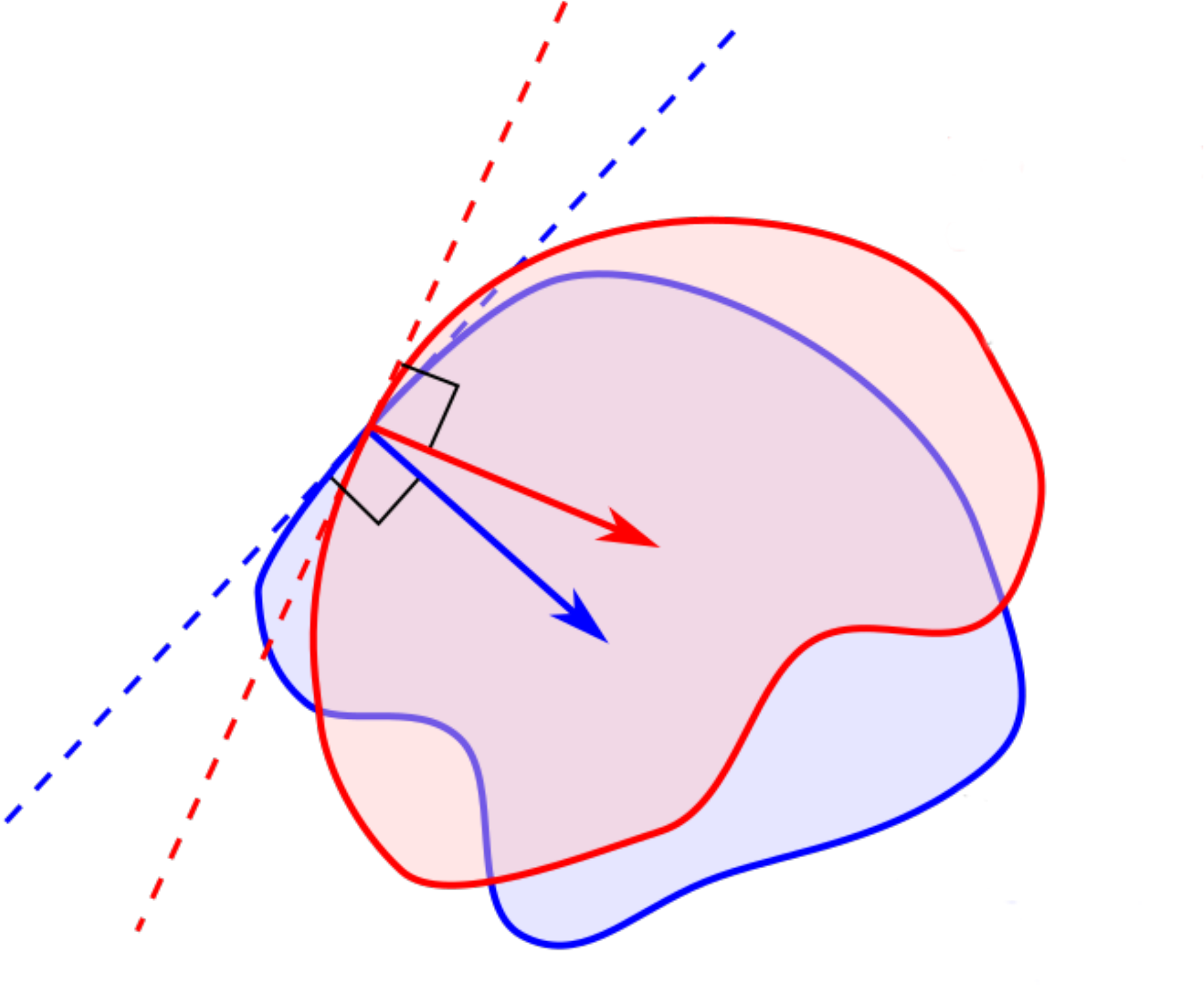}};
\begin{scope}[x={(image1.south east)},y={(image1.north west)}]
\node[red,thick] at (0.6,0.6) {$\nabla m_1$};
\node[red, text width=1.5cm] at (1.1,0.8) {level set of $m_1$};
\node[blue,thick] at (0.5,0.35) {$\nabla m_2$};
\node[blue, text width=1.5cm] at (1.0,0.1) {level set of $m_2$};
\end{scope}
\end{tikzpicture}
\caption{Sketch of a level set of the parameter fields $m_1$ (red) and
  $m_2$ (blue), with their respective gradients at a point.}
\label{fig:levelsets}
\end{figure}

In this section, we introduce the cross-gradient term and its normalized
version.  The main idea behind both of these structural similarity
terms is to express the preference that the
level sets of the inversion parameter fields $m_1$ and $m_2$ align. As 
illustrated in figure~\ref{fig:levelsets}, alignment of the level sets is
equivalent to the alignment of the gradients $\nabla m_1$ and
$\nabla m_2$ at each point.
By definition of the cross-product of two vectors, the vectors
$\nabla m_1$ and $\nabla m_2$ are aligned when $| \nabla m_1
\times \nabla m_2 |^2$ vanishes.

\subsection{The cross-gradient term}
\label{sec:crossgradient}

The cross-gradient term $\hRs_{\cg}$, defined as
\begin{equation} \label{eq:crossgrad}
  \hRs_{\cg}(m_1, m_2) \coloneqq \frac12 \int_\Omega
| \nabla m_1 \times \nabla m_2 |^2 \, dx, \end{equation}
was introduced in~\cite{GallardoMeju03} and has become a popular choice in
geophysical applications, particularly in seismic imaging.
Although the formulation~\eqref{eq:crossgrad} is intuitive, it is
inconvenient for discretization and computation of derivatives. Hence,
using vector calculus, we re-write \eqref{eq:crossgrad} as
\begin{equation} \label{eq:crossgrad_alt} \hRs_{\cg}(m_1, m_2) = \frac{1}2
\int_\Omega |\nabla m_1 |^2 | \nabla m_2 |^2 - (\nabla m_1 \cdotp \nabla m_2 )^2
\, dx .  \end{equation}
Combining the cross-gradient term~\eqref{eq:crossgrad_alt} with independent TV
regularizations for $m_1$ and $m_2$, we obtain the joint regularization
\begin{equation} \label{eq:rscg}
\Rs(m_1,m_2) = 
\gamma_1 \Rs_\tve(m_1) + \gamma_2 \Rs_\tve(m_2) + \gamma \hRs_\cg(m_1,m_2), 
\end{equation}
where and $\gamma,\gamma_1,\gamma_2>0$, and here and in the remainder
of this paper, we use the notation
\begin{equation}\label{eq:TVreg}
  \Rs_\tve(m) \coloneqq \int_\Omega \sqrt{|\nabla m|^2 + \varepsilon}
  \,dx \text{ for } \varepsilon > 0.
\end{equation}
In~\cite{HaberHoltzman13} the authors
propose a different formulation, in which each independent TV regularization is
weighted by a non-linear function of the gradient of the other parameter. The
goal of this weighting is to apply TV regularization only for points in the
parameter space where the cross-gradient term by itself is not sufficient to prevent
oscillatory solutions.  Such oscillations may occur where the gradient of one
parameter is very small, resulting in an (almost) vanishing cross-gradient term.
Because this formulation further increases the nonlinearity of the problem,
we instead use~\eqref{eq:rscg}.

Next, we derive first and second derivatives of the cross-gradient
regularization, interpret these derivatives as PDE operators, and draw
analogies with the derivatives of the TV functional $\Rs_{\tv}$ or its
regularized version \eqref{eq:TVreg}.  For this purpose, we first
derive the first and second variation of the TV functional as follows:
\begin{align*}
\delta_{m} \Rs_{\tv}(m;\tm) & =  \int_\Omega
|\nabla m|^{-1} (\nabla m \cdotp \nabla \tm) \, dx ,\\
\delta^2_{m} \Rs_{\tv}(m;\hm,\tm) & =  \int_\Omega 
|\nabla m|^{-1} (\nabla \hm \cdotp \nabla \tm) 
- |\nabla m|^{-3} (\nabla m \cdotp \nabla \tm) (\nabla m \cdotp \nabla \hm)  \, dx,
\end{align*}
where $\tm$ and $\hm$ are arbitrary directions.
Using integration by parts, the fact that $\tm$ in the expression
for $\delta^2_m$ is arbitrary, and the vector identity $(a\cdot
b)(c\cdot d) = b \cdot (a c^T) \cdot d$, one finds that the Hessian $\Hess$ is
the following second-order elliptic PDE operator
\begin{equation*}
\Hess_\tv \hm := -\nabla \cdotp (A_{\tv}(m) \nabla \hm),
\end{equation*}
with the anisotropic coefficient tensor
\begin{equation} \label{eq:difftv}
A_\tv(m) = \frac1{|\nabla m|} \left( I - \frac{\nabla m \nabla m^T}{|\nabla
  m|^2} \right).
\end{equation}
This interpretation as diffusion operator shows that $\Hess_\tv$ acts
very differently at different points $x\in \Omega$.  In particular,
let us consider a point $x$ where the norm of $\nabla m$ is large, e.g., $x$
is located at an interface in the parameter field $m$.  Then,
in directions orthogonal to $\nabla m$ (i.e., directions normal to
an interface), $A_\tv$ vanishes and thus the elliptic operator does not
smooth the reconstruction $m$ in these directions. In contrast, in
directions that are orthogonal to $\nabla m$ (i.e., directions that
are tangent to interfaces), $\A_\tv$ does not vanish, 
thus smoothing the reconstruction $m$ along interfaces.
This explains the anisotropic smoothing properties of the
TV functional and, in particular, its ability to recover sharp interfaces in
parameter fields. Away from interfaces, where $\nabla m$ is small,
$\A_\tv$ behaves like a scaled identity, thus smoothing $m$ in all
directions, much as $H^1$ norm-based Tikhonov regularization does.

We now turn to the derivation of the derivatives of the cross-gradient
term $\hRs_{\cg}$. Following similar arguments as for the scalar TV
regularization above, this will provide us with insight regarding the
regularization properties. Additionally, these derivatives are useful
for devising a Newton-type algorithm for the inverse problem solution
and for preconditioning the linear systems that arise.

Starting from \eqref{eq:crossgrad_alt}, we now compute the gradient, and
the action of the Hessian in a given direction for the cross-gradient
term. We perform the computations using weak forms
and then use integration by parts to derive the corresponding strong
forms. The directional derivative at $m \coloneqq (m_1, m_2)$ in a
direction~$\tm \coloneqq (\tm_1, \tm_2)$ is given by
\begin{align*}
\delta_{m_1} \hRs_{\cg}(m;\tm_1) & =  \int_\Omega
|\nabla m_2|^2 (\nabla \tm_1 \cdotp \nabla m_1) - (\nabla m_1 \cdotp \nabla m_2)
(\nabla \tm_1 \cdotp \nabla m_2) \, dx , \\
\delta_{m_2} \hRs_{\cg}(m;\tm_2) & =  \int_\Omega 
|\nabla m_1|^2 (\nabla \tm_2 \cdotp \nabla m_2) - (\nabla m_1 \cdotp \nabla m_2)
(\nabla \tm_2 \cdotp \nabla m_1) \, dx .
\end{align*}
Taking another variation, we find that the action of the Hessian of the
cross-gradient term in a direction~$\hm=(\hm_1, \hm_2)$ is
given by
\begin{align*}
\delta^2_{m_1} \hRs_{\cg}(m;\hm_1,\tm_1) & =  \int_\Omega 
|\nabla m_2|^2 (\nabla \tm_1 \cdotp \nabla \hm_1) 
-  (\nabla \tm_1 \cdotp \nabla m_2) (\nabla m_2 \cdotp \nabla \hm_1) \, dx,  \\
\delta^2_{m_1,m_2} \hRs_{\cg}(m;\hm_2,\tm_1) & =  \int_\Omega 
 2 (\nabla \tm_1 \cdotp \nabla m_1) (\nabla m_2 \cdotp \nabla \hm_2)
- (\nabla m_1 \cdotp \nabla m_2) (\nabla \tm_1 \cdotp \nabla \hm_2) \\
& \hspace{2in}
-  (\nabla \tm_1 \cdotp \nabla m_2) (\nabla m_1 \cdotp \nabla \hm_2) \, dx, \\
\delta^2_{m_2} \hRs_{\cg}(m;\hm_2,\tm_2) & =  \int_\Omega 
|\nabla m_1|^2 (\nabla \tm_2 \cdotp \nabla \hm_2) 
-  (\nabla \tm_2 \cdotp \nabla m_1) (\nabla m_1 \cdotp \nabla \hm_2) \, dx .
\end{align*}
In strong form and neglecting boundary conditions, the Hessian~$\Hess$
acts, in a direction~$\hm$, like an anisotropic vector diffusion operator,
i.e.,
\[ \Hess \hm = - \nabla \cdotp (A_{\cg}(m) \nabla \hm) , \]
where $A_{\cg}$ is a diffusion tensor
given by
\begin{equation} \label{eq:difftensorcg}
A_{\cg}(m) = \begin{bmatrix}
D(m_2) & B(m) \\
B(m)^T & D(m_1)
\end{bmatrix},
\end{equation}
with, for $i=1,2$,
\[ \begin{aligned}
D(m_i) & \coloneqq |\nabla m_i|^2 I - \nabla m_i \nabla m_i^T, \\
B(m) & \coloneqq 2 \nabla m_1 \nabla m_2^T - (\nabla m_1 \cdotp \nabla m_2) I -
\nabla m_2 \nabla m_1^T.
\end{aligned} \]
The block-diagonal part of~$A_\cg$ indicates a
TV-like behavior but where parameter~$m_1$ (resp.~$m_2$) preserves
interfaces in directions where parameter~$m_2$ (resp. $m_1$) presents
an interface; this illustrates the coupling between both parameters.
As we show numerically
in figure~\ref{fig:eigncg}, the Hessian of the cross-gradient term can be
indefinite.  The TV regularization being a convex functional, its Hessian is
guaranteed to be positive semidefinite.  Therefore, the Hessian obtained by
retaining the block diagonal parts of the diffusion
tensor~\eqref{eq:difftensorcg}, i.e., $\Hess_d \hm \coloneqq - \nabla \cdotp
(A_{\cg,d}(m) \nabla \hm)$, with
\begin{equation} \label{eq:acgdiag}
A_{\cg,d}(m) \coloneqq \begin{bmatrix} D(m_2) & 0 \\ 0 & D(m_1) \end{bmatrix},
\end{equation}
is also guaranteed to be positive semidefinite.  For this reason, when
using the cross-gradient paired with two independent TV
regularizations, we precondition the Newton system with a
block-diagonal matrix containing the Hessian of the TV
regularizations, combined with a small multiple of the identity in each block,
and the block-diagonal part of the Hessian
of the cross-gradient term~\eqref{eq:acgdiag}.

\newcommand{\wcb}{0.18\textwidth}
\begin{figure}[h]
\centering
\begin{tabular}{c@{}c@{}c@{\hspace{-.1in}}c@{}c}
\includegraphics[height=\wcb, trim=50 260 1150 260,
clip=true]{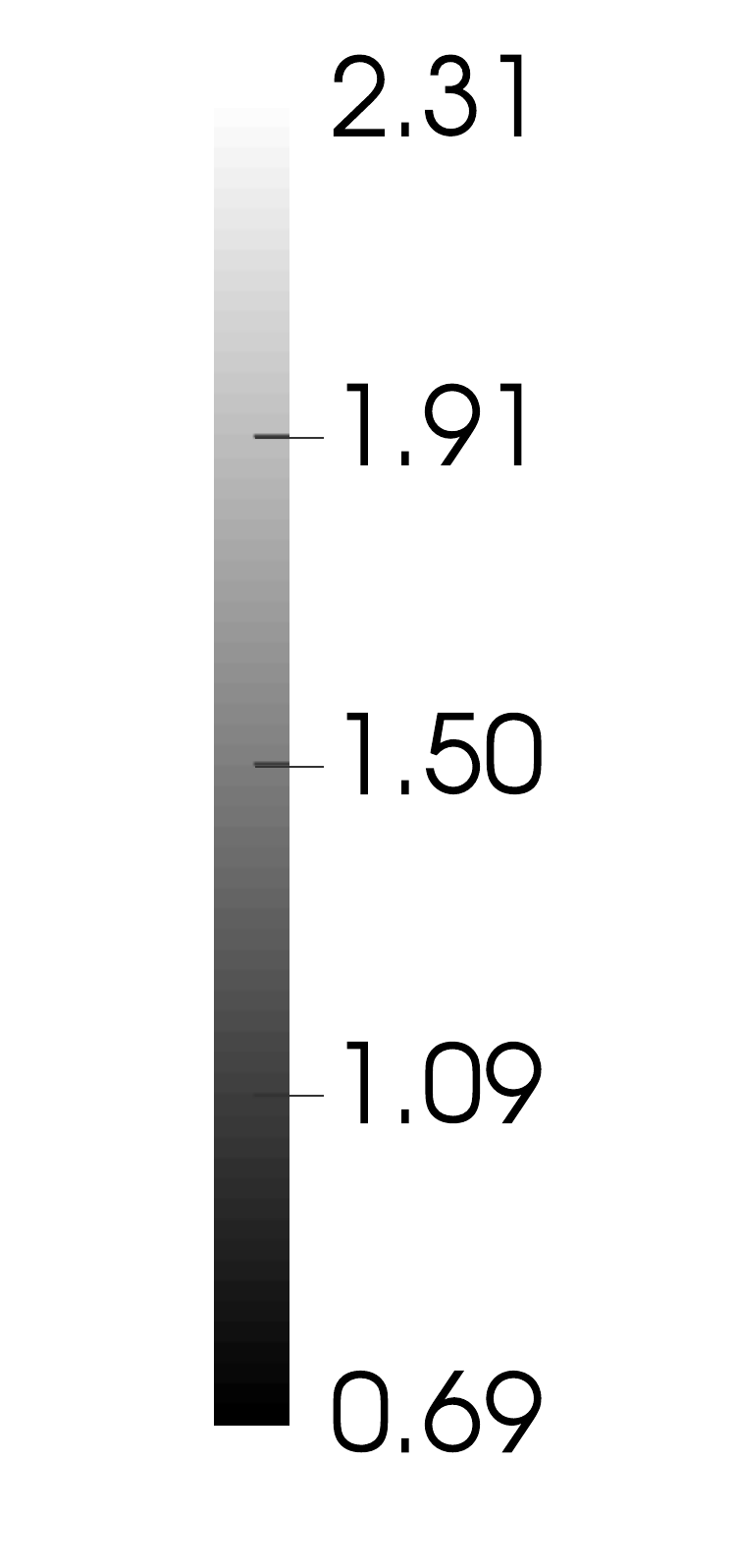}
&
\raisebox{.15\height}{
\includegraphics[height=\wcb, trim=290 70 290 70, clip=true]{fig/ncg/a1c1}
}
&
\raisebox{.15\height}{
\includegraphics[height=\wcb, trim=290 70 290 70, clip=true]{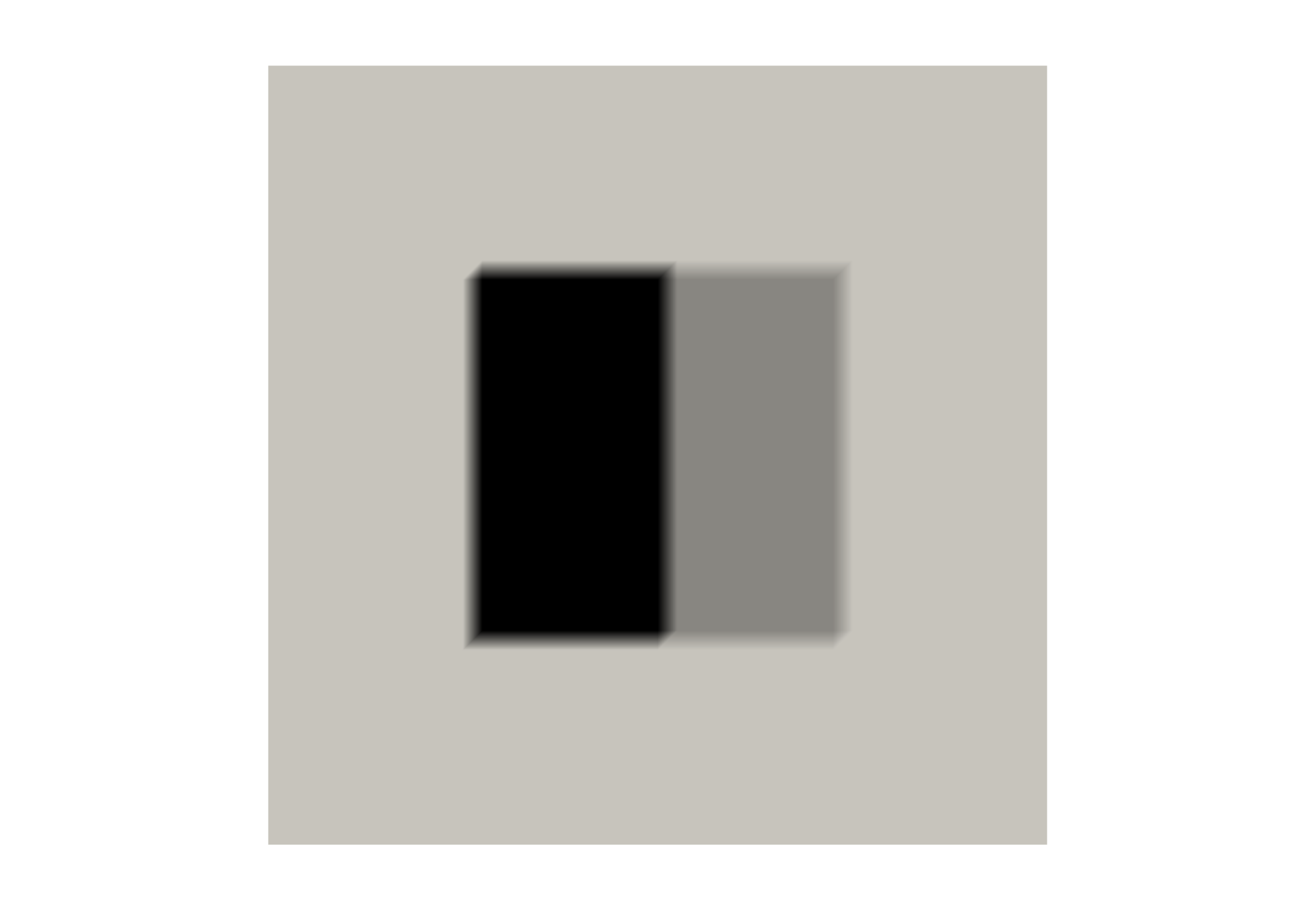}
}
&
\tikzsetnextfilename{cgH1}
\begin{tikzpicture}
\begin{axis}[height=\wcb, scale only axis,
legend style={font=\small}, legend pos=north east,
xtick pos=left,
xtick distance=2000,
ytick pos=bottom,
ytick={0, 10000},
scaled y ticks = true,
ymin=-2100, ymax=17000]
\addplot [color=blue, mark=none, thick] table [x expr=\coordindex, y index=0] {fig/ncg/eigcgH1.txt};
\addlegendentry{cg};
\addplot [color=red, mark=none, very thick, densely dashed] table [x expr=\coordindex, y index=0] {fig/ncg/eigcgHprecond1.txt};
\addlegendentry{cg precond};
\draw[thin] (axis cs:\pgfkeysvalueof{/pgfplots/xmin},0) -- (axis cs:\pgfkeysvalueof{/pgfplots/xmax},0);
\legend{};
\end{axis}
\end{tikzpicture}
&
\tikzsetnextfilename{ncgH1}
\begin{tikzpicture}
\begin{axis}[height=\wcb, scale only axis,
legend style={font=\small}, legend pos=north east,
xtick pos=left,
xtick distance=2000,
xlabel style={align=center},
ytick pos=bottom,
ytick={-30000,-20000,-10000, 0},
scaled y ticks = true,
ymin=-33000, ymax=1000]
\addplot [color=blue, mark=none, thick] table [x expr=\coordindex, y index=0] {fig/ncg/eigncgH1.txt};
\addlegendentry{ncg};
\addplot [color=red, mark=none, very thick, densely dashed] table [x expr=\coordindex, y index=0] {fig/ncg/eigncgHprecond1.txt};
\addlegendentry{ncg precond};
\draw[thin] (axis cs:\pgfkeysvalueof{/pgfplots/xmin},0) -- (axis cs:\pgfkeysvalueof{/pgfplots/xmax},0);
\legend{};
\end{axis}
\end{tikzpicture}
\\
&
\raisebox{.4\height}{
\includegraphics[height=\wcb, trim=290 70 290 70, clip=true]{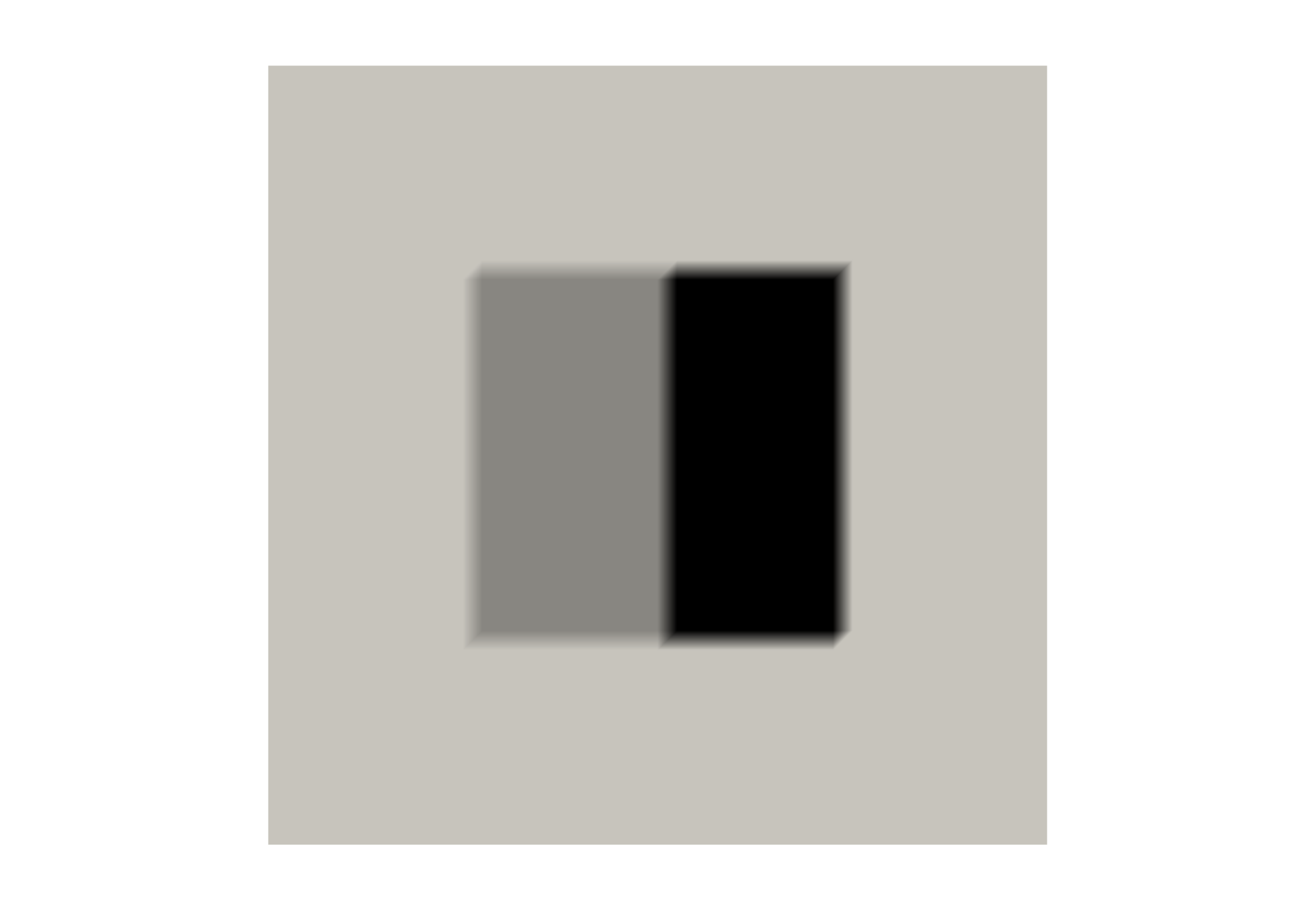}
}
&
\raisebox{.4\height}{
\includegraphics[height=\wcb, trim=290 70 290 70, clip=true]{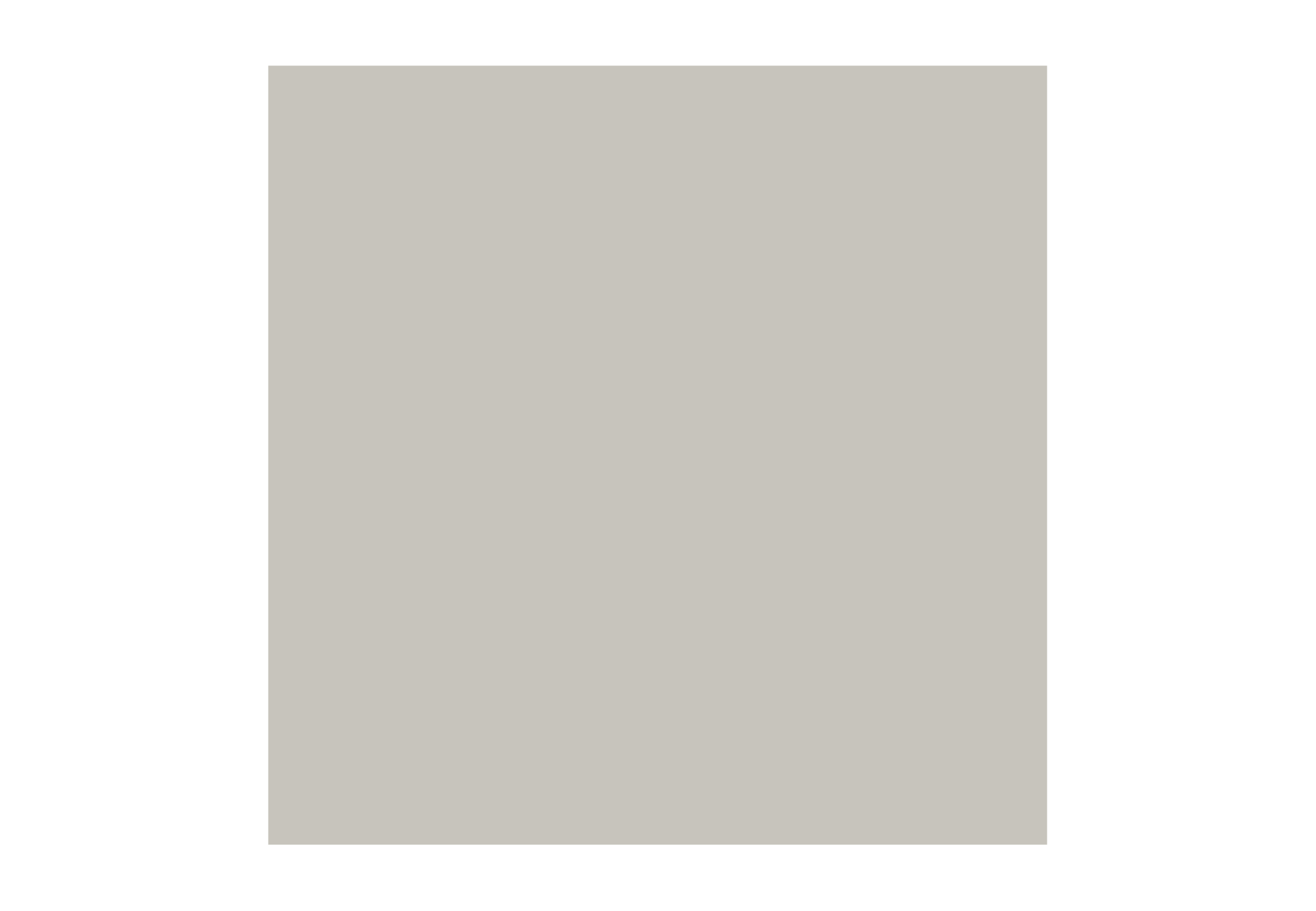}
}
&
\tikzsetnextfilename{cgH3}
\begin{tikzpicture}
\begin{axis}[height=\wcb, scale only axis,
legend style={font=\small}, legend pos=north east,
xtick pos=left,
xtick distance=2000,
ytick pos=bottom,
ytick={0, 10000},
scaled y ticks = true,
xlabel={eigenvalue},
ymin=-2100, ymax=17000]
\addplot [color=blue, mark=none, thick] table [x expr=\coordindex, y index=0] {fig/ncg/eigcgH3.txt};
\addlegendentry{cg};
\addplot [color=red, mark=none, very thick, densely dashed] table [x expr=\coordindex, y index=0] {fig/ncg/eigcgHprecond3.txt};
\addlegendentry{cg precond};
\draw[thin] (axis cs:\pgfkeysvalueof{/pgfplots/xmin},0) -- (axis cs:\pgfkeysvalueof{/pgfplots/xmax},0);
\legend{};
\end{axis}
\end{tikzpicture}
&
\tikzsetnextfilename{ncgH3}
\begin{tikzpicture}
\begin{axis}[height=\wcb, scale only axis,
legend style={font=\small}, legend pos=north east,
xtick pos=left,
xtick distance=2000,
xlabel style={align=center},
xlabel={eigenvalue},
ytick pos=bottom,
ytick={-30000,-20000,-10000, 0},
scaled y ticks = true,
ymin=-33000, ymax=1000]
\addplot [color=blue, mark=none, thick] table [x expr=\coordindex, y index=0] {fig/ncg/eigncgH3.txt};
\addlegendentry{ncg};
\addplot [color=red, mark=none, very thick, densely dashed] table [x expr=\coordindex, y index=0] {fig/ncg/eigncgHprecond3.txt};
\addlegendentry{ncg precond};
\draw[thin] (axis cs:\pgfkeysvalueof{/pgfplots/xmin},0) -- (axis cs:\pgfkeysvalueof{/pgfplots/xmax},0);
\legend{};
\end{axis}
\end{tikzpicture}
\\
 & (i) $m_1$ & (ii) $m_2$ & (iii) cross-gradient & (iv) norm.cross-gd
\end{tabular}
\caption{Eigenvalues of the Hessian operator (blue) and block-diagonal part of
the Hessian operator (red) for the (iii) cross-gradient
term~\eqref{eq:crossgrad} and the (iv) normalized cross-gradient
term~\eqref{eq:ncg} with $\varepsilon=10^{-4}$, for two combinations of truth
parameter fields (i) $m_1$ and (ii) $m_2$.  The domain is a unit square
discretized by a $40 \times 40$ mesh of squares subdivided
into triangles, and the parameter fields $m_1$ and $m_2$
are discretized using continuous piecewise linear finite elements.}
\label{fig:eigncg}
\end{figure}

\subsection{Normalized cross-gradient}
\label{sec:normcg}

A disadvantage of the cross-gradient term~\eqref{eq:crossgrad_alt} is
that it vanishes where one of the inversion parameter fields is
constant, hence potentially ignoring sharp discontinuities in the
other.
A remedy, proposed in the context of image
registration in~\cite{HaberModersitzki06b}, is to normalize the gradient of both
inversion parameters in the formulation of the cross-gradient.  The normalized
cross-gradient is given by
\[ \hRs(m_1,m_2) = \frac12 \int_\Omega \left| \frac{\nabla
m_1}{|\nabla m_1|} \times \frac{\nabla m_2}{|\nabla m_2|} \right|^2 \, dx
= \frac12 \int_\Omega 1 - \left( \frac{\nabla m_1 \cdotp
\nabla m_2}{|\nabla m_1| |\nabla m_2|} \right)^2 \, dx. \]
Since this formulation is non-differentiable where $|\nabla m_1|=0$
or $|\nabla m_2|=0$, we use the modified normalized cross-gradient,
\begin{equation} \label{eq:ncg} 
\hRs_\ncg(m_1,m_2) \coloneqq \frac12 \int_\Omega 1 - \left( \frac{\nabla m_1 \cdotp
\nabla m_2}{\sqrt{|\nabla m_1|^2 + \varepsilon} \sqrt{|\nabla m_2|^2 +
\varepsilon}} \right)^2 \, dx,
\end{equation}
with $\varepsilon > 0$. In the rest of this paper, we refer to~\eqref{eq:ncg}
when discussing the normalized cross-gradient.
Combining the normalized cross-gradient term~\eqref{eq:ncg} with two TV
regularizations, we obtain the joint regularization
\begin{equation} \label{eq:rsncg}
\Rs(m_1,m_2) = 
\gamma_1 \Rs_\tve(m_1) + \gamma_2 \Rs_\tve(m_2) + \gamma \hRs_\ncg(m_1,m_2),
\end{equation}
where $\gamma,\gamma_1,\gamma_2>0$.
Compared to the cross-gradient term, the derivatives of the normalized
cross-gradient term give less obvious insight into its regularization
behavior.  Instead, we illustrate numerically that the normalized
cross-gradient often behaves as a concave operator. In
figure~\ref{fig:eigncg}, we plot the eigenvalues of its Hessian
and of the block-diagonal part of its Hessian for different parameter
fields $m_1$ and $m_2$, and observe that most eigenvalues are
negative.  The main
practical consequence of this observation is that the Hessian of the
joint regularization~\eqref{eq:rsncg} may be indefinite.  For this
reason, the preconditioner for the Newton system is formed by the
Hessians of the TV regularizations alone.

\section{Vectorial total variation}
\label{sec:vtv}

The vectorial total variation functional~\cite{BressonChan08}, or color
TV~\cite{BlomgrenChan98}, is the multi-parameter equivalent of the total
variation functional. It was first introduced for multi-channel imaging
applications~\cite{BlomgrenChan98,BressonChan08}, and later applied to joint
inverse problems~\cite{HaberHoltzman13}.
The VTV functional is convex, and unlike the
cross-gradient and normalized cross-gradient, it serves as a regularization by
itself, i.e., it does not require additional regularization terms.  It is given by
\begin{equation} \label{eq:vtvorig}
\Rs(m_1, m_2) = \gamma \int_\Omega \sqrt{|\nabla m_1|^2 + |\nabla m_2|^2} \, dx,
\end{equation}
with $\gamma > 0$. Since this formulation is
non-differentiable where 
$|\nabla m_1| = |\nabla m_2|=0$, we introduce a modified VTV regularization given by
\begin{equation} \label{eq:VTV}
\Rs_{\vtv}(m_1, m_2) \coloneqq \gamma \int_\Omega \sqrt{|\nabla m_1|^2 + |\nabla
m_2|^2 + \varepsilon} \, dx ,
\end{equation}
with $\varepsilon, \gamma > 0$.
Whereas the cross-gradient terms (see section~\ref{sec:cg}) work by
aligning the level sets of the inversion parameter fields, VTV 
favors superimposition of discontinuities. An
intuitive way to explain this, given the understanding of the TV
regularization~\cite{ChambolleCasellesCremersEtAl10}, is sketched in
figure~\ref{fig:vtvexplain}.  Given two parameter fields with a single
jump of same amplitude, the VTV functional is minimum when both jumps
occur at the same location.
\begin{figure}[h]
\centering
\begin{tabular}{r@{}ccc}
&
\tikzsetnextfilename{vtv}
\begin{tikzpicture}
\begin{axis}[width=0.2\textwidth, scale only axis,
xtick pos=left, ytick pos=bottom]
\addplot+[const plot, no marks, thick, blue] coordinates {(0,0) (0.7,1) (2,1)}
node[above, pos=0.3, blue] {$\color{blue} m_1$};
\addplot+[const plot, no marks, thick, red, dashed] coordinates {(0,0) (1.3,1)
(2,1)} node[below, pos=0.8, red] {$\color{red} m_2$};
\end{axis}
\end{tikzpicture}
& &
\tikzsetnextfilename{vtv2}
\begin{tikzpicture}
\begin{axis}[width=0.2\textwidth, scale only axis,
xtick pos=left, ytick pos=bottom]
\addplot+[const plot, no marks, thick, blue] coordinates {(0,0) (0.9,1) (2,1)};
\addplot+[const plot, no marks, thick, red, dashed] coordinates {(0,0) (0.9,1) (2,1)};
\end{axis}
\end{tikzpicture}
\\
$\Rs_\vtv(m_1,m_2) = $ &
\hspace{.2in}
''$2 \int_\Omega \sqrt{|\nabla m_{1}|^2} \, dx$'' &
$>$ &
''$\sqrt{2} \int_\Omega \sqrt{|\nabla m_{1}|^2} \, dx$''
\end{tabular}
\caption{Values of the VTV
regularization~\eqref{eq:vtvorig}, for two parameter fields $m_1$ and $m_2$
defined over $\Omega = [0,2]$, with both parameter fields having a single
jump of the same amplitude, and $\Rs_\tv(m_1) = \Rs_\tv(m_2)$.
This informal argument can be made rigorous by using piecewise linear functions for
$m_1$ and $m_2$.}
\label{fig:vtvexplain}
\end{figure}

The derivatives of the VTV regularization resemble those of the TV
regularization. For simplicity, we set $\gamma \equiv 1$ in the rest
of this section.  The directional derivative at a point~$m = (m_1,
m_2)$ in a direction~$\tm = (\tm_1, \tm_2)$ is given by
\begin{equation} \label{eq:VTVgrad} \begin{aligned}
\delta_{m_i} \Rs_\vtv(m;\tm_i) & = \int_\Omega \frac{\nabla m_i \cdotp \nabla \tm_i}{\sqrt{|\nabla m_1|^2 +
|\nabla m_2|^2 + \varepsilon}} \, dx, \quad \text{ for } i=1,2. %
\end{aligned} \end{equation}
We again interpret the Hessian of the VTV as a diffusion tensor to study
its anisotropic diffusion behavior.
In strong form (see section~\ref{sec:crossgradient}), it is given by
\begin{equation} \label{eq:AVTV}
 A_\vtv(m) := 
\frac1{|\nabla m|_\varepsilon}
\begin{bmatrix}
I - \frac{\nabla m_1 \nabla m_1^T}{|\nabla m|^2_\varepsilon} &
- \frac{\nabla m_1 \nabla m_2^T}{|\nabla m|^2_\varepsilon} \\
- \frac{\nabla m_2 \nabla m_1^T}{|\nabla m|^2_\varepsilon} &
I - \frac{\nabla m_2 \nabla m_2^T}{|\nabla m|^2_\varepsilon}
\end{bmatrix} ,
\end{equation}
where $|\nabla m|^2_\varepsilon = |\nabla m_1|^2 + |\nabla m_2|^2 +
\varepsilon$.  Comparing with the diffusion tensor for the Hessian of the TV
regularization~\eqref{eq:difftv}, we find similar terms along the
block diagonal, with the exception of the normalization factor in the denominator.
It is $|\nabla m_i|$ in the case of TV, and $|\nabla m|_\varepsilon$ in the case
of VTV, i.e., it involves the gradient of both parameters, hence introducing
coupling between the parameter fields.
The eigen-decomposition of the diffusion tensor of the Hessian
provides further insights.
For simplicity, we use $\varepsilon=0$ in this analysis.
Skipping details that can be found in~\cite{Crestel17}, the eigenpairs for the diffusion
tensor~$A_\vtv$ are
\[ \left( \begin{bmatrix} \nabla m_1 \\ \nabla m_2 \end{bmatrix} , 0 \right), 
\left( \begin{bmatrix} (\nabla m_1)^\perp \\ 0 \end{bmatrix} , \frac1{|\nabla \mm|} \right) , 
\left( \begin{bmatrix} 0 \\ (\nabla m_2)^\perp \end{bmatrix} , \frac1{|\nabla \mm|} \right) , 
\left( \begin{bmatrix} \nabla m_2 \\ - \nabla m_1
\end{bmatrix} , \frac1{|\nabla \mm|} \right) . \]
The kernel of the diffusion tensor contains parameter field directions
that are not smoothed out by the regularization. Reconstructions in
these directions can
display sharp edges.
It is informative to compare the eigenpairs of the diffusion tensor
arising from the VTV Hessian with those arising from %
$%
\Rs_\tv(m_1) + \Rs_\tv(m_2)$,
the sum of two independent TV regularizations.
In this case, the eigenpairs are
\[ \left( \begin{bmatrix} \nabla m_1 \\ 0 \end{bmatrix} , 0 \right), 
\left( \begin{bmatrix} 0 \\ \nabla m_2 \end{bmatrix} , 0 \right), 
\left( \begin{bmatrix} (\nabla m_1)^\perp \\ 0 \end{bmatrix} , \frac1{|\nabla m_1|} \right) ,
\left( \begin{bmatrix} 0 \\ (\nabla m_2)^\perp \end{bmatrix} , \frac1{|\nabla m_2|} \right) . \]
The sum of independent TV regularizations acts in the direction of each
parameter~$m_i$ independently from the other parameters, analogously
to the TV functional for a single inverse problem.  That is, it
preserves sharp interfaces in the parameter fields (large values
of~$|\nabla m_i|$) but smoothes along interfaces.
This is in contrast with the kernel of the diffusion tensor of VTV,
which favors parameter fields with sharp variations occurring at the
same physical locations.

The use of TV regularization in PDE-constrained inverse problems
increases the
nonlinearity of the problem, and requires the use of customized
solvers.  Due to the similarity between TV and VTV, a similar
challenging numerical behaviour can be expected for VTV.
In~\cite{Crestel17} we tailor a primal-dual Newton method
\cite{HintermullerStadler06} for the efficient, scalable solution of
PDE-constrained joint inverse problems regularized with VTV.  Since
the focus of the current paper is on a qualitative comparison of
several joint regularization terms, we skip details of this solver
here and instead refer to~\cite{Crestel17}.

\section{Nuclear norm joint regularization}
\label{sec:nn}

The nuclear norm joint regularization
seeks to promote gradient alignment by minimizing the rank of the Jacobian of the gradients of the
parameter fields.
Different versions of that idea have been used in various imaging applications.
In color image denoising, this approach is often referred to as total nuclear
variation~\cite{Holt14}; the unified framework to discuss
VTV and the total nuclear variation in \cite{Holt14} shows that the
nuclear norm-based functional is a regularizer in itself.
This can be simply justified by the equivalence of all norms in
finite dimensions (here, on the space of matrices).
In~\cite{KnollHollerKoestersEtAl17}, the authors
propose the pointwise nuclear norm of a matrix field as regularization to
express a preference for alignment of image edges. 
Building on~\cite{Holt14,KnollHollerKoestersEtAl17}, we propose a 
nuclear norm joint regularization suitable for large-scale PDE constrained optimization.

As for the methods discussed in section~\ref{sec:cg}, this term seeks to
promote alignment of parameter level sets by attaining its minimum
value when gradients align.  
Let us introduce the matrix-valued function~$\Gm:
\Omega \rightarrow \R^{d \times 2}$, with $\Omega \subset \R^d$ the
physical domain, defined by
\[ \Gm(x) \coloneqq \big[ \nabla m_1 | \nabla m_2 \big]
= \begin{bmatrix}
\partial_{x_1} m_1 & \partial_{x_1} m_2 \\
\vdots & \vdots \\
\partial_{x_d} m_1 & \partial_{x_d} m_2 
\end{bmatrix} . \]
The gradients $\nabla m_1$ and $\nabla m_2$ are aligned at $x \in
\Omega$ if the columns of~$\Gm(x)$ are multiples of each other, in
which case the rank of $\Gm(x)$ is~1.  One could seek to promote
gradient alignment by minimizing $\int_\Omega \text{rank} (\Gm(x)) \,
dx$.  However, in practice, minimization of the rank of a matrix is
notoriously difficult. The nuclear norm of a matrix, defined as the
$\ell_1$-norm of its singular values and denoted by
$\| \cdotp \|_\nuc$, is often a good proxy for the
rank~\cite{Fernandez-Granda16}.
We therefore define, with $\gamma>0$, the nuclear norm joint
regularization as
\begin{equation}
\label{eq:nucl}
\Rs_\nuc(m_1,m_2) \coloneqq \gamma \int_\Omega \| \Gm(x) \|_\nuc \,
dx.
\end{equation}

\subsection{Gradient of the nuclear norm joint regularization}

We now compute derivatives of~\eqref{eq:nucl} using the chain rule.
Let us introduce the notation $f(M) \coloneqq \| M \|_\nuc$, for
arbitrary $M \in \R^{d \times 2}$.  Thus, $\Rs_\nuc(m_1,m_2) = \gamma
\int_\Omega f(\Gm(x)) \, dx$.  Denoting the gradient of~$f$ with
respect to the entries of matrix~$M$ by $\nabla f(M) \in \R^{d \times
  2}$, the first directional derivatives of~\eqref{eq:nucl} with
respect to the inversion parameters $m_i$, $i=1,2$, in a
direction~$\tm_i$, are given by
\begin{equation} \label{eq:derivnn}
\partial_{m_i} \Rs_\nuc(m_1, m_2) \tm_i  = \gamma \int_\Omega ( \nabla f(\Gm),
\partial_{m_i} \Gm(x) \tm_i ) \, dx, \end{equation}
where
\[ \partial_{m_1} \Gm(x) \tilde{m}_1  
= \begin{bmatrix}
\partial_{x_1} \tm_1 & 0 \\
\vdots & \vdots \\
\partial_{x_d} \tm_1 & 0 
\end{bmatrix}  \text{ and } 
 \partial_{m_2} \Gm(x) \tilde{m}_2  
= \begin{bmatrix}
0 & \partial_{x_1} \tm_2  \\
\vdots & \vdots \\
0 & \partial_{x_d} \tm_2  
\end{bmatrix} ,
\]
and the inner product for matrices $M = (m_{ij})_{ij},N=(n_{ij})_{ij} \in \R^{d
\times 2}$ is defined as $(M,N) = \sum_{i=1}^d \sum_{j=1}^2 m_{ij} n_{ij}$.

We next compute the gradient of the nuclear norm~$\nabla f(M)$.  Given
a full-rank matrix $M \in \R^{n \times m}$, i.e., $r \coloneqq
\text{rank}(M) = \min(m,n)$, and singular values $\{ \sigma_k
\}_{k=1}^r$, we define its (reduced) singular value
decomposition~(SVD) by $M = U \Sigma V^T$, with $U \in R^{n \times
  r}$, $V \in \R^{m \times r}$, and $\Sigma \in \R^{r \times r}$ a
diagonal matrix containing the singular values of~$M$, i.e.,
$\Sigma_{kk}=\sigma_k > 0$, $k=1, \ldots, r$. The $(i,j)$-entry of the
gradient of the nuclear norm is given by
\[ \left( \nabla f(M) \right)_{ij} 
= \sum_{k=1}^r \frac{\partial \sigma_k}{\partial m_{ij}} 
= \sum_{k=1}^r u_{ik}v_{jk} , \]
where the second equality uses the singular value
sensitivity~\cite{PapadopouloLourakis00}.  The gradient of the nuclear norm with
respect to the entries of~$M$ is then given by
\[ \nabla f(M) = U V^T . \]

\subsection{Modified nuclear norm joint regularization}

The nuclear norm~$f(M)$ is non-differentiable when the matrix~$M$ is not full-rank,
corresponding to the case where at least one of the singular values
vanishes.
To make it differentiable, similar to the treatment of TV
regularization, we define 
the modified nuclear norm by
\begin{equation} \label{eq:nuclreg}
f_\varepsilon(M) \coloneqq
\| M \|_{\nuc,\varepsilon} = \sum_{k=1}^{\min(m,n)} \sqrt{\sigma_k^2 +
\varepsilon} , \end{equation}
where $\varepsilon > 0$. For $\gamma>0$, we define the modified
nuclear norm joint regularization as
\begin{equation}
\label{eq:nucl2}
\Rs_{\nuc,\varepsilon}(m_1,m_2) \coloneqq \gamma \int_\Omega f_\varepsilon(
\Gm(x) ) \, dx.
\end{equation}

The $(i,j)$-entry of the gradient of the modified
nuclear norm~\eqref{eq:nuclreg} is given by
\[ \left( \nabla f_{\varepsilon}(M) \right)_{ij}  = 
\frac{\partial}{\partial m_{ij}} \sum_{k=1}^{\min(m,n)} \sqrt{\sigma_k^2 +
\varepsilon} 
 = \sum_{k=1}^{r} \frac{\sigma_k}{\sqrt{\sigma_k^2 + \varepsilon}} 
\frac{\partial \sigma_k}{\partial m_{ij}} , \]
where in the last expression the sum is up to~$r$ since, by definition of the
rank of a matrix, $\sigma_k=0$ for all $k > r$.  Let us now introduce the
diagonal matrix~$W^{}_\varepsilon \in \R^{r \times r}$, with entries
$(W^{}_\varepsilon)_{ii} = \sigma_i/\sqrt{\sigma_i^2 + \varepsilon}$.  Using the
expression for the sensitivity of the singular
values~\cite{PapadopouloLourakis00}, the gradient of the modified nuclear
norm is then given by
\begin{equation} \label{eq:gradmnn}
\nabla f_{\varepsilon}(M) = U W^{}_\varepsilon V^T . \end{equation}
The first directional
derivatives of~\eqref{eq:nucl2} with respect to the inversion parameters $m_i$,
$i=1,2$, in a direction~$\tm_i$, are given by
\begin{equation} \label{eq:derivnn2}
\partial_{m_i} \Rs_{\nuc,\varepsilon}(m_1, m_2) \tm_i  = \gamma \int_\Omega ( \nabla
f_\varepsilon(\Gm), \partial_{m_i} \Gm(x) \tm_i ) \, dx, \end{equation}

The modified nuclear norm~\eqref{eq:nuclreg}, however, is not twice
differentiable when two singular values are equal (crossing singular
values). This is because the second derivative requires the
sensitivity of the individual singular vectors, which are not
differentiable where singular values cross.  We have not found a
practical workaround for this singularity, and thus proceed with a
gradient-based method to solve joint inverse problems regularized with
the nuclear norm joint regularization; the solver is detailed in \ref{sec:solvernn}.
In the rest of this paper, when using ``nuclear norm joint regularization'', we
refer to the modified nuclear norm joint regularization~\eqref{eq:nucl2}.

\section{Numerical examples}
\label{sec:numerics}

In this section, we present a comprehensive numerical comparison of the four
joint regularization approaches introduced in sections~\ref{sec:cg}--\ref{sec:nn},
i.e., the cross-gradient \eqref{eq:rscg}, the normalized cross-gradient
\eqref{eq:rsncg}, the vectorial
total variation~\eqref{eq:VTV}, and the nuclear norm
\eqref{eq:nuclreg} regularization.  Reconstructions obtained with these joint
regularization terms are compared with each other, and with the reconstructions
obtained by solving a joint inverse problem with independent TV
regularizations. 
The parameters for all joint regularization terms are selected empirically as
leading to the best reconstructions. The values of~$\varepsilon$ are chosen small
enough to provide reconstructions with sharp interfaces, but large enough to
avoid numerical difficulties (see for instance the discussion
in~\cite{AscherHaberHuang06}).

The different regularizations are compared using three examples
covering both types of joint inverse problems~\eqref{eq:joint1}
and~\eqref{eq:joint2}. 
In section~\ref{sec:poissonpoisson}, we combine two uncoupled Poisson inverse
problems to form the joint inverse problem~\eqref{eq:joint2}, where we 
invoke prior knowledge that the two truth
parameter fields have similar structure.  
In section~\ref{sec:acoustic}, we compare 
the ability of the 
joint regularization
terms to improve the reconstruction of the bulk modulus and the density in
an acoustic wave equation, an example of a joint inverse
problem~\eqref{eq:joint1}.
Finally, in section~\ref{sec:poissonacoustic}, we formulate a multi-physics
joint inverse problem~\eqref{eq:joint2}, which combines an inverse
problem governed by the Poisson equation with one governed by
the acoustic wave equation.
Here again, the Poisson parameter and the wave speed fields are assumed to
have similar structure.

In all examples, the domain is a 2D unit square, 
with a uniform mesh of
isosceles right
triangles obtained by cutting in half $N \times N$ squares; we define the
mesh size parameter~$h \coloneqq 1/N$.
All data are generated synthetically from the truth parameter fields, and
then polluted by adding independent and identically distributed Gaussian noise;
the noise level is specific to each example.  We use continuous Galerkin finite elements
to discretize all field variables, with the state,
adjoint, incremental state, and incremental adjoint variables using quadratic
elements, and the parameter fields using linear elements.  
All examples are implemented in Python and build on the finite element
library FEniCS~\cite{LoggMardalGarth12,LoggWells10}.  For the examples
in section~\ref{sec:poissonpoisson} and~\ref{sec:poissonacoustic}, we
used the optimization routines from
hIPPYlib~\cite{VillaPetraGhattas16}, a Python library for
deterministic and Bayesian inverse problems. A short description of
the numerical methods used for the solution of these problems can be
found in \ref{sec:numopt}. For details regarding the computation of
the adjoint-based derivatives we refer to~\cite{Crestel17}.

\subsection{Joint Poisson inverse problems with different observation points}
\label{sec:poissonpoisson}

Here, we solve a joint inverse problem of the form~\eqref{eq:joint2}
for the two coefficient fields $m_1$ and $m_2$. Considered separately,
$m_1$ and $m_2$ are solutions to the (almost identical) TV-regularized
inverse problems governed by the Poisson equation, i.e.,
\begin{equation} \label{eq:poisson} \begin{gathered}
m_i \coloneqq \arg\min_m \left\{ \frac12 | B_iu - \dd_i |^2  + \gamma_i
\int_\Omega \sqrt{ |\nabla m|^2 + \varepsilon} \, dx
\right\}, \quad \text{where}\\
\left\{ \begin{aligned}
- \nabla \cdotp ( e^m \nabla u) & = 1, \, \text{ in } \Omega, \\
u & = 0 , \, \text{ on } \partial \Omega.
\end{aligned} \right.
\end{gathered} \end{equation}
The operators~$B_i$ represent pointwise observation operators, and the
data~$\dd_i$ are synthetic observations polluted with 2\%~Gaussian noise.  
The domain $\Omega$ is discretized with a mesh of 8192
triangles (i.e., $h=1/64$).  
In all experiments presented in this section, the initial guesses for
both parameter fields are constant zero over the domain, i.e.,
$m_1^0 \equiv 0$ and $m_2^0 \equiv 0$.

The differences between the inverse problems for~$m_1$ and~$m_2$
reside in the truth parameter fields, and in the observation operators
$B_i$.  In the first example (section~\ref{sec:coincide}), the truth
parameter fields differ but have interfaces at the same spatial
locations.  In the second example (section~\ref{sec:coincide2}), some
interfaces in the truth parameter field for~$m_2$ are not present in the
truth parameter field for~$m_1$.
In both examples, the observation locations defined by $B_1$ only
cover the top-right quadrant of the domain, whereas the observation
locations defined by $B_2$ are distributed over the entire domain; see
figures~\ref{fig:coincide-target} and~\ref{fig:coincide2-target}.

\subsubsection{Truth parameter fields having identical interface locations}
\label{sec:coincide}

\newcommand{\wca}{0.23\textwidth}
\pgfplotstableread[col sep=comma]{fig/coincide/targetcs.csv}\targetcs
\begin{figure}%
\centering
\begin{tabular}{c@{\hspace{-.1in}}c@{}c@{}c@{}c}
\tikzsetnextfilename{coincide-m1}
\begin{tikzpicture}
\node[anchor=south west, inner sep=0] (image1) at (0,0)
{\includegraphics[width=\wca, trim=50 65 50 70, clip=true]
{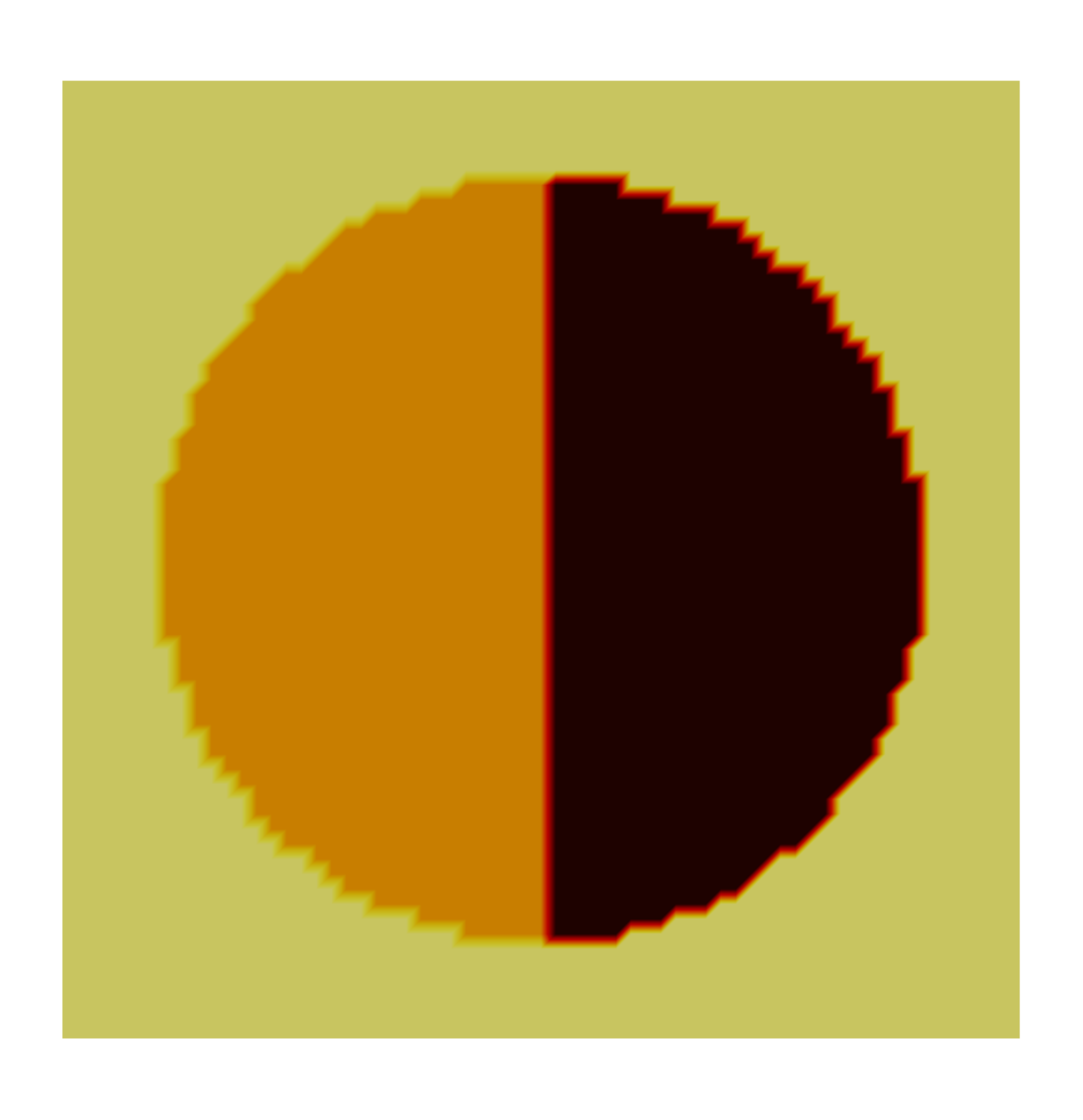}};
\begin{scope}[x={(image1.south east)},y={(image1.north west)}]
\foreach \x in {26,27,...,50} {
\foreach \y in {26,27,...,50} {
\filldraw[fill=black,draw=white] (0.01961*\x,0.01961*\y) circle (0.002);}}
\end{scope}
\end{tikzpicture}
&
\tikzsetnextfilename{coincide-m2}
\begin{tikzpicture}
\node[anchor=south west, inner sep=0] (image1) at (0,0)
{\includegraphics[width=\wca, trim=50 65 50 70, clip=true]
{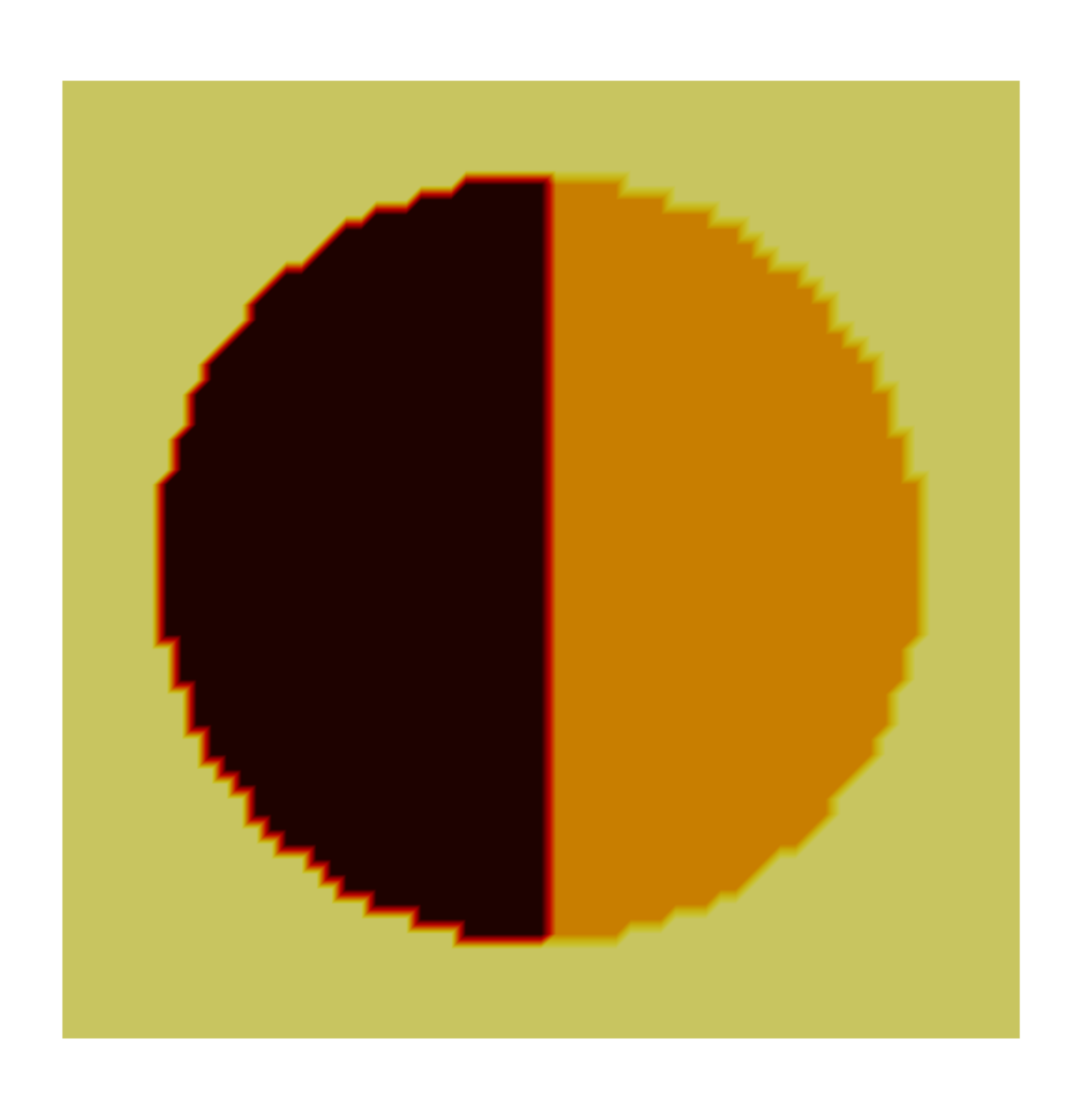}};
\begin{scope}[x={(image1.south east)},y={(image1.north west)}]
\foreach \x in {1,2,...,50} {
\foreach \y in {1,2,...,50} {
\filldraw[fill=black,draw=white] (0.01961*\x,0.01961*\y) circle (0.002);}}
\end{scope}
\end{tikzpicture}
&
\tikzsetnextfilename{coincide-legend}
\begin{tikzpicture}
\node[anchor=south west, inner sep=0] (image1) at (0,0)
{\includegraphics[height=\wca, trim=460 320 390 220, clip=true]
{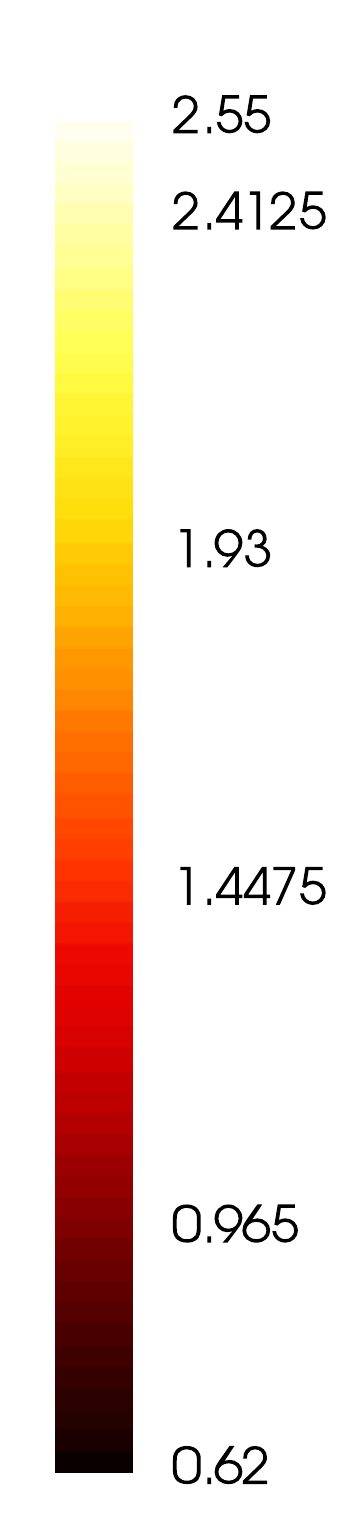}};
\begin{scope}[x={(image1.south east)},y={(image1.north west)}]
\end{scope}
\end{tikzpicture}
&
\includegraphics[width=\wca, trim=50 65 50 70, clip=true]
{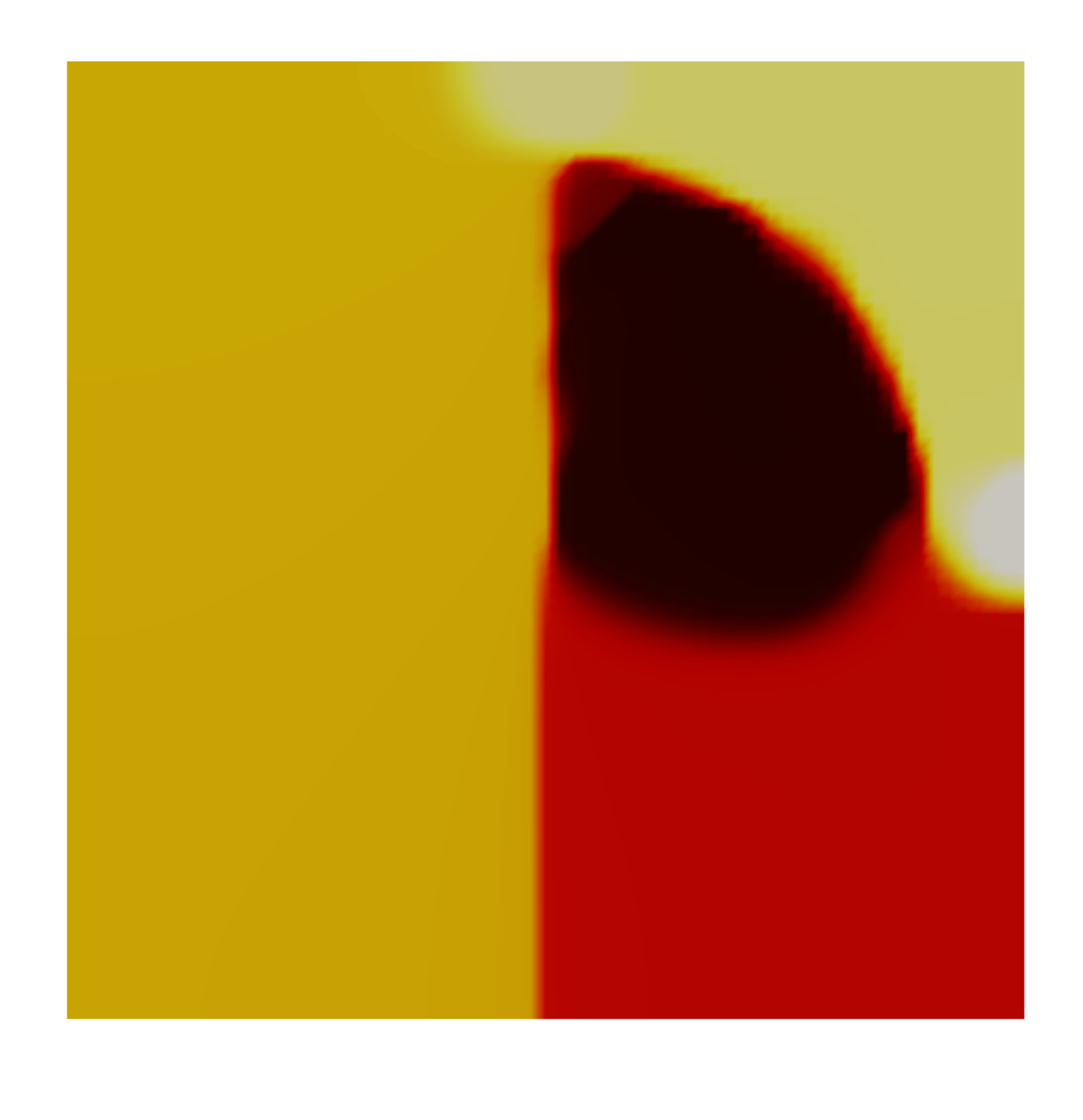}
&  
\includegraphics[width=\wca, trim=50 65 50 70, clip=true]
{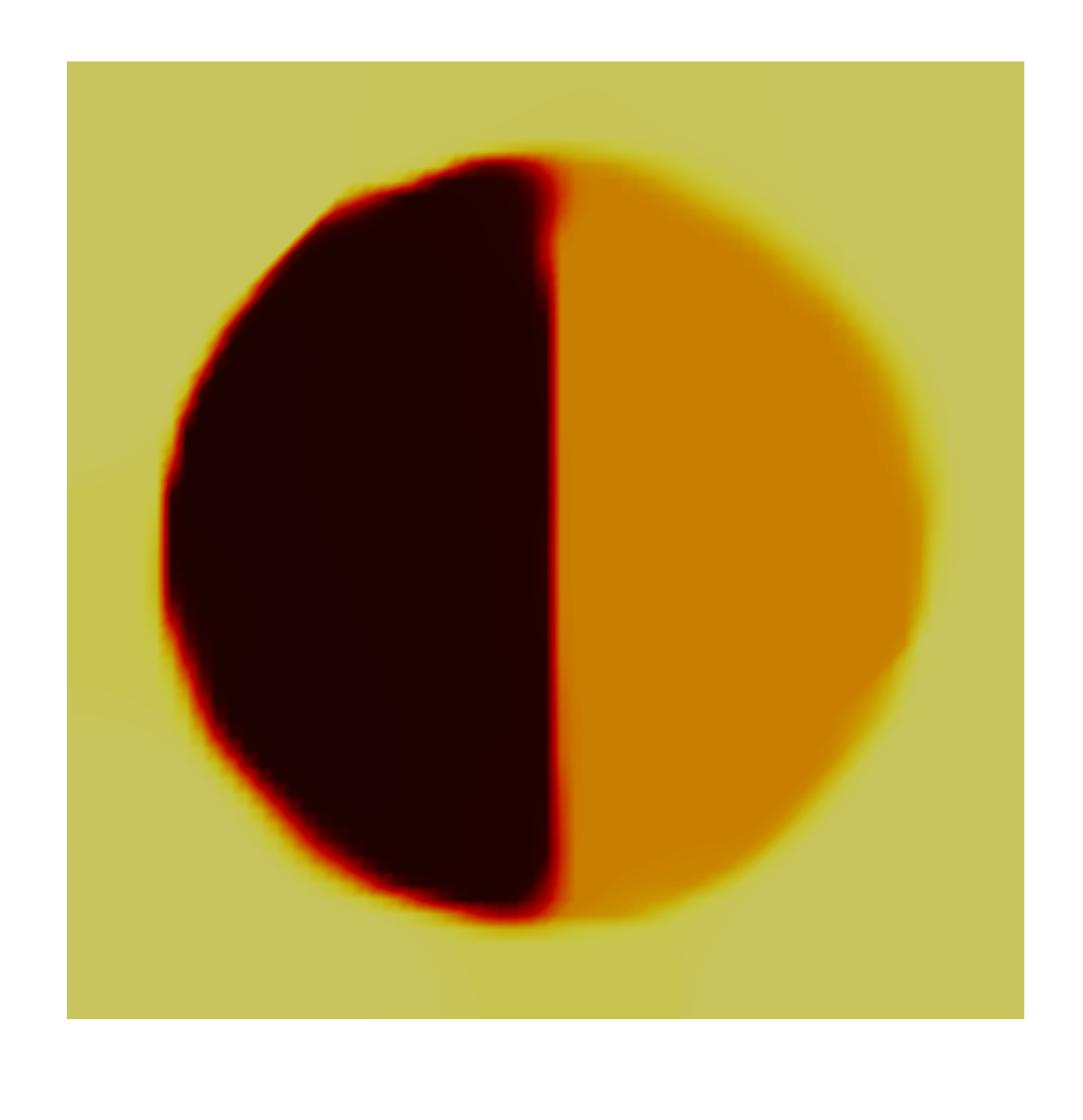} 
\\
(a) Truth $m_1$ &  (b) Truth $m_2$ & 
& (c) Indep $m_1$ & (d) Indep $m_2$
\end{tabular}
\caption{ Parameter fields $m_1$ and $m_2$ in the example of
section~\ref{sec:coincide}: truth parameter fields (a,b) and reconstructions
(c,d) obtained by solving the inverse problem~\eqref{eq:poisson} with
$\varepsilon = 10^{-3}$, $\gamma_1=3 \cdotp 10^{-7}$, $\gamma_2=4 \cdotp
10^{-7}$, and initial guesses $m_1^0 = m_2^0 = 0$.  
White dots in (a) and (b) indicate the location of the pointwise
observations.  The observation points defined through $B_1$ are a lattice of $25
\times 25$ points that cover only the top-right quadrant of the domain.  The
observation points for $B_2$ are a square lattice of~$50 \times 50$ points
distributed over the entire domain.}
\label{fig:coincide-target}
\end{figure}
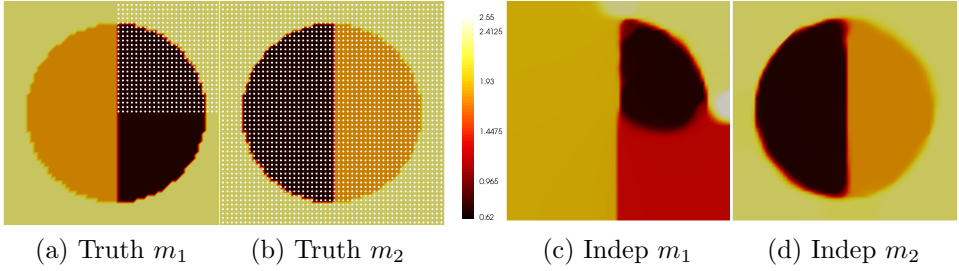
In the first example, the parameter fields
have interfaces at the same locations.
In figure~\ref{fig:coincide-target}, we show the truth parameter
fields $m_1$ and $m_2$ and their reconstructions
obtained by solving the inverse
problems~\eqref{eq:poisson} independently.
The reconstructions obtained with the four regularization methods are
shown in figure~\ref{fig:coincide-reconstruct}, and the corresponding values of
the relative medium misfit are given in table~\ref{tab:all-med}.
\begin{figure}%
\centering
\begin{tabular}{l@{\hspace{.08in}}c@{\hspace{.01in}}c@{\hspace{.01in}}c@{\hspace{.01in}}c}
\rotatebox[origin=l]{90}{\hspace{0.4in}(a) $m_1$}
&
\includegraphics[height=\wca, trim=50 65 50 70, clip=true]
{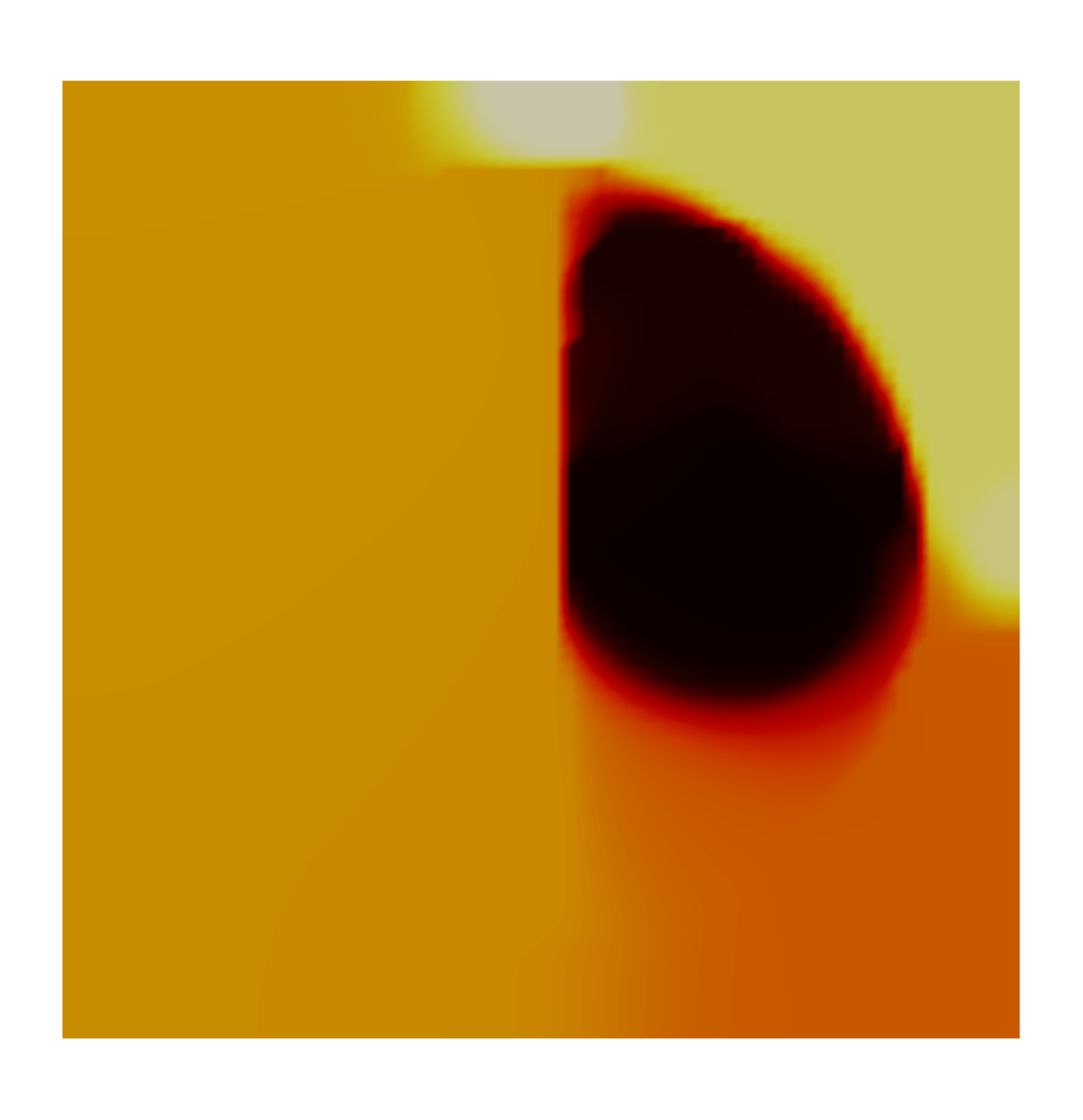}
& 
\includegraphics[height=\wca, trim=330 60 360 90, clip=true]
{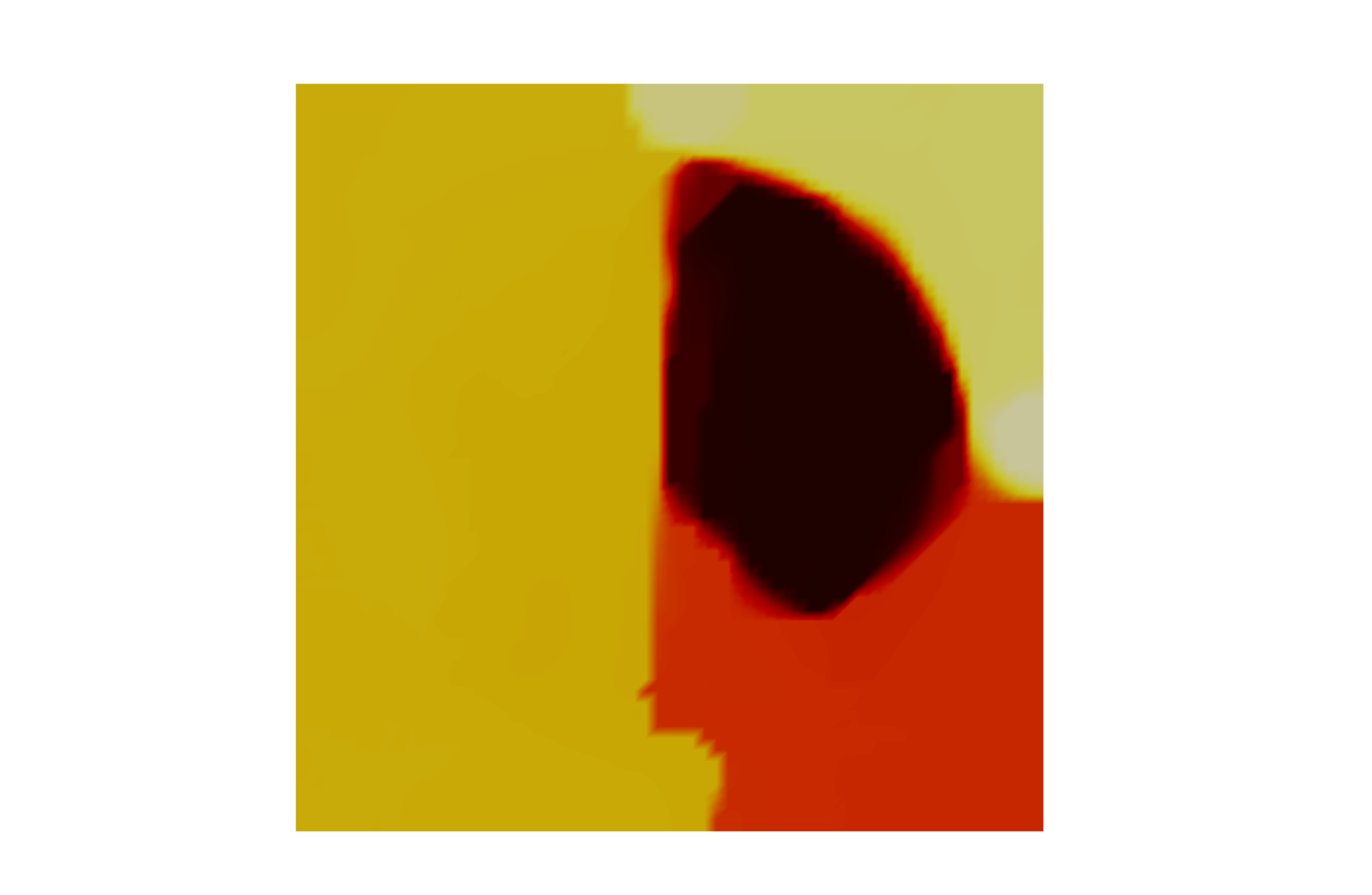}
&
\includegraphics[height=\wca, trim=50 65 50 70, clip=true]
{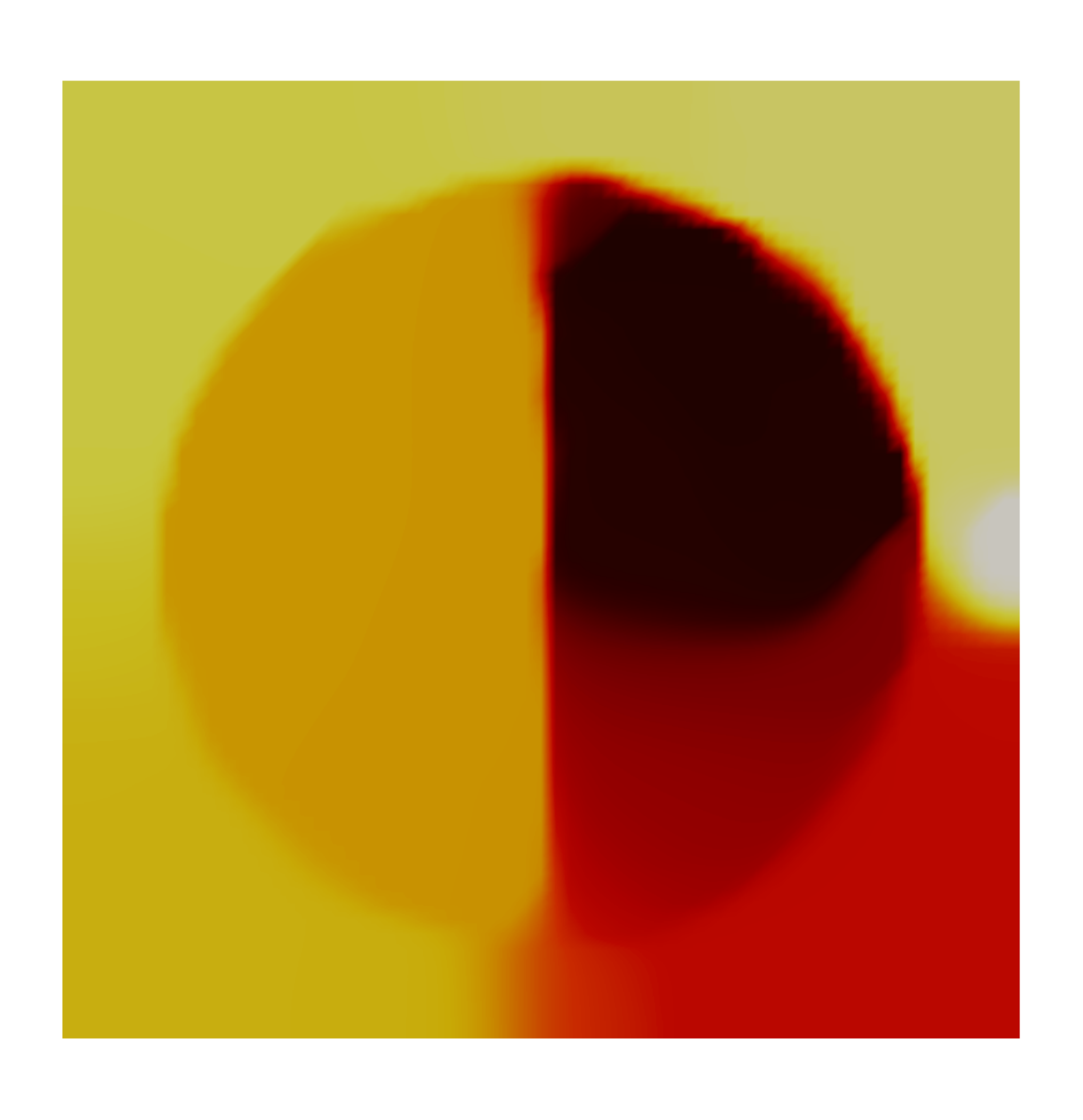}
& 
\includegraphics[height=\wca, trim=50 65 50 70, clip=true]
{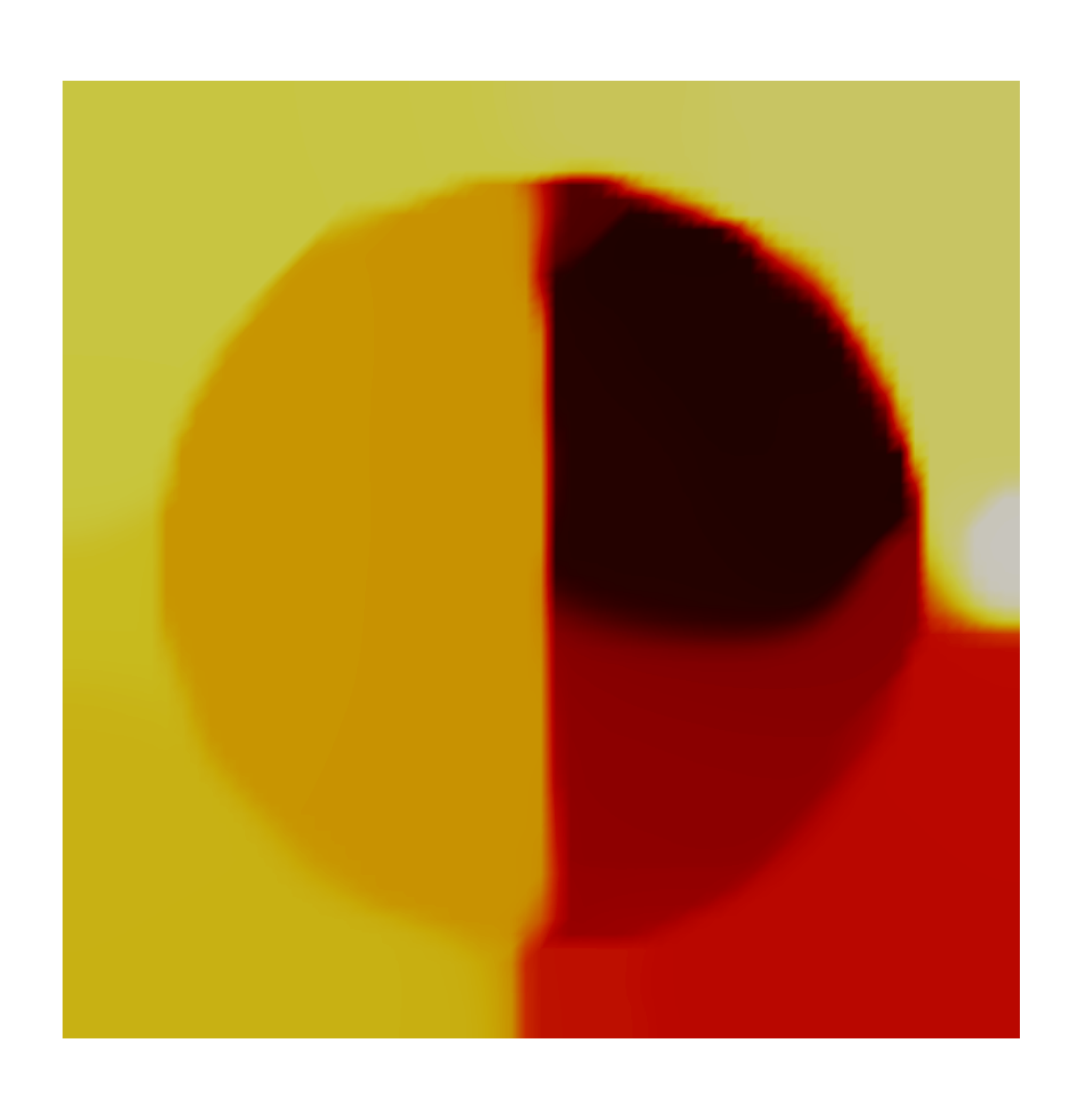}
\\
\rotatebox[origin=l]{90}{\hspace{0.4in}(b) $m_2$}
&
\includegraphics[height=\wca, trim=50 65 50 70, clip=true]
{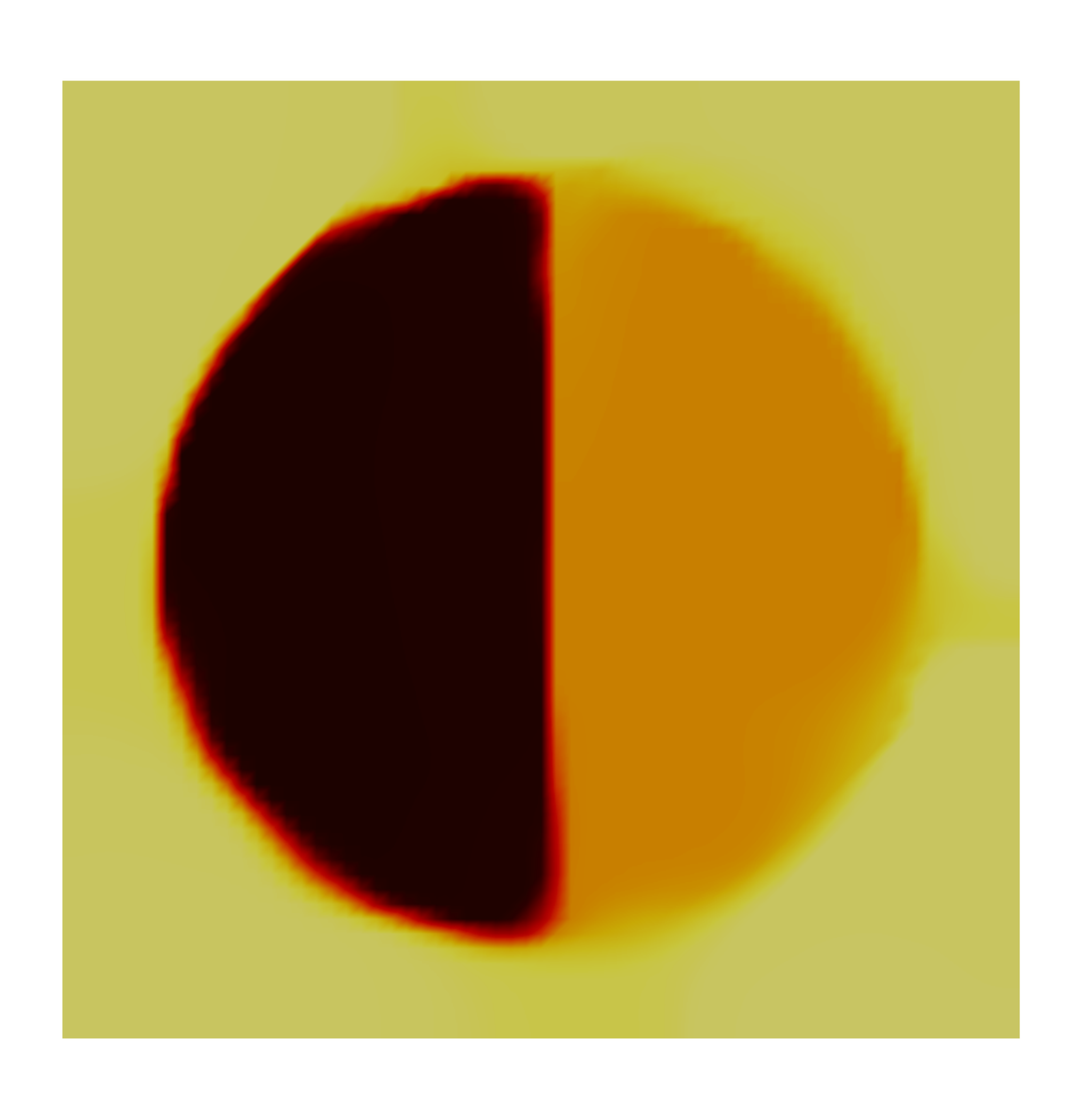} 
& 
\includegraphics[height=\wca, trim=330 60 360 90, clip=true]
{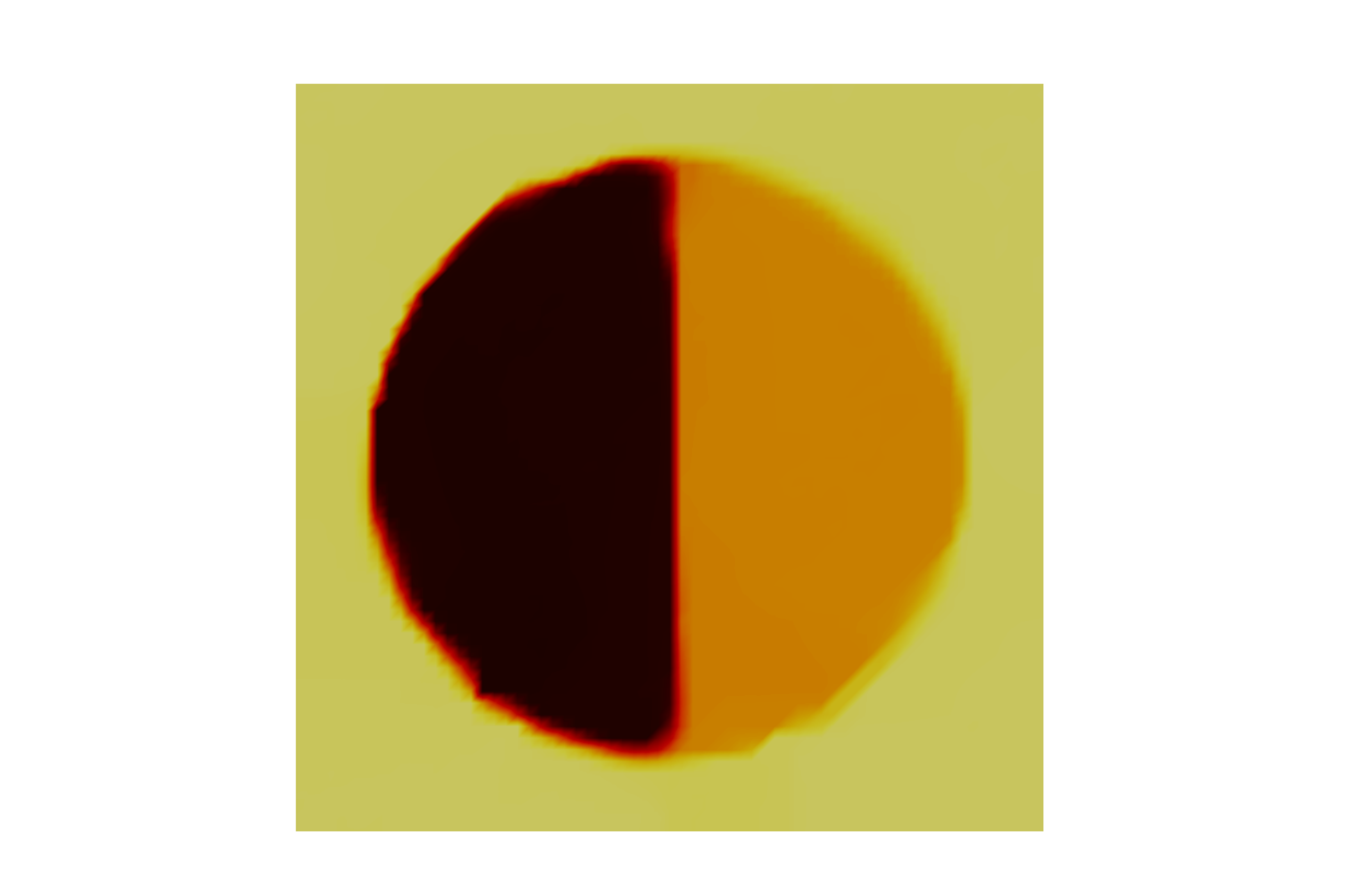} 
&
\includegraphics[height=\wca, trim=50 65 50 70, clip=true]
{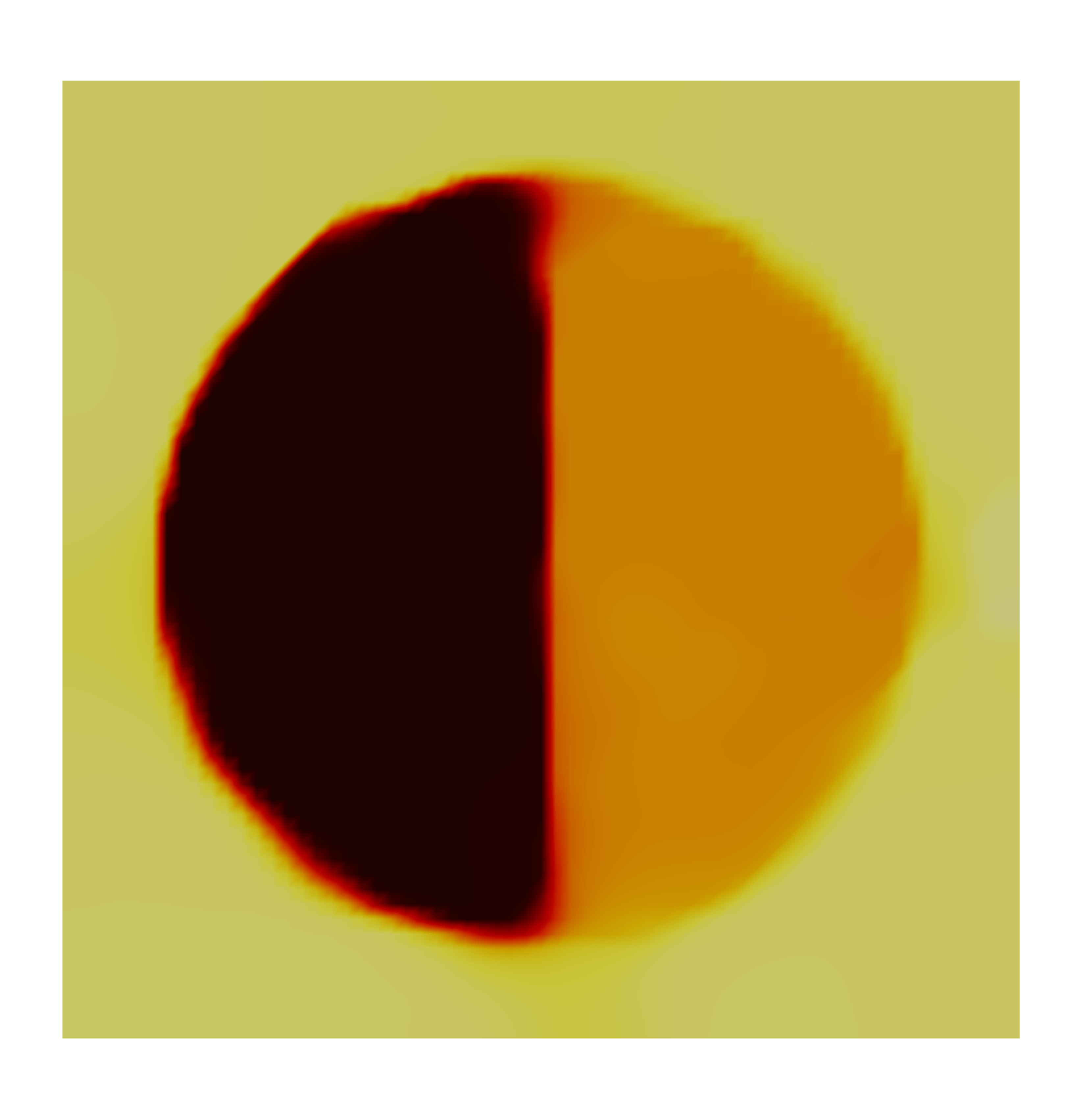}
& 
\includegraphics[height=\wca, trim=50 65 50 70, clip=true]
{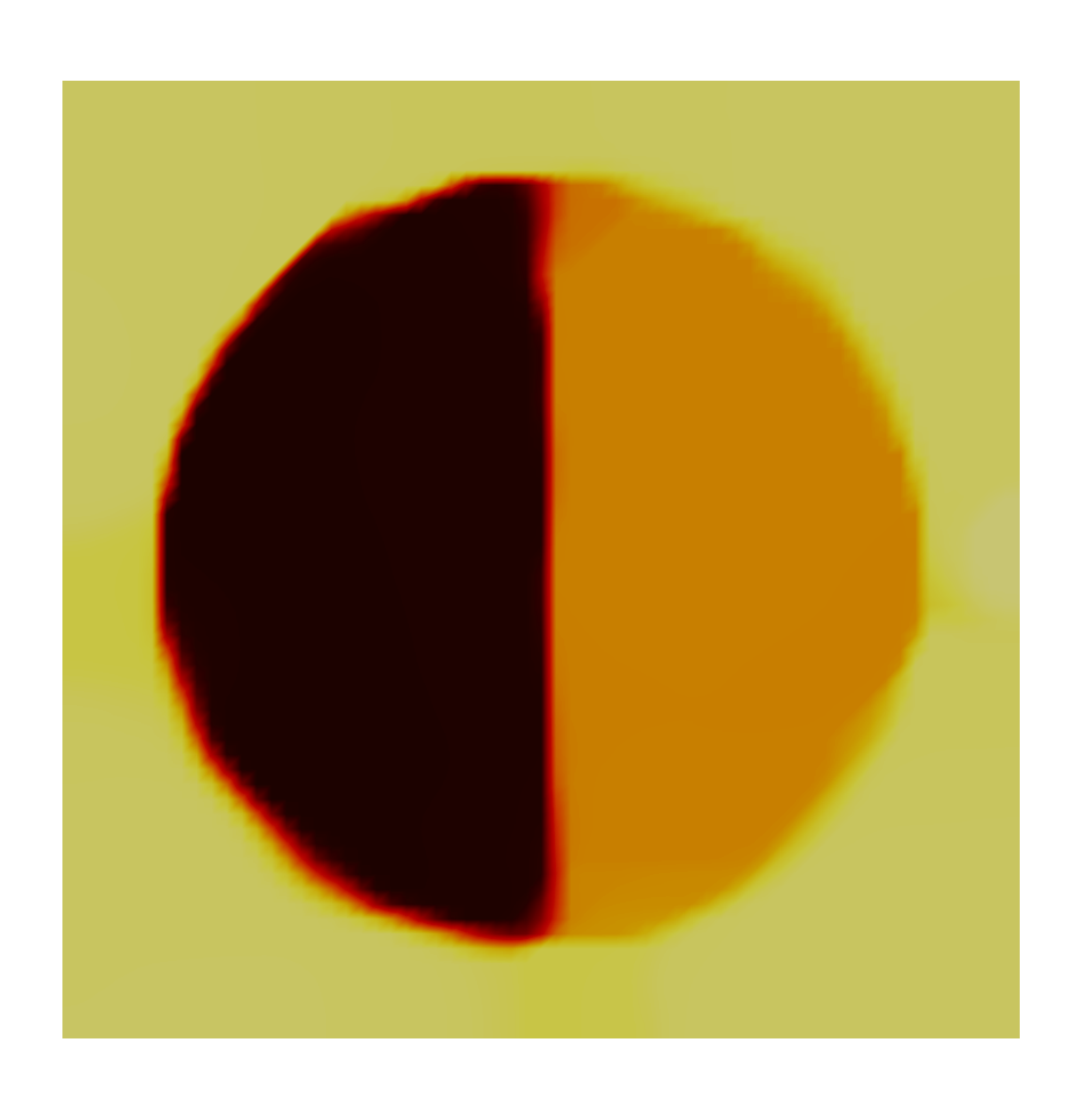}
\\
& (i) cross-gradient & (ii) norm. cross-gd &
(iii) vectorial TV & (iv) nuclear norm 
\end{tabular}
\caption{Reconstructions for the parameter fields (a) $m_1$ and (b) $m_2$,
obtained by solving a joint inverse problem~\eqref{eq:joint2} regularized with
(i) the cross-gradient ($\gamma=2 \cdotp 10^{-8}$) combined with two  independent
TV regularizations, (ii) the normalized cross-gradient ($\gamma=6 \cdotp
10^{-6}$ and $\varepsilon=10^{-3}$) combined with the same independent TV
regularizations, (iii) the VTV joint regularization ($\gamma=3 \cdotp 10^{-7}$
and $\varepsilon=10^{-3}$), and (iv) the nuclear norm joint regularization
($\gamma=3 \cdotp 10^{-7}$ and $\varepsilon=10^{-3}$).  The parameters for the
independent TV regularizations 
and the initial guesses for all problems
are as for the independent inverse problems (see
caption of figure~\ref{fig:coincide-target}).  The legend for all plots is 
as in figure~\ref{fig:coincide-target}.}
\label{fig:coincide-reconstruct}
\end{figure}

The reconstructions for parameter~$m_2$ do not differ significantly
(figure~\ref{fig:coincide-reconstruct}b).
Due to the large number of observation points, this parameter is already 
well reconstructed in an independent inverse problem
(figure~\ref{fig:coincide-target}d).  We observe an improvement in the
reconstruction of parameter~$m_1$ for all four joint inverse problems
compared to the
independent reconstruction shown in figure~\ref{fig:coincide-target}c.
Using the cross-gradient only marginally improves the
reconstruction for parameter~$m_1$, most likely because the independent
reconstruction for~$m_1$ shows large areas of constant values, where the
cross-gradient term vanishes; these areas 
therefore cannot be improved by the cross-gradient.  The normalized cross-gradient
improves over the cross-gradient but fails to recover the
circular interface.  Both the VTV joint regularization and the
nuclear norm joint regularization perform better in this example, and lead to
reconstructions that contain all sharp interfaces in the target image.

\subsubsection{Truth parameter fields having different interface locations}
\label{sec:coincide2}

\begin{figure}
\centering
\begin{tabular}{c@{\hspace{-.02in}}c@{\hspace{-.05in}}c@{\hspace{-.05in}}c@{\hspace{.02in}}c}
\tikzsetnextfilename{coincide2-m1}
\begin{tikzpicture}
\node[anchor=south west, inner sep=0] (image1) at (0,0)
{\includegraphics[width=\wca, trim=180 70 180 70, clip=true]
{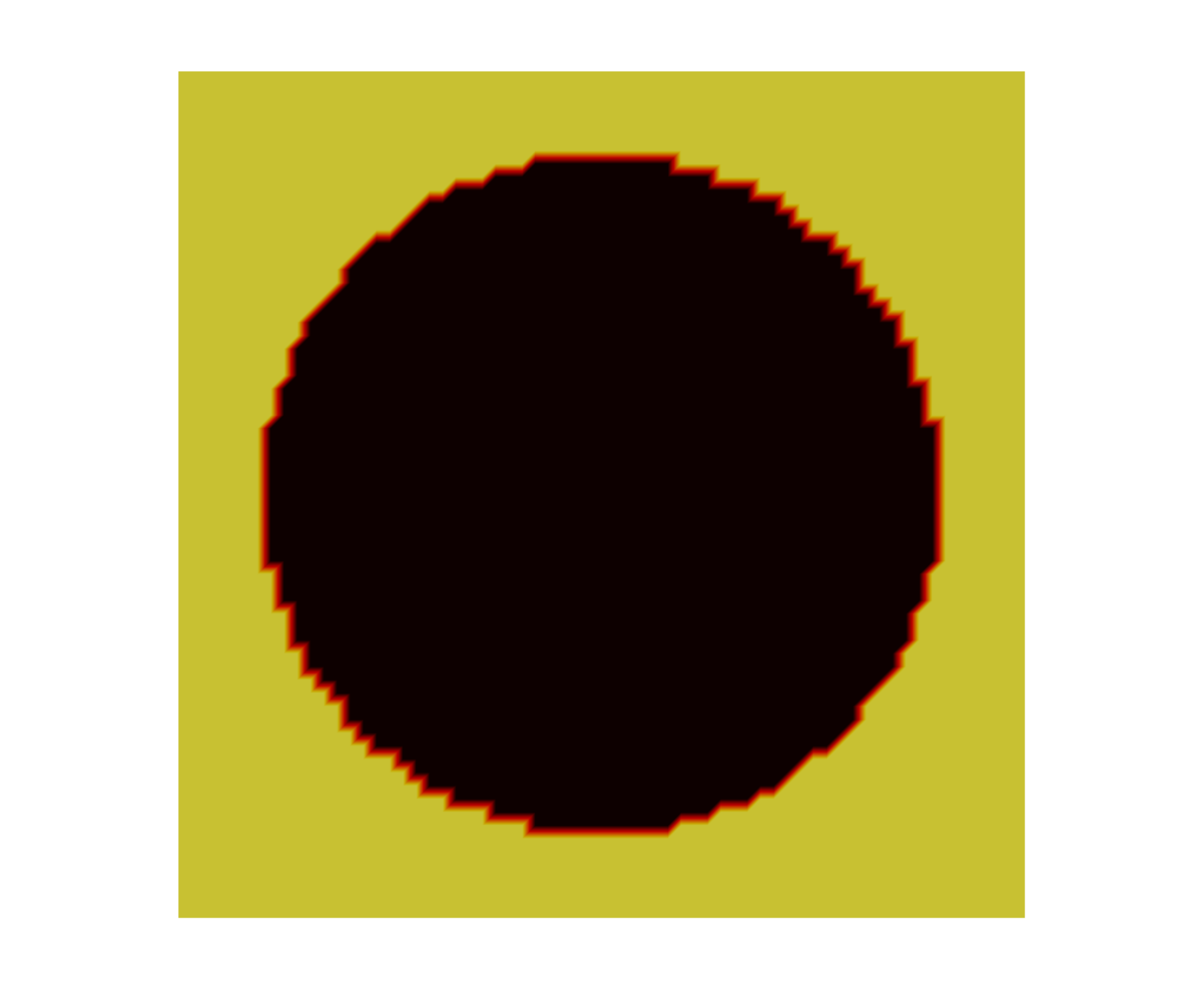}};
\begin{scope}[x={(image1.south east)},y={(image1.north west)}]
\foreach \x in {26,27,...,50} {
\foreach \y in {26,27,...,50} {
\filldraw[fill=black,draw=white] (0.01961*\x,0.01961*\y) circle (0.002);}}
\end{scope}
\end{tikzpicture}
&
\tikzsetnextfilename{coincide2-m2}
\begin{tikzpicture}
\node[anchor=south west, inner sep=0] (image1) at (0,0)
{\includegraphics[width=\wca, trim=180 70 180 70, clip=true]
{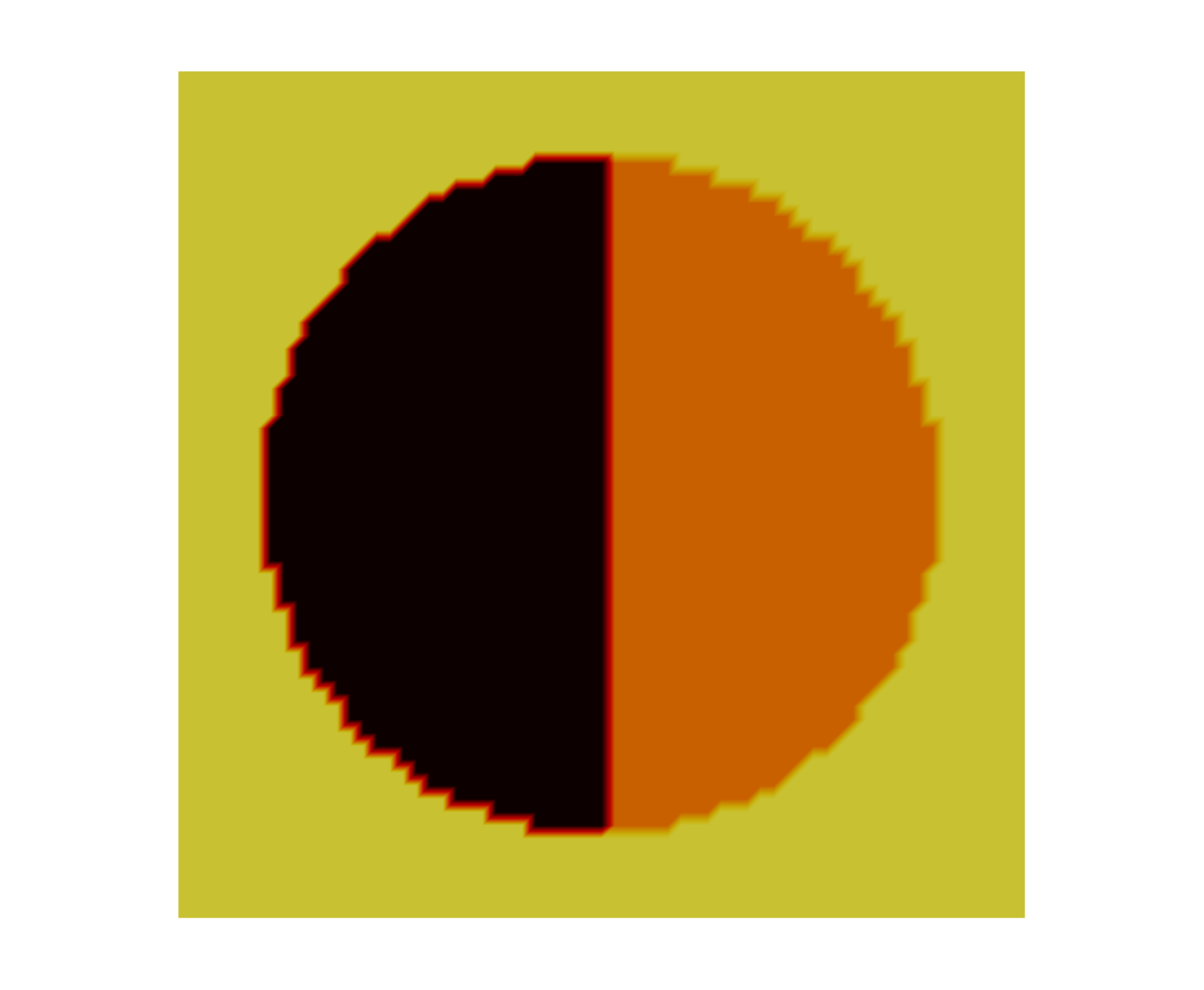}};
\begin{scope}[x={(image1.south east)},y={(image1.north west)}]
\foreach \x in {1,2,...,50} {
\foreach \y in {1,2,...,50} {
\filldraw[fill=black,draw=white] (0.01961*\x,0.01961*\y) circle (0.002);}}
\end{scope}
\end{tikzpicture}
&
\tikzsetnextfilename{coincide2-legend}
\begin{tikzpicture}
\node[anchor=south west, inner sep=0] (image1) at (0,0)
{\includegraphics[height=\wca, trim=515 270 550 280, clip=true]
{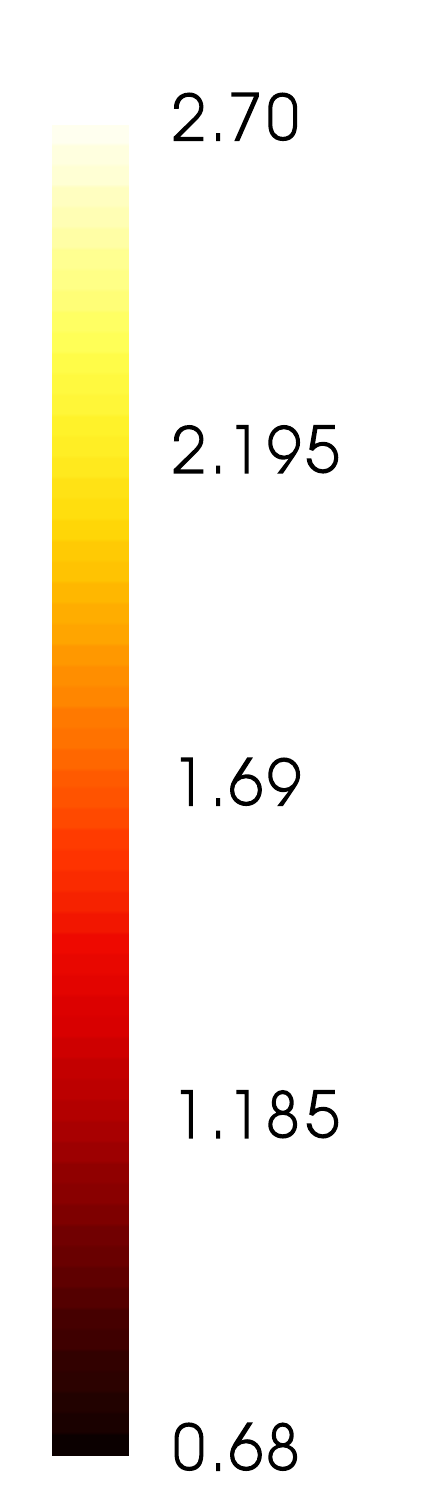}};
\begin{scope}[x={(image1.south east)},y={(image1.north west)}]
\end{scope}
\end{tikzpicture}
&
\includegraphics[width=\wca, trim=180 70 180 70, clip=true]
{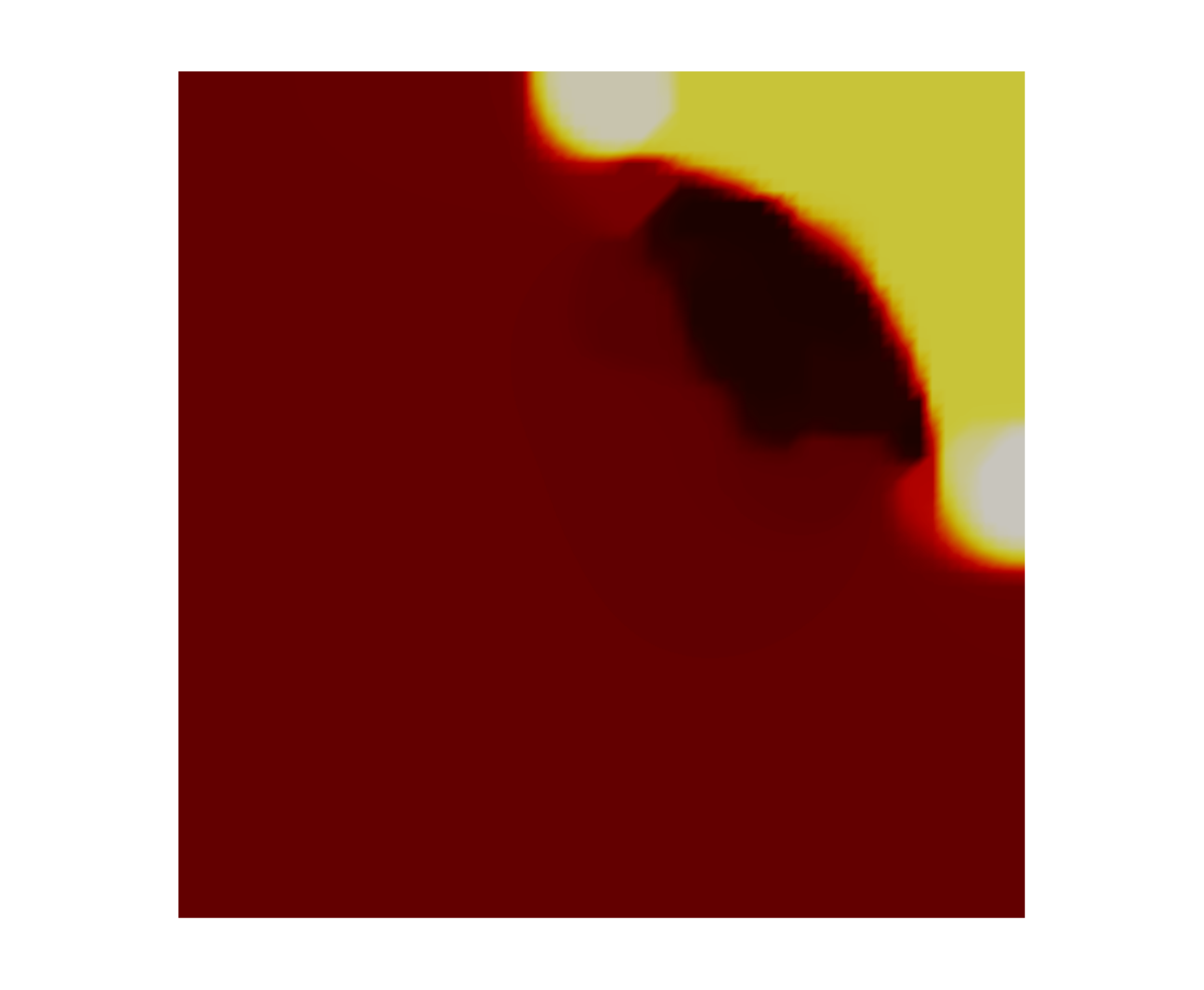}
&  
\includegraphics[width=\wca, trim=180 70 180 70, clip=true]
{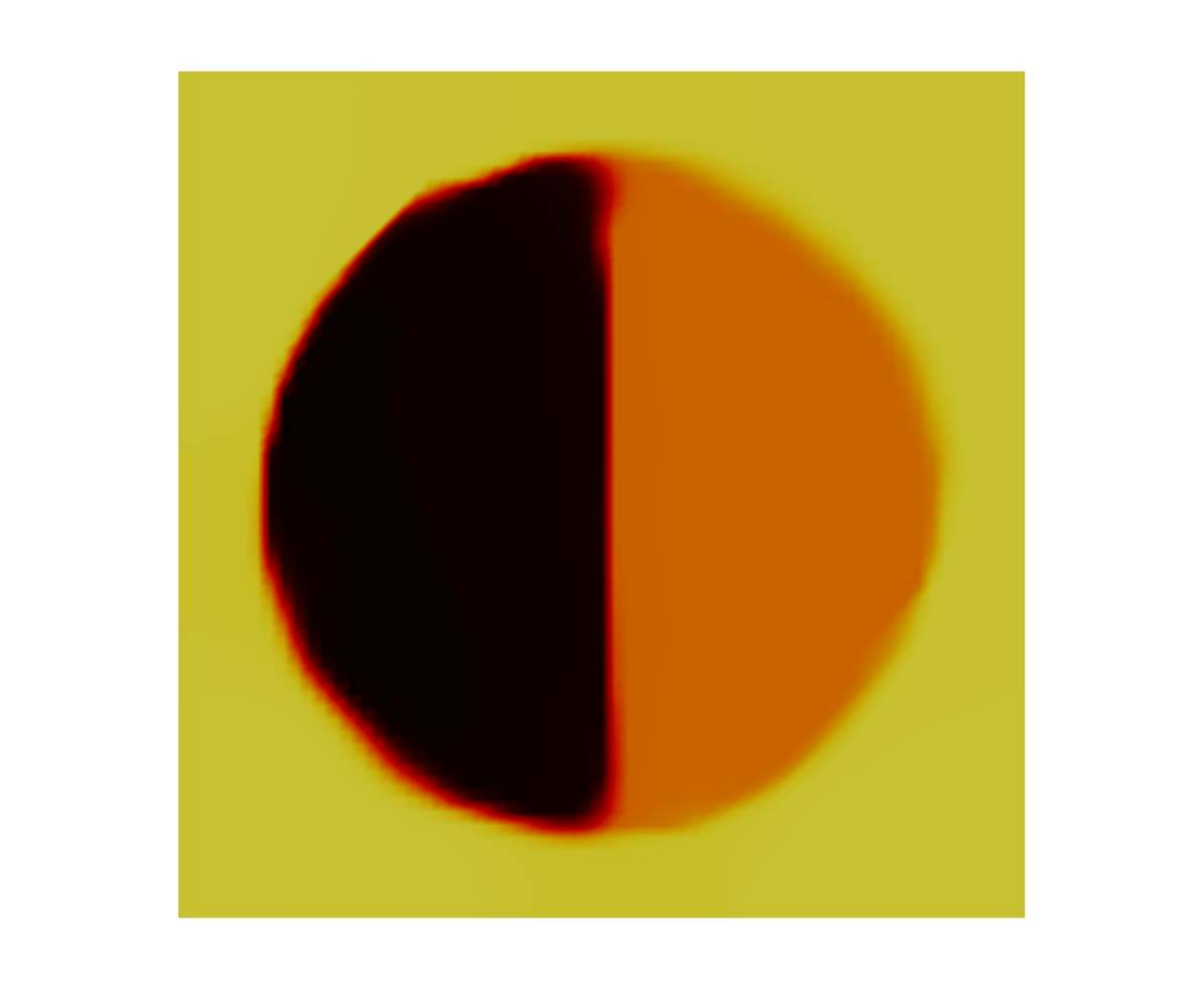} 
\\
(a) Truth $m_1$ & (b) Truth $m_2$ &
& (c) Indep $m_1$ & (d) Indep $m_2$
\end{tabular}
\caption{
Parameter fields for $m_1$ and $m_2$ in the example of
section~\ref{sec:coincide2}:
truth parameter field (a,b) and 
reconstructions (c,d)
obtained by solving the inverse problem~\eqref{eq:poisson} with $\varepsilon =
10^{-3}$, $\gamma_1=4 \cdotp 10^{-7}$, $\gamma_2=4 \cdotp 10^{-7}$, and
initial guesses $m_1^0 = m_2^0 = 0$.  
 White dots in (a) and (b) indicate the location of
the pointwise observations, as detailed in figure~\ref{fig:coincide-target}.
}
\label{fig:coincide2-target}
\end{figure}
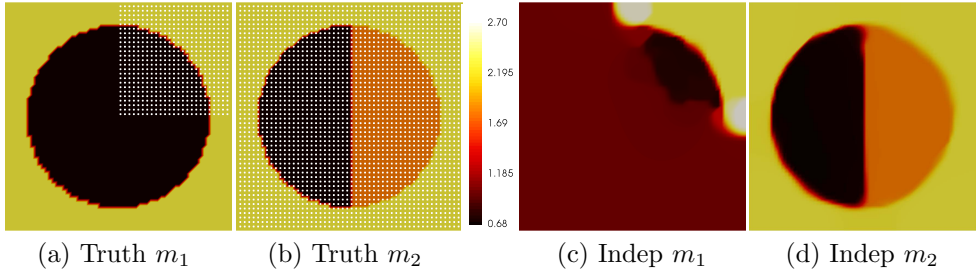

\begin{figure}%
\centering
\begin{tabular}{l@{\hspace{.08in}}c@{\hspace{.01in}}c@{\hspace{.01in}}c@{\hspace{.01in}}c}
\rotatebox[origin=l]{90}{\hspace{0.4in}(a) $m_1$}
&
\includegraphics[height=\wca, trim=180 70 180 70, clip=true]
{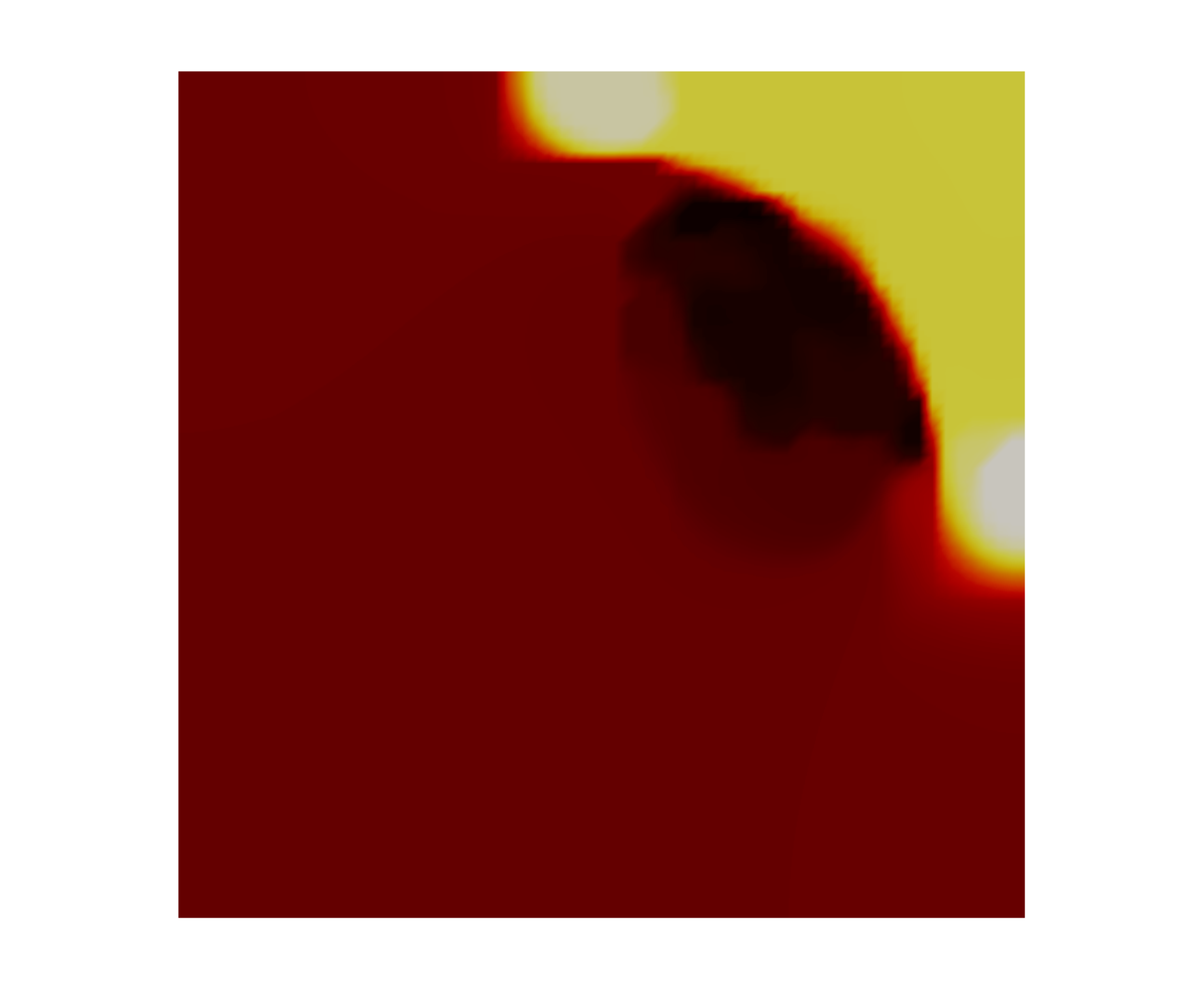}
& 
\includegraphics[height=\wca, trim=290 70 290 70, clip=true]
{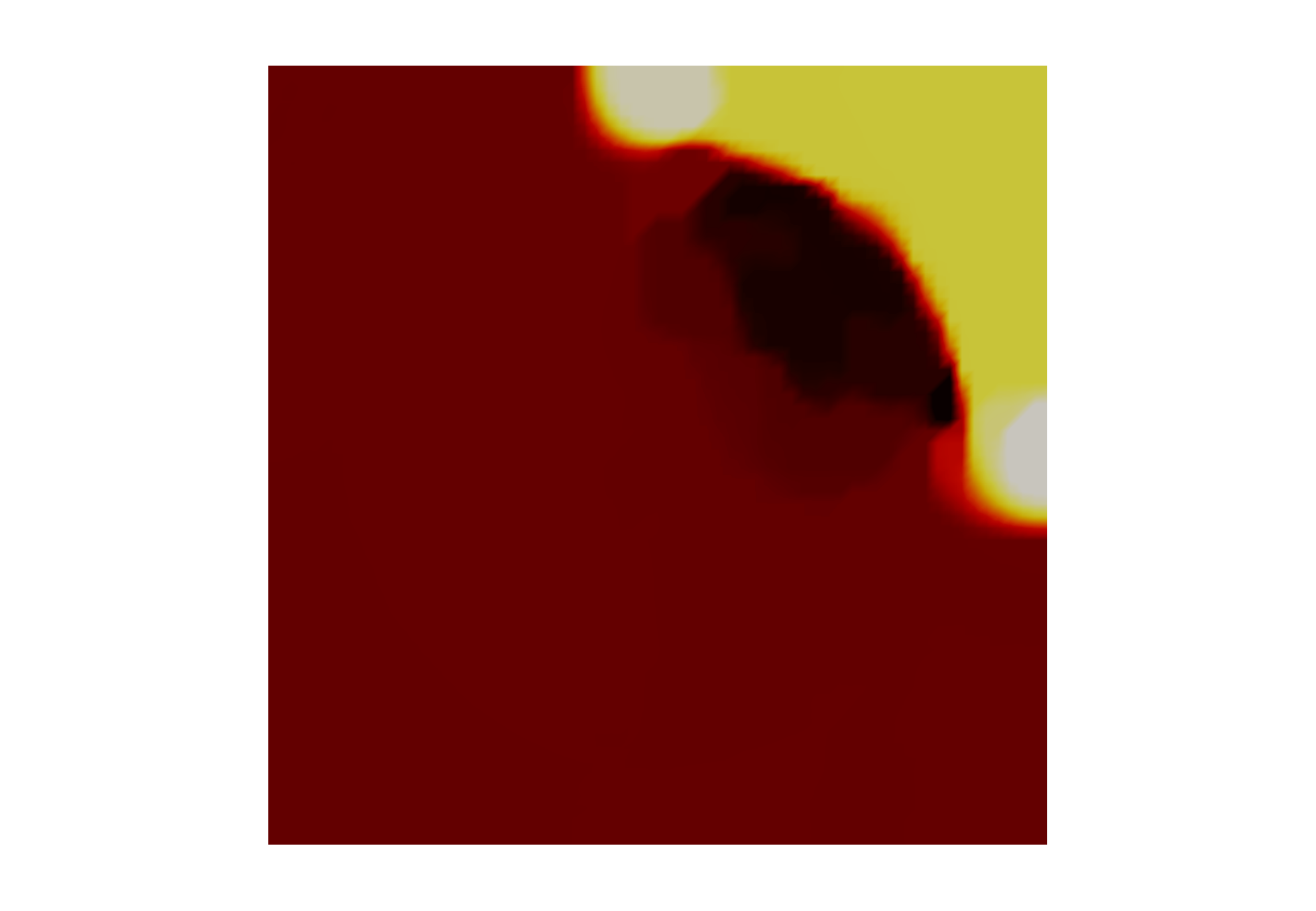}
& 
\includegraphics[height=\wca, trim=180 70 180 70, clip=true]
{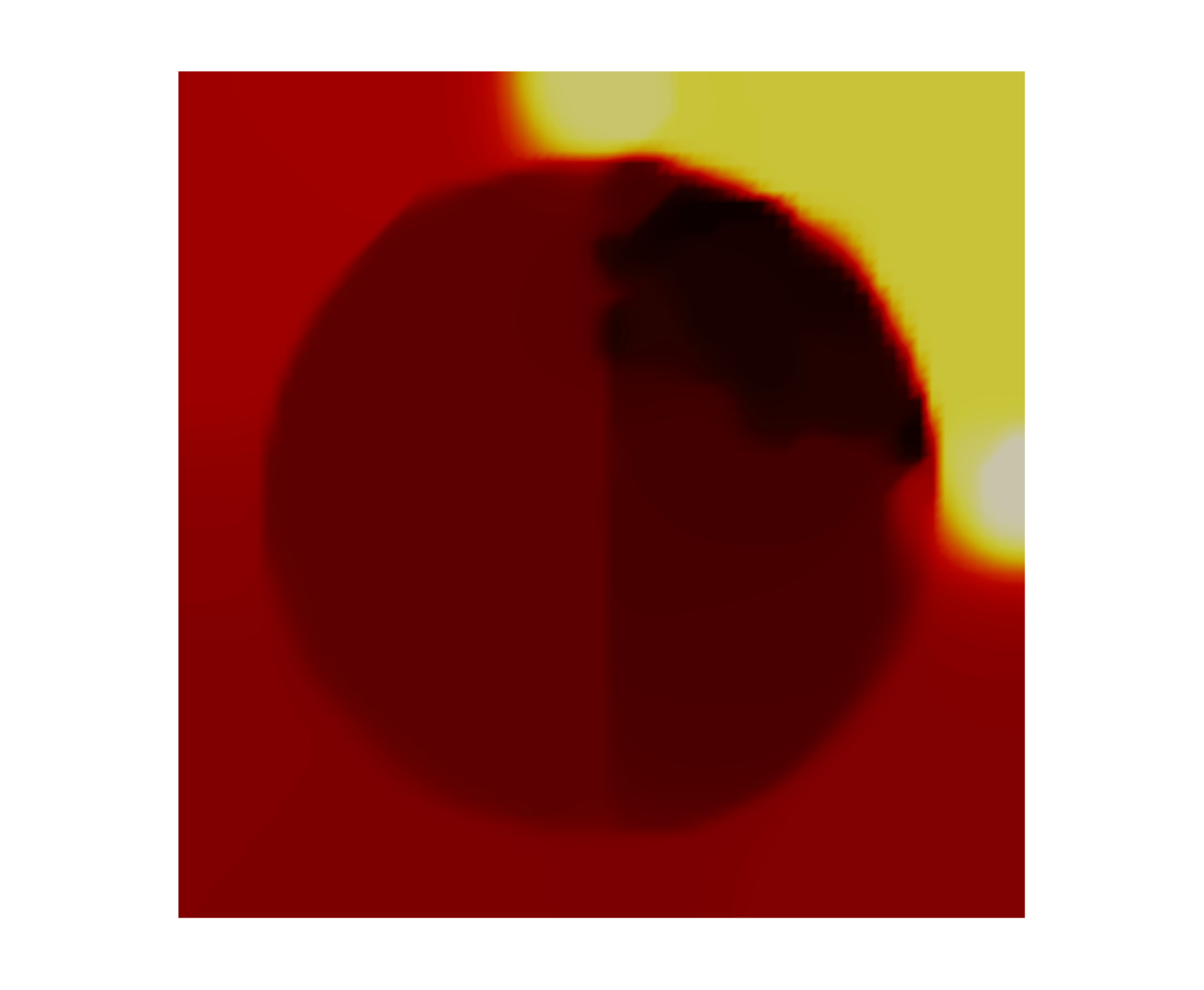}
& 
\includegraphics[height=\wca, trim=180 70 180 70, clip=true]
{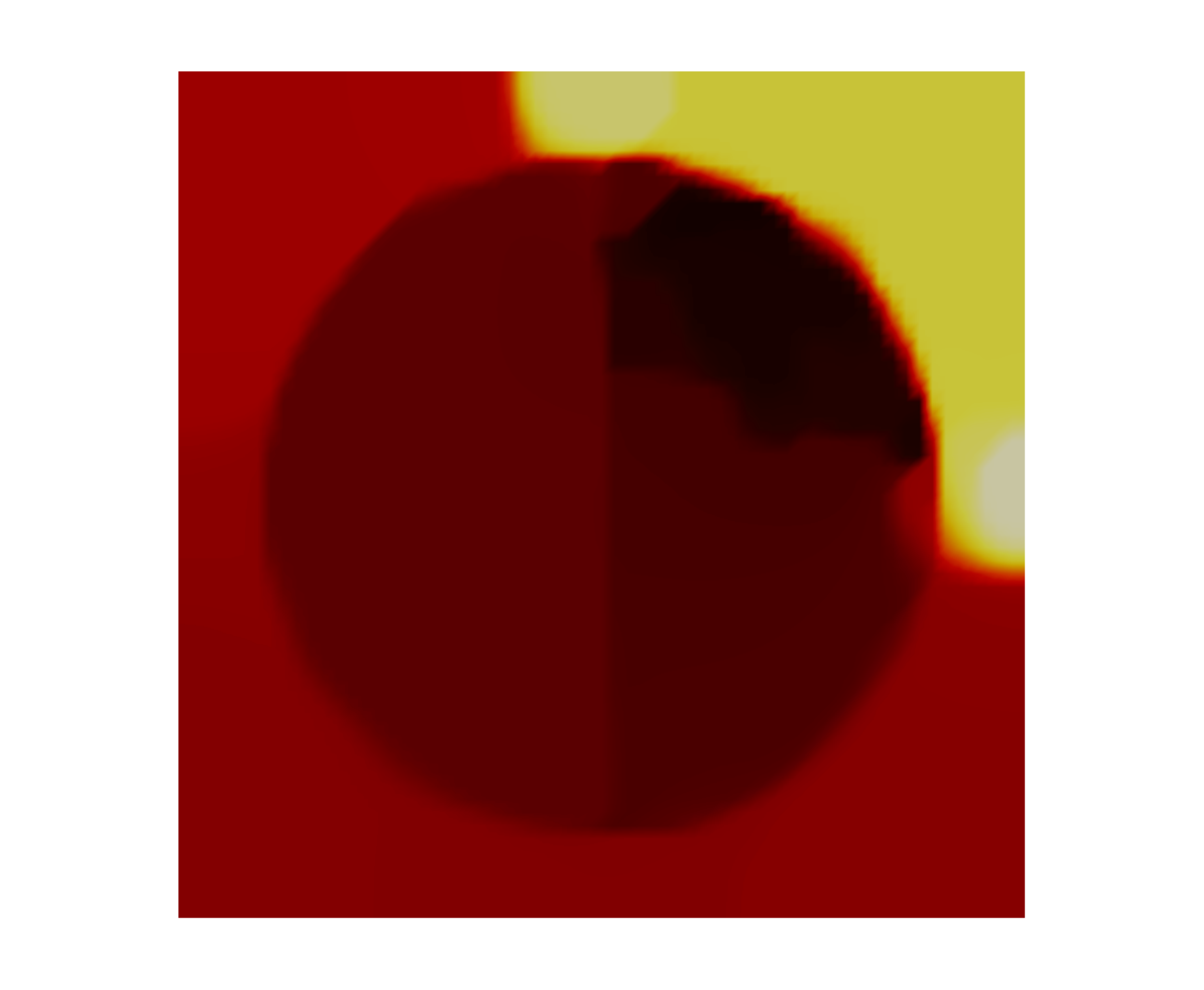}
\\
\rotatebox[origin=l]{90}{\hspace{0.4in}(b) $m_2$}
&
\includegraphics[height=\wca, trim=180 70 180 70, clip=true]
{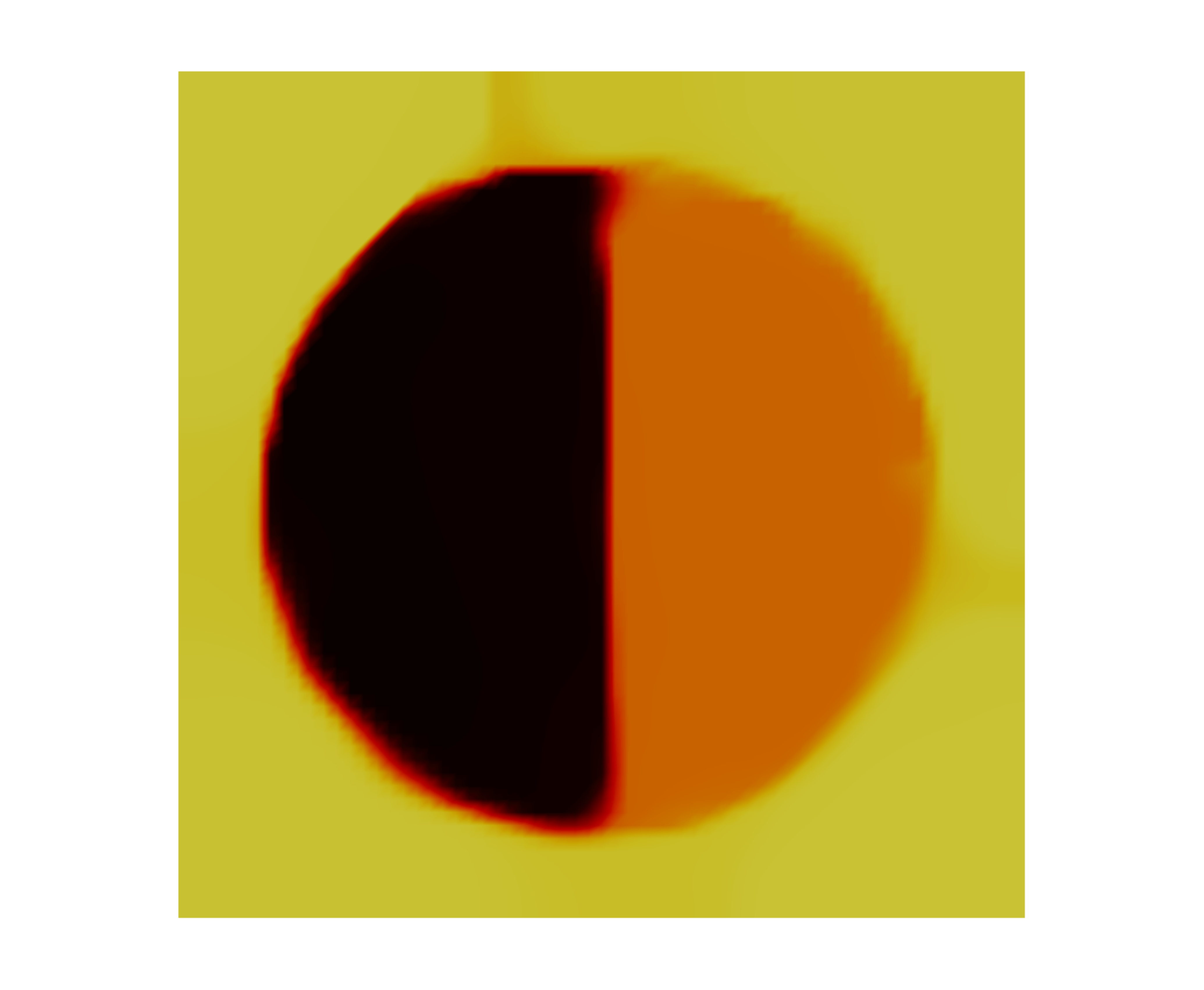} 
& 
\includegraphics[height=\wca, trim=290 70 290 70, clip=true]
{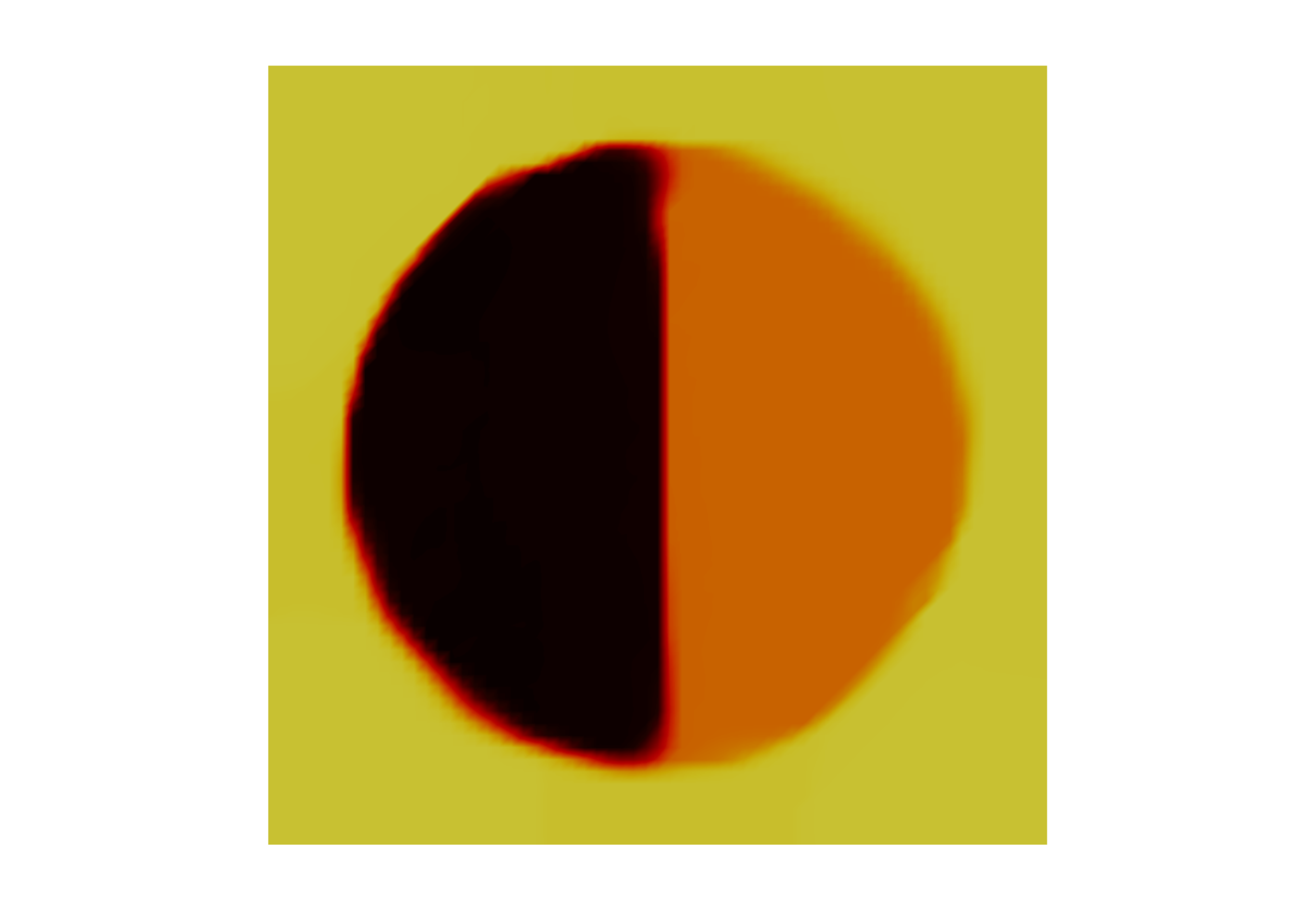}
&
\includegraphics[height=\wca, trim=180 70 180 70, clip=true]
{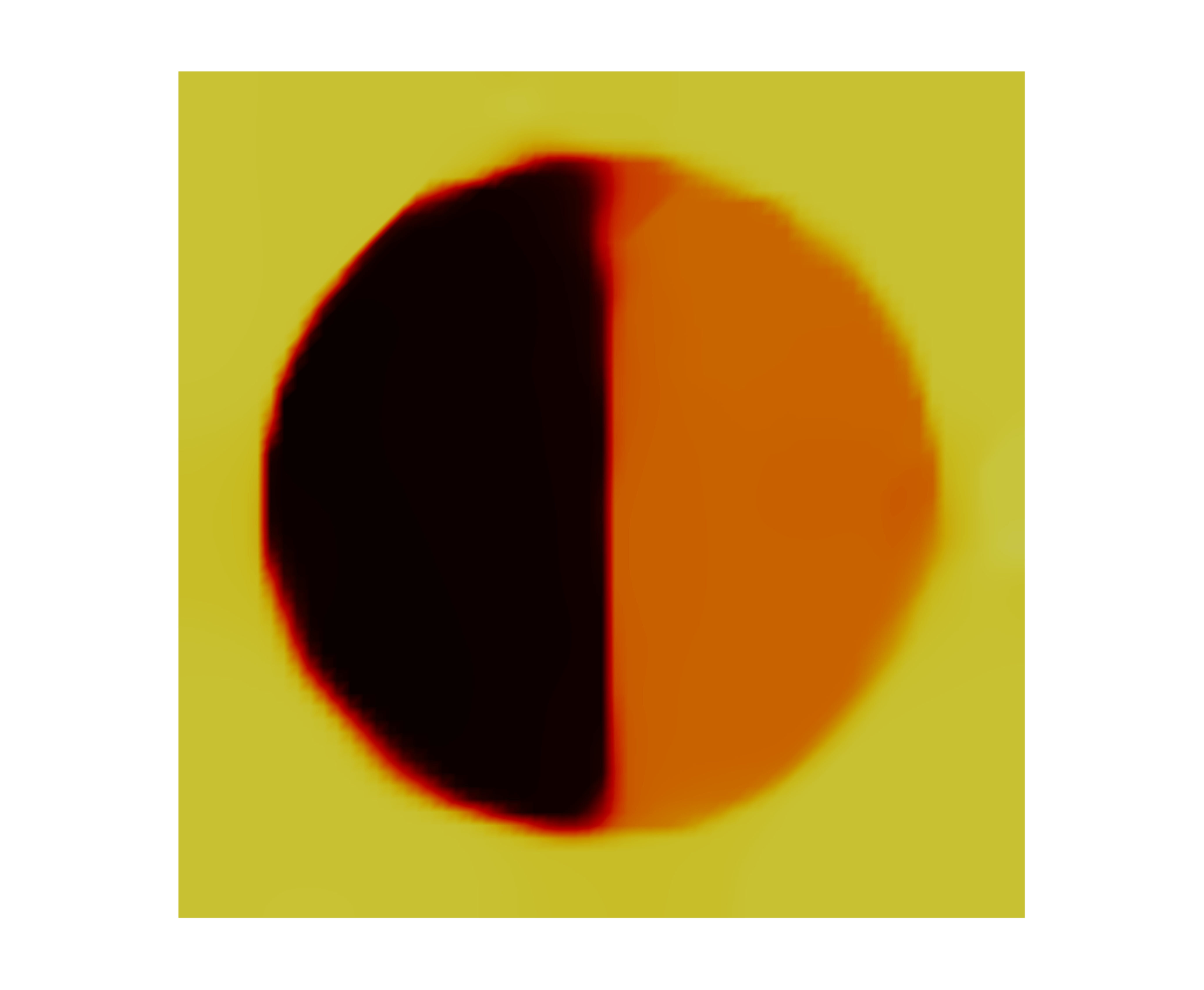}
& 
\includegraphics[height=\wca, trim=180 70 180 70, clip=true]
{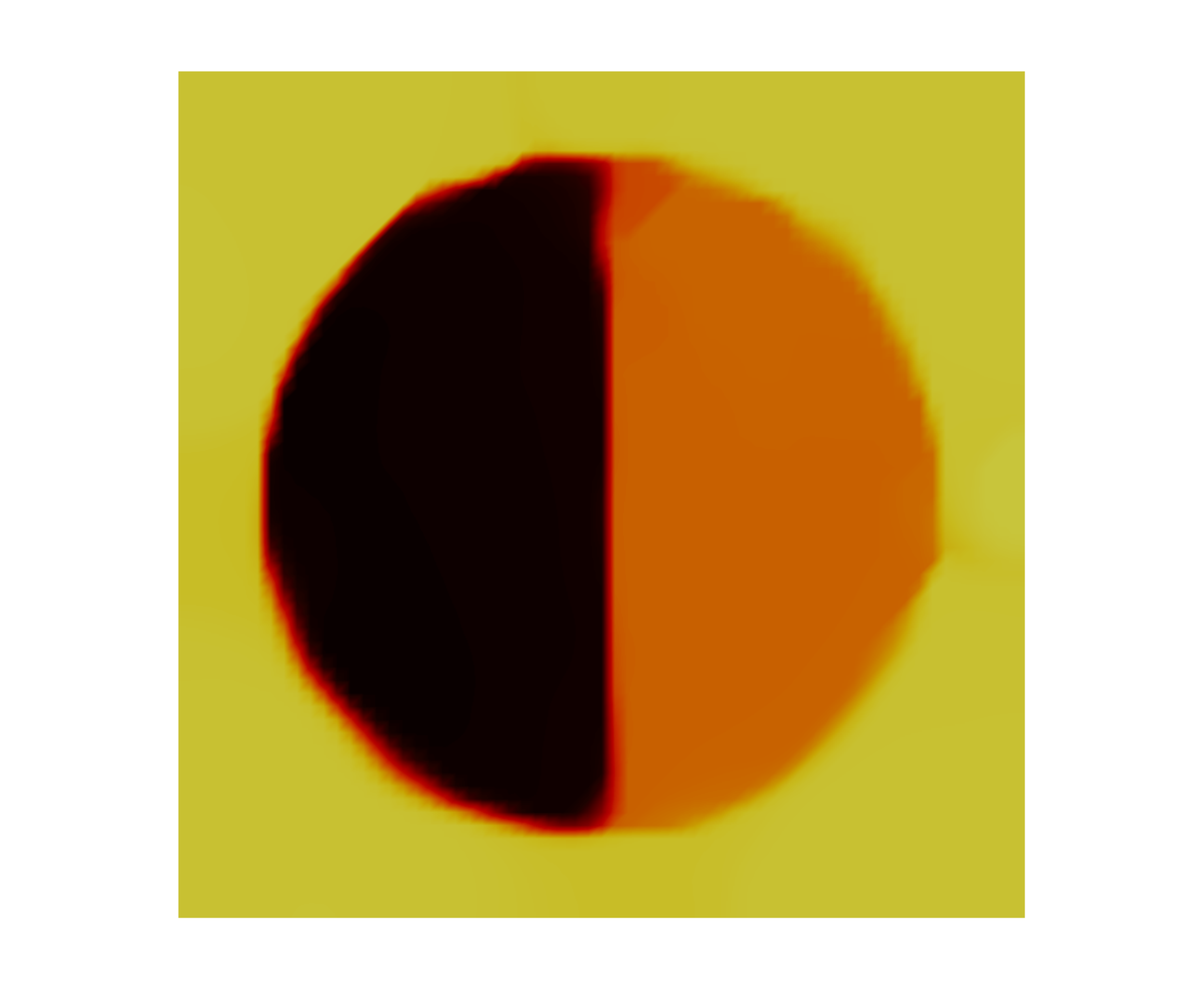}
\\
& (i) cross-gradient & (ii) norm. cross-gd &
(iii) vectorial TV & (iv) nuclear norm 
\end{tabular}
\caption{Reconstructions for the parameter fields (a) $m_1$ and (b) $m_2$, obtained
by solving a joint inverse problem~\eqref{eq:joint2} regularized with 
(i) the cross-gradient combined with 2 independent TV regularizations ($\gamma=5 \cdotp 10^{-9}$),
(ii) the normalized cross-gradient combined with the same independent TV regularizations
($ \gamma=7 \cdotp 10^{-7}$ and $\varepsilon=10^{-3}$),
(iii) the VTV joint regularization ($\gamma=4 \cdotp 10^{-7}$ and $\varepsilon=10^{-3}$),
and (iv) the nuclear norm joint regularization ($\gamma=4 \cdotp 10^{-7}$ and $\varepsilon=10^{-3}$).
The parameters for the independent TV regularizations 
and all initial guesses %
are the same as used for
the independent inverse problems (see caption in figure~\ref{fig:coincide2-target}).
The legend  is as in figure~\ref{fig:coincide2-target}.}
\label{fig:coincide2-reconstruct}
\end{figure}

Here, the only difference with the previous example is that
the truth parameter field for~$m_1$ no longer has a
vertical discontinuity along the line~$x=0.5$ (see
figure~\ref{fig:coincide2-target}). 
\pgfplotstableread[col sep=comma]{fig/coincide2/centercs.csv}\centercsb
\pgfplotstableread[col sep=comma]{fig/coincide2/bottcs.csv}\bottcsb
In figure~\ref{fig:coincide2-target}(c-d), we again show the reconstructions
for parameters~$m_1$ and $m_2$ obtained by solving two independent inverse
problems~\eqref{eq:poisson}.  
The reconstructions for the four joint inverse problems are
shown in figure~\ref{fig:coincide2-reconstruct}, and the corresponding values of
the relative medium misfit are given in table~\ref{tab:all-med}.

As in the previous example, for $m_2$ the reconstructions obtained
with the different joint
inverse problems do not differ significantly (see
figure~\ref{fig:coincide2-reconstruct}b).
However, we observe differences among the reconstructions for parameter~$m_1$.
Using the cross-gradient only marginally improves the
reconstruction for parameter~$m_1$.  The use of the normalized
cross-gradient does not show improvement over the cross-gradient. As in the
first example, both the VTV joint regularization and the nuclear norm joint
regularization perform the best, and their corresponding reconstructions contain
all sharp interfaces present in the true image. However, in
figures~\ref{fig:coincide2-reconstruct}a (iii) and (iv) we also see
a vertical discontinuity not present in the true
image~\ref{fig:coincide2-target}c. This ghost interface in $m_1$ is due to the presence of
such a discontinuity in $m_2$, and highlights the tendency
of the VTV joint regularization and nuclear norm joint regularization to
superimpose discontinuities in both parameters. 
Note, however, that the amplitude of this ghost interface is small compared to
the amplitudes of the correctly recovered interfaces.

\subsection{Joint inversion of bulk modulus and density in the
  acoustic wave equation}
\label{sec:acoustic}

We now study a joint inverse problem of the
form~\eqref{eq:joint1}, i.e., both parameters enter the same equation,
namely the acoustic wave equation.

\subsubsection{Problem description}

We start by defining the forward problem, i.e., the acoustic wave
PDE.  The propagation of acoustic waves depends on the
bulk modulus~$\kappa$ 
and the density~$\rho$ of the medium
of propagation.  Let us define the acoustic pressure, $u(\xx,t) \coloneqq  -
\kappa(\xx) \nabla \cdotp \mathbf{u}(\xx,t)$, with $\mathbf{u}(\xx,t)$ the
displacement vector at location~$\xx$ and time~$t$.  The time-domain
acoustic wave equation
with first order absorbing boundary
condition~\cite{EngquistMajda77} and initial conditions at rest is
given by
\begin{equation} \label{eq:acoustic}
\begin{aligned} 
 \frac1\kappa \ddot{u} - \nabla \cdotp \left( \frac1\rho \nabla u
\right)  = f, \quad & \text{in } \Omega \times (0,T),  \\ 
 u(\xx,0) = \dot{u}(\xx,0)  = 0, \quad & \text{in } \Omega,  \\
 \frac1\rho \nabla u \cdotp \nn = 0, \quad & \text{on } \partial \Omega_n \times
(0,T), \\
 \frac1\rho \nabla u \cdotp \nn = -\frac1{\sqrt{\kappa \rho}} \dot{u}, \quad &
\text{on } \partial \Omega_a \times (0,T),
\end{aligned} \end{equation}
where $f$ is a forcing term, $\dot{u}$ and $\ddot{u}$ are the first
and second time derivatives of~$u$, and the boundary of the domain $\partial
\Omega$ is
partitioned as $\partial \Omega = \partial \Omega_a \cup
\partial \Omega_n$.
The acoustic wave velocity of the medium is given by $c$, with the
relation $\kappa = \rho c^2$.  The PDE in~\eqref{eq:acoustic} is the
variable density form of the acoustic wave equation; when the
density~$\rho$ is assumed constant, equation~\eqref{eq:acoustic} reduces to
$\frac1{c^2} \ddot{u} - \Delta u = \tilde{f}$.

Here, we assume that both the bulk modulus~$\kappa$ and the
density~$\rho$ are unknown.  Since they both appear
in~\eqref{eq:acoustic} through their inverse, we introduce the
parameters $\alpha \coloneqq 1/\kappa$ and $\beta \coloneqq 1/\rho$,
and formulate the inverse problem in terms of $\alpha$ and $\beta$.
As common in seismic inversion, we consider $N_s$~multiple experiments,
characterized by their forcing terms~$f_i$ and datasets $\dd_i$, which
corresponds to pointwise observations in space, recorded continuously in time.  The
acoustic wave inverse problem is then formulated as
\begin{equation} \label{eq:acousticinv}
\min_{\alpha, \beta>0} \left\{ \frac1{2N_s} \sum_{i=1}^{N_s} \int_0^T |B u_i(t) -
\dd_i(t)|^2 \, dt + \Rs(\alpha, \beta) \right\},
\end{equation}
where each $u_i$ solves the forward problem~\eqref{eq:acoustic} with forcing
term~$f_i$, 
\begin{equation*}
\begin{aligned}
  \alpha \ddot{u}_i - \nabla \cdotp ( \beta \nabla u_i ) = f_i,
\quad & \text{in } \Omega \times (0,T),  \\ 
  u_i(\xx,0) = \dot{u}_i(\xx,0)  = 0, \quad & \text{in } \Omega,  \\
 \beta \nabla u_i \cdotp \nn  = 0, \quad & \text{on } \partial \Omega_n \times
(0,T), \\
 \beta \nabla u_i \cdotp \nn = -\sqrt{\alpha \beta} \dot{u}_i, \quad &
\text{on } \partial \Omega_a \times (0,T).
\end{aligned} \end{equation*}
In our experiments,
the physical constraints $\alpha, \beta > 0$ are never active, and therefore
not enforced explicitly.

\subsubsection{Solution of the acoustic wave joint inverse problem}

Because 
the solution of the acoustic wave equation couples the
parameters~$\alpha$ and $\beta$,
the inverse
problem~\eqref{eq:acousticinv} could be regularized by two independent
TV regularizations, i.e., $\Rs(\alpha, \beta) = \Rs_\tve(\alpha) +
\Rs_\tve(\beta)$~\cite{EpanomeritakisAkccelikGhattasEtAl08}.  However,
the resulting problem can be difficult to solve and does not
incorporate 
the structural correlation that usually exists
between these
parameters due to the types of rock occurring in the subsurface.
Going beyond the use of ad-hoc methods to handle both parameters at
once, some researchers have addressed~\eqref{eq:acousticinv}
as a joint inverse
problem~\cite{ManukyanMaurerNuber16,LiLiangAbubakarEtAl13}.  Previous
attempts have used the cross-gradient term, but not its normalized
version, the VTV or the nuclear norm regularization.  In this section,
we study whether the use of joint regularization can improve reconstructions
for~$\alpha$ and~$\beta$.

In our numerical tests, we use 6~independent sources, $f_i(\xx,t)$,
located on the top boundary of the domain at 0.1, 0.25, 0.4, 0.6,
0.75, and 0.9 from the left boundary (yellow stars in
figure~\ref{fig:targetab}a); each source is a point source in space,
and a Ricker wavelet in time with a central frequency of 2~Hz.  The data
are recorded at 20~locations equally spaced along the top boundary
(green triangles in figure~\ref{fig:targetab}b), and polluted by
independent Gaussian noise with zero mean and variance corresponding
to a signal-to-noise ratio of 20~dB.  The boundary conditions are a
homogeneous Neumann boundary condition along the top boundary
$\partial \Omega_n = [0,1] \times \{1\}$, and an absorbing boundary
condition along the left, bottom, and right boundaries $\partial
\Omega_a = \{0,1\} \times [0,1] \cup [0,1] \times \{0\}$.  The truth
parameter fields for~$\alpha$ and~$\beta$ are shown in
figure~\ref{fig:targetab}; they correspond to an acoustic wave
velocity varying from 2km/s to 3km/s\footnote{The following units are
  used: distance in km, velocity in km/s, density in g/cm$^3$, and
  bulk modulus in GPa.},  typical values for a shallow subsurface 
(see for instance~\cite{MartinWileyMarfurt06,
  VirieuxOperto09}).
The finite-element mesh consists of 800 triangles ($h=1/20$).
The initial guesses for parameters~$\alpha$ and~$\beta$ are smoothed
versions of the truth parameters fields (see figure~\ref{fig:targetab}ii).
\newcommand{\wcleg}{0.2\textwidth}
\newcommand{\wcc}{0.22\textwidth}
\begin{figure}%
\centering
\begin{tabular}{lr@{\hspace{-.05in}}l@{\hspace{.01in}}c@{\hspace{.01in}}c}
\rotatebox[origin=l]{90}{\hspace{0.4in}(a) $\alpha$}
&
\includegraphics[height=\wcleg, trim=85 60 770 80, clip=true]
{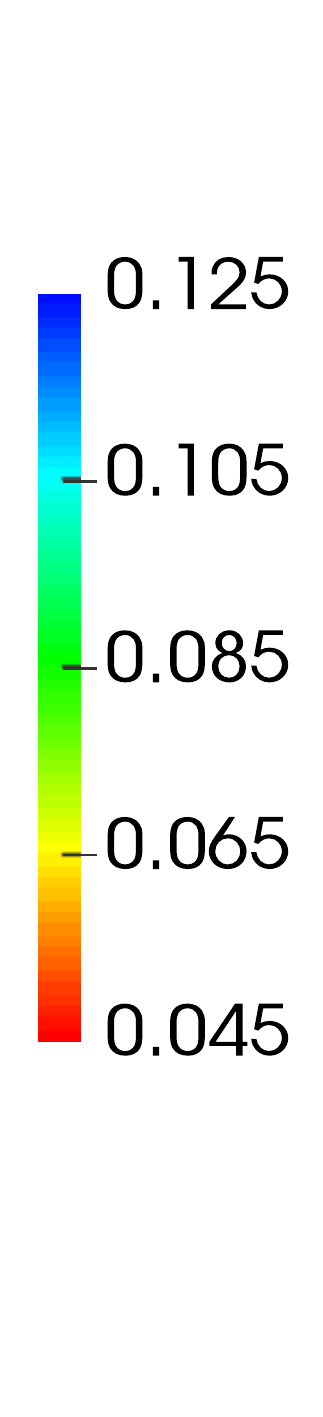}
&
\tikzsetnextfilename{targeta}
\begin{tikzpicture}
\node[anchor=south west, inner sep=0] (image1) at (0,0)
{\includegraphics[width=\wcc, trim=240 40 240 40, clip=true]
{fig/acoustic/a_target}};
\begin{scope}[x={(image1.south east)},y={(image1.north west)}]
\draw[draw=black] (0.1,1.0) node[star,fill=yellow,star points=9,scale=0.4,draw] {};
\draw[draw=black] (0.25,1.0) node[star,fill=yellow,star points=9,scale=0.4,draw] {};
\draw[draw=black] (0.4,1.0) node[star,fill=yellow,star points=9,scale=0.4,draw] {};
\draw[draw=black] (0.6,1.0) node[star,fill=yellow,star points=9,scale=0.4,draw] {};
\draw[draw=black] (0.75,1.0) node[star,fill=yellow,star points=9,scale=0.4,draw] {};
\draw[draw=black] (0.9,1.0) node[star,fill=yellow,star points=9,scale=0.4,draw] {};
\end{scope}
\end{tikzpicture}
&
\includegraphics[height=\wcc, width=\wcc, trim=300 80 300 80,
clip=true]{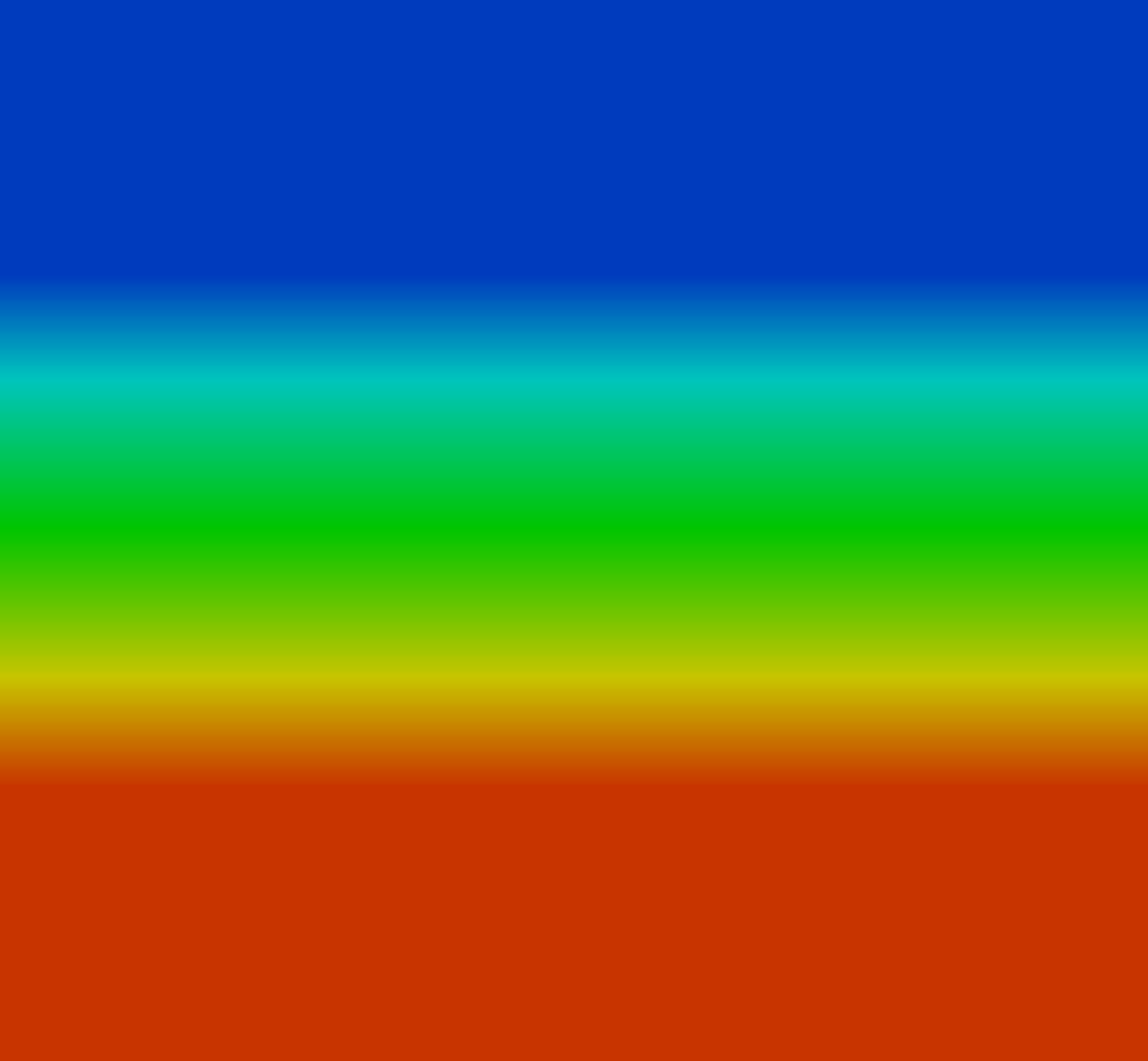}
&
\includegraphics[width=\wcc, trim=240 40 240 40, clip=true]
{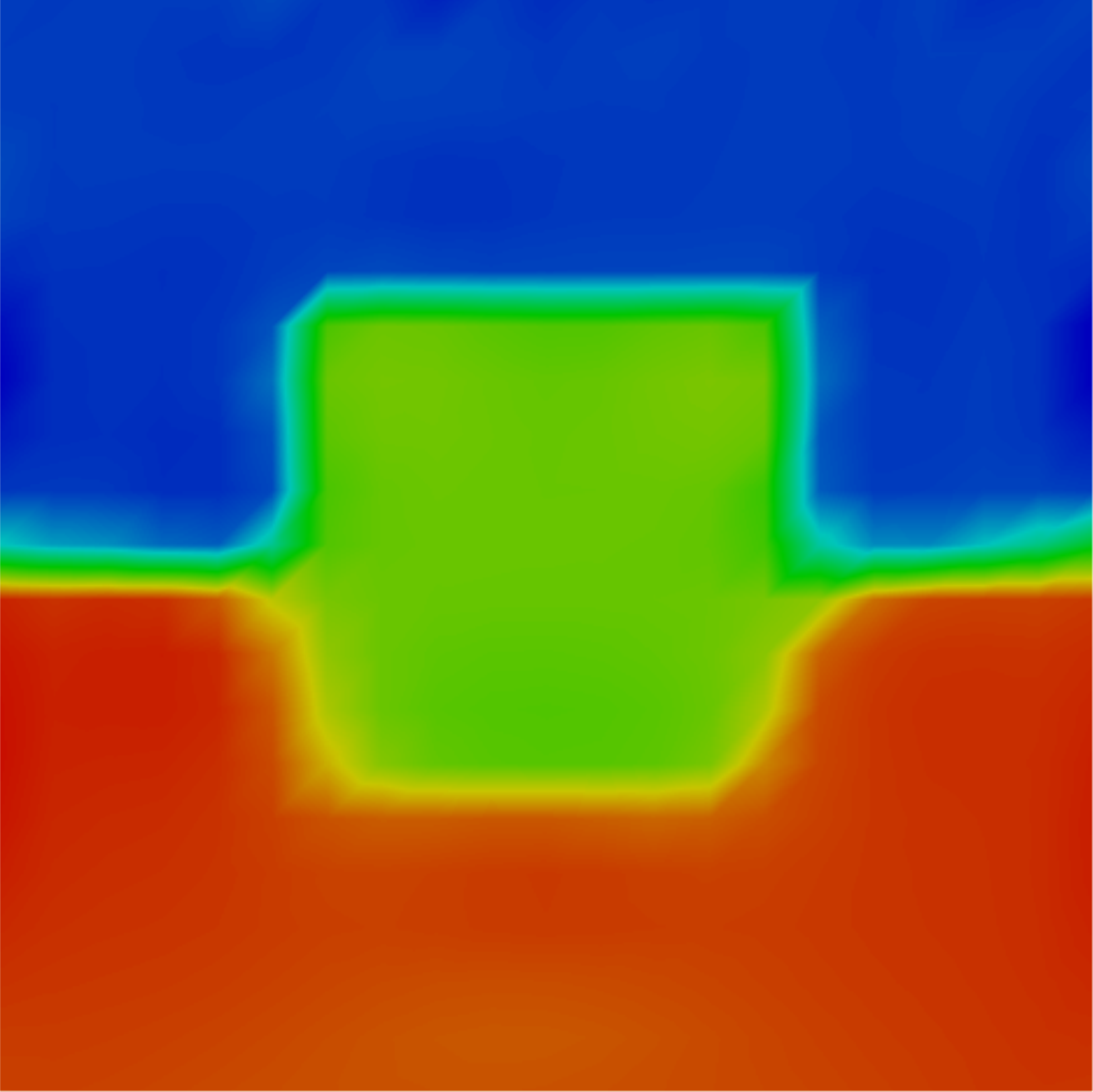}
\\
\rotatebox[origin=l]{90}{\hspace{0.4in}(b) $\beta$}
&
\includegraphics[height=\wcleg, trim=95 60 750 80, clip=true]
{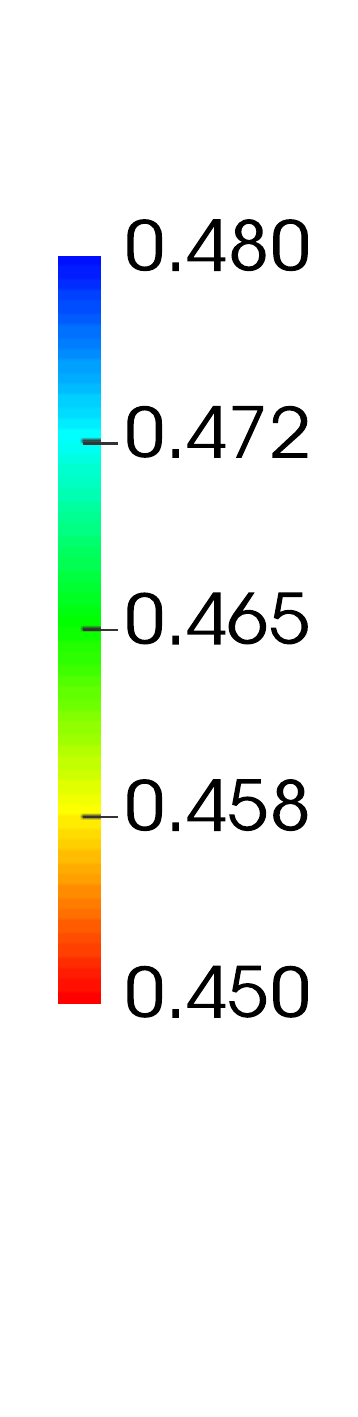}
&
\tikzsetnextfilename{targetb}
\begin{tikzpicture}
\node[anchor=south west, inner sep=0] (image1) at (0,0)
{\includegraphics[width=\wcc, trim=240 40 240 40, clip=true]
{fig/acoustic/b_target}};
\begin{scope}[x={(image1.south east)},y={(image1.north west)}]
\foreach \x in {1,2,...,20} {
\draw[draw=black] (0.047619*\x,1.03) node[regular polygon,regular polygon sides=3,scale=0.3,fill=green,rotate=180,draw] {};}
\end{scope}
\end{tikzpicture}
&
\includegraphics[height=\wcc, width=\wcc, trim=300 80 300 80, clip=true]
{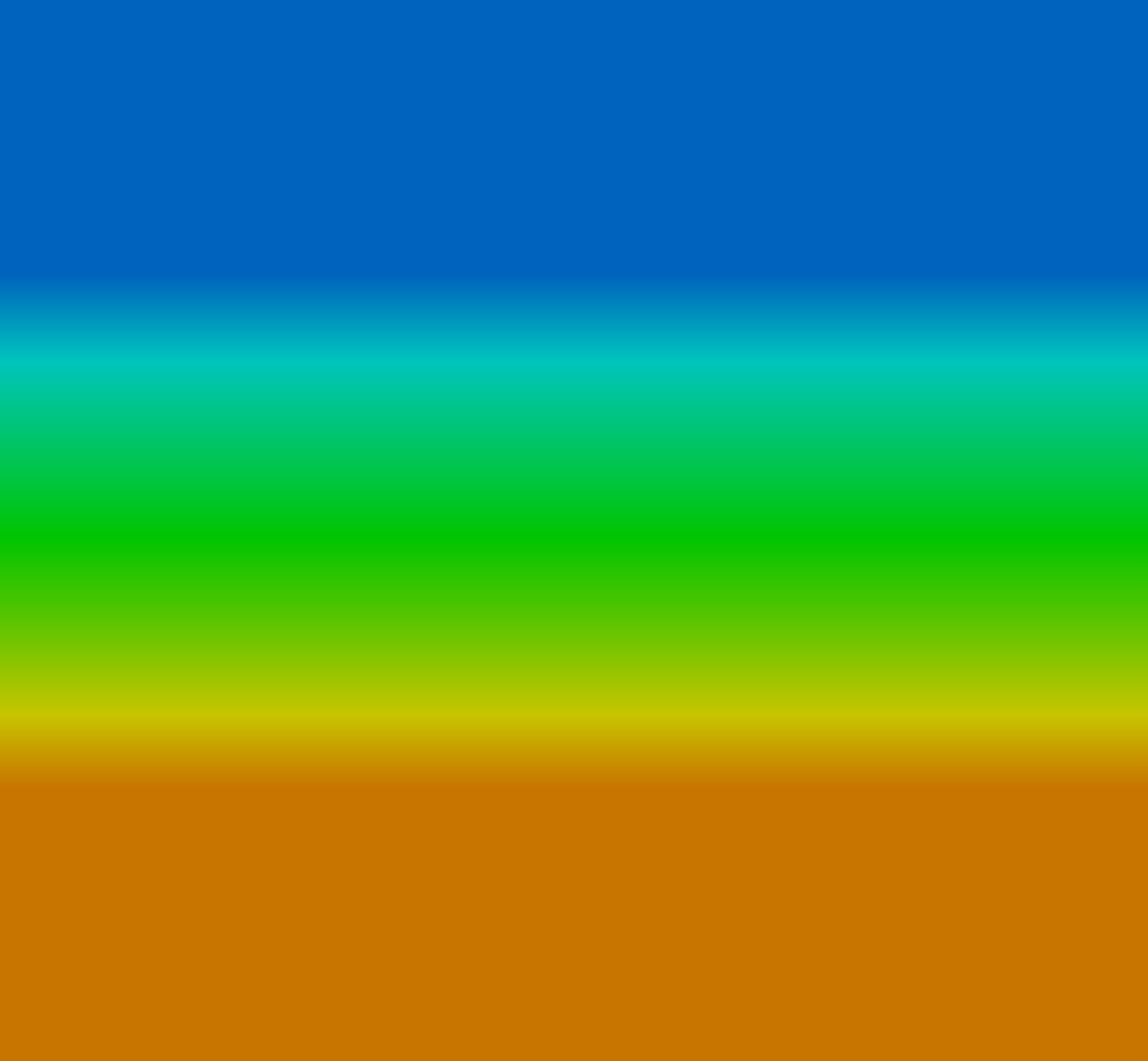} 
&
\includegraphics[width=\wcc, trim=240 40 240 40, clip=true]
{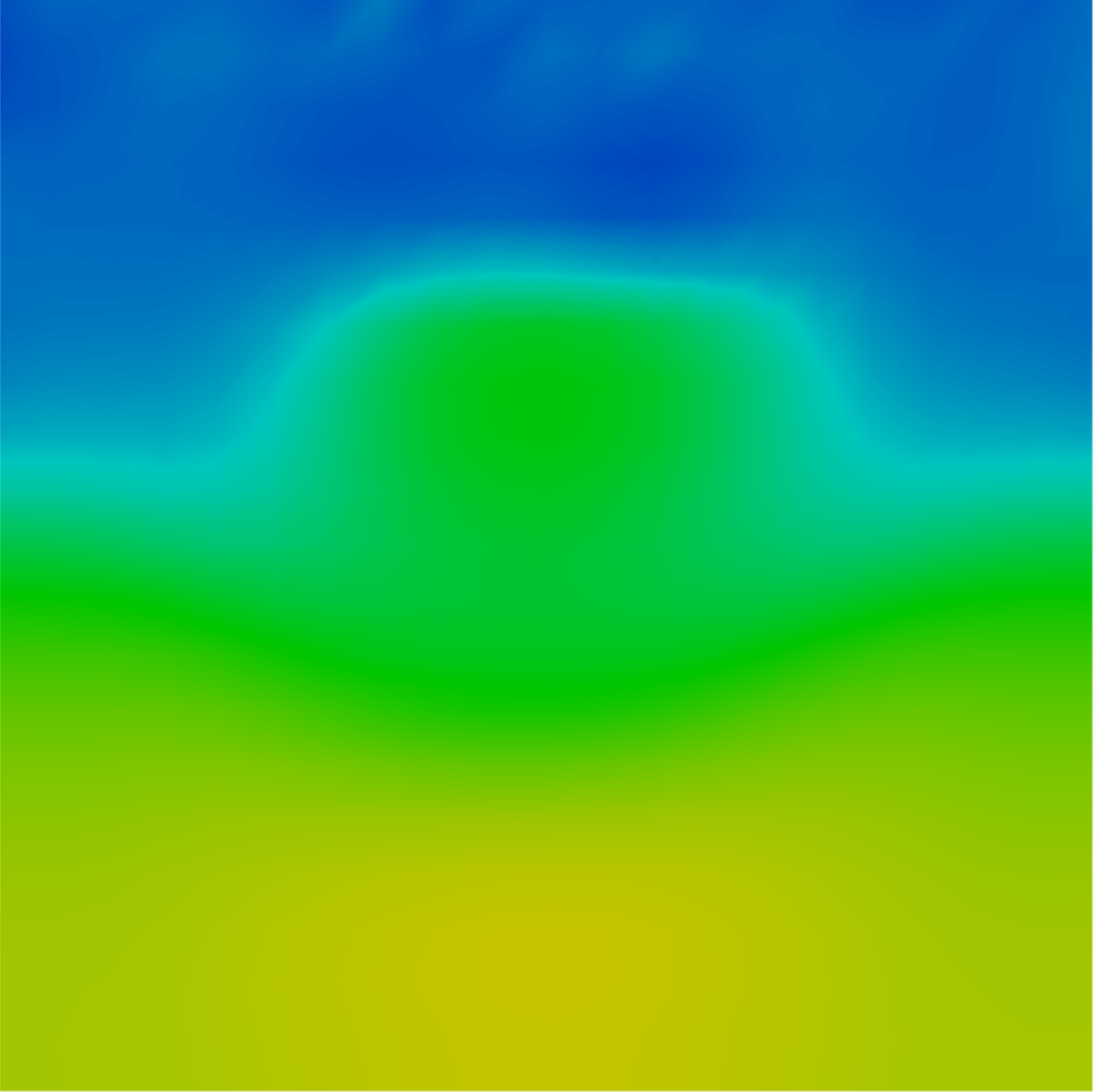} 
\\
& \multicolumn{2}{c}{(i) truth} 
& (ii) initial
&  (iii) independent 
\end{tabular}
\caption{Parameter fields (a) $\alpha$ and (b) $\beta$ in the joint acoustic inverse
problem~\eqref{eq:acousticinv}: (i) truth parameter fields, (ii) initial guesses,
and (iii) reconstructions when solving~\eqref{eq:acousticinv} regularized with
two independent TV regularizations ($\varepsilon = 10^{-3}$, $\gamma_\alpha=5
\cdotp 10^{-6}$, and $\gamma_\beta=9 \cdotp 10^{-6}$).  The yellow stars in (a-i)
and the green triangles in (b-i) indicate the locations of the point sources and
observations, respectively.}
\label{fig:targetab}
\end{figure}

\begin{figure}%
\centering
\begin{tabular}{l@{\hspace{.08in}}c@{\hspace{.01in}}c@{\hspace{.01in}}c@{\hspace{.01in}}c}
\rotatebox[origin=l]{90}{\hspace{0.4in}(a) $\alpha$}
&
\includegraphics[width=\wca, trim=290 70 290 70 , clip=true]
{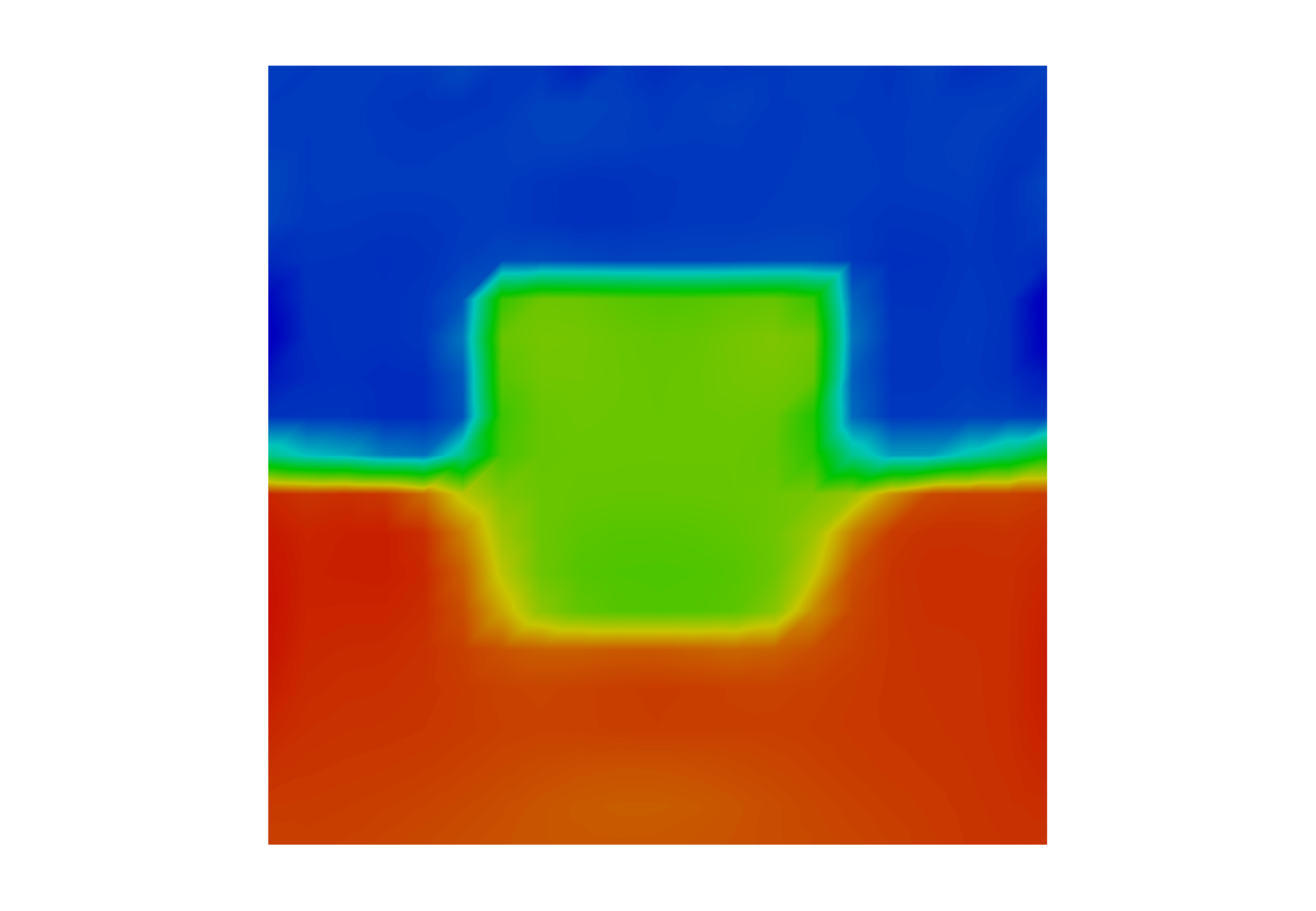}
& 
\includegraphics[width=\wca, trim=290 70 290 70, clip=true]
{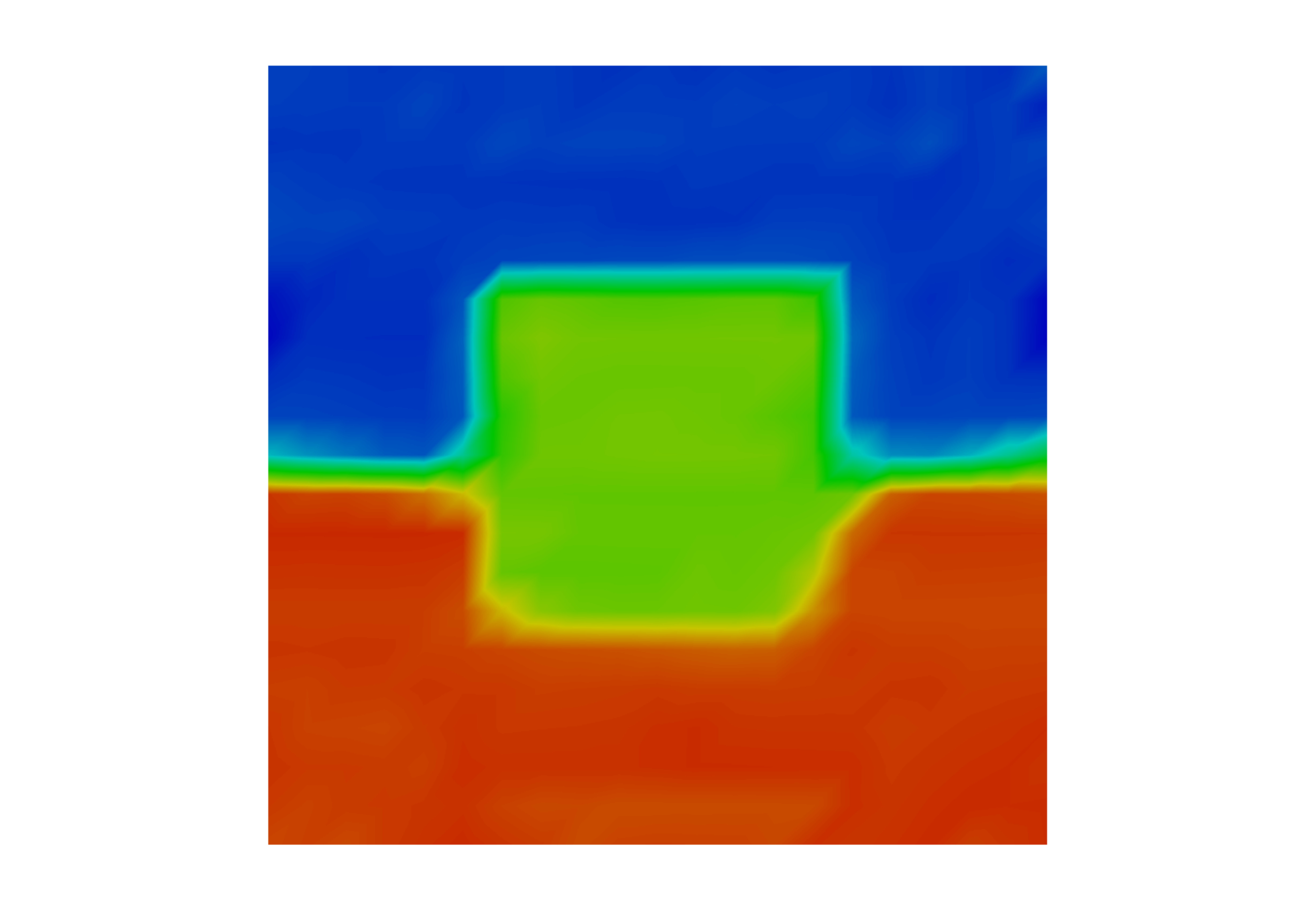}
& 
\includegraphics[width=\wca, trim=240 40 240 40, clip=true]
{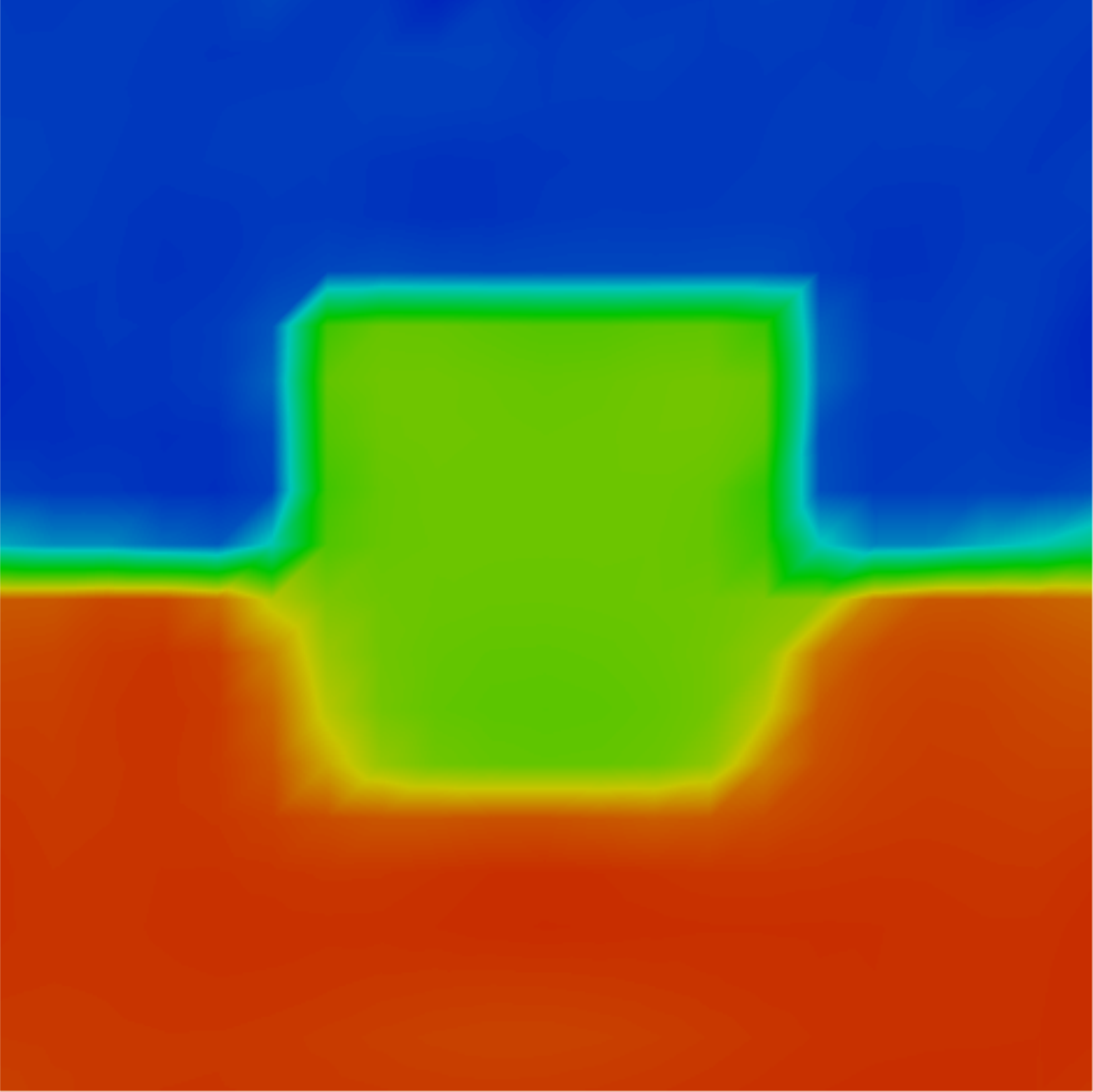}
& 
\includegraphics[width=\wca, trim=240 40 240 40, clip=true]
{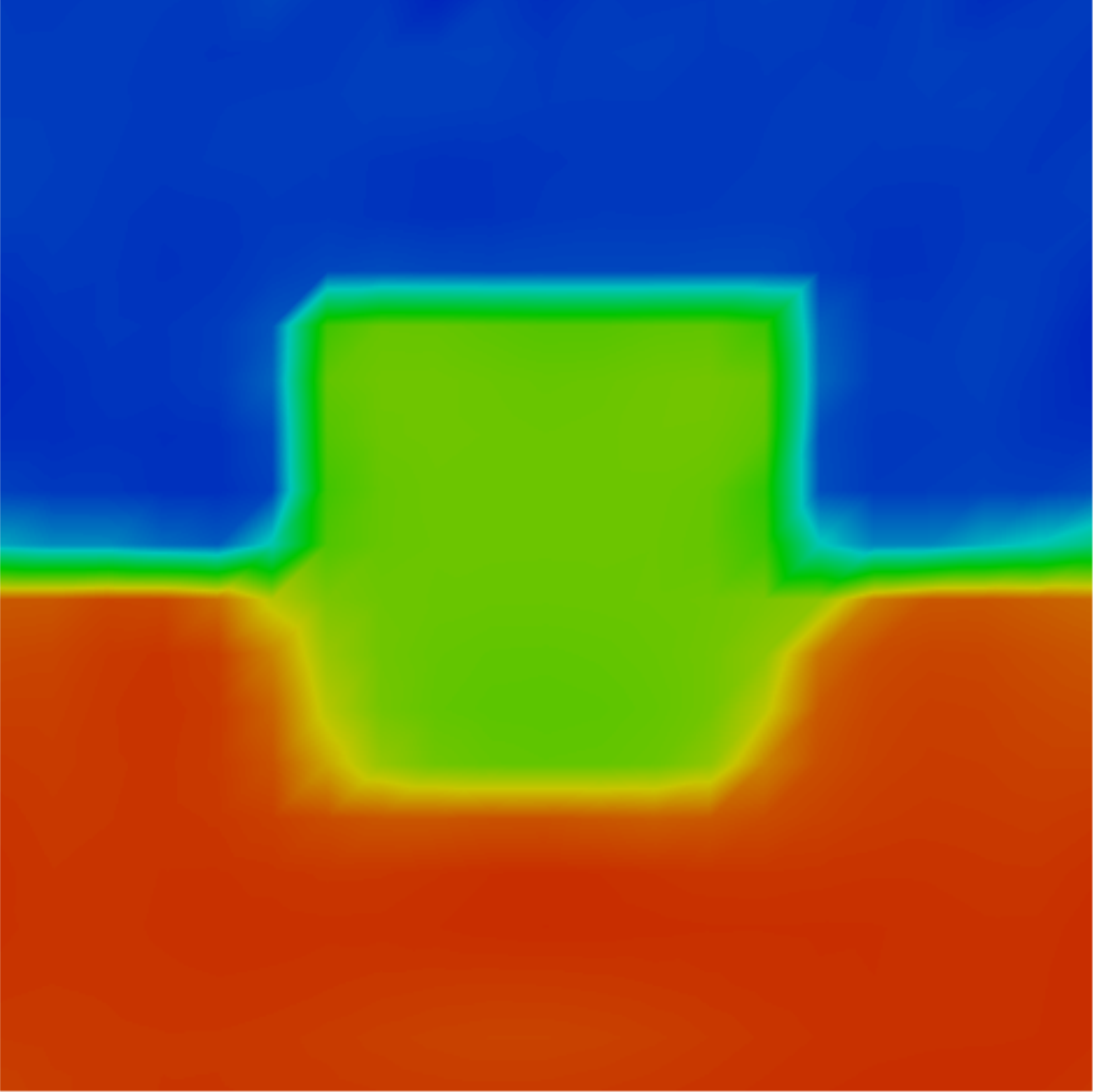}
\\
\rotatebox[origin=l]{90}{\hspace{0.4in}(b) $\beta$}
&
\includegraphics[width=\wca, trim=290 70 290 70, clip=true]
{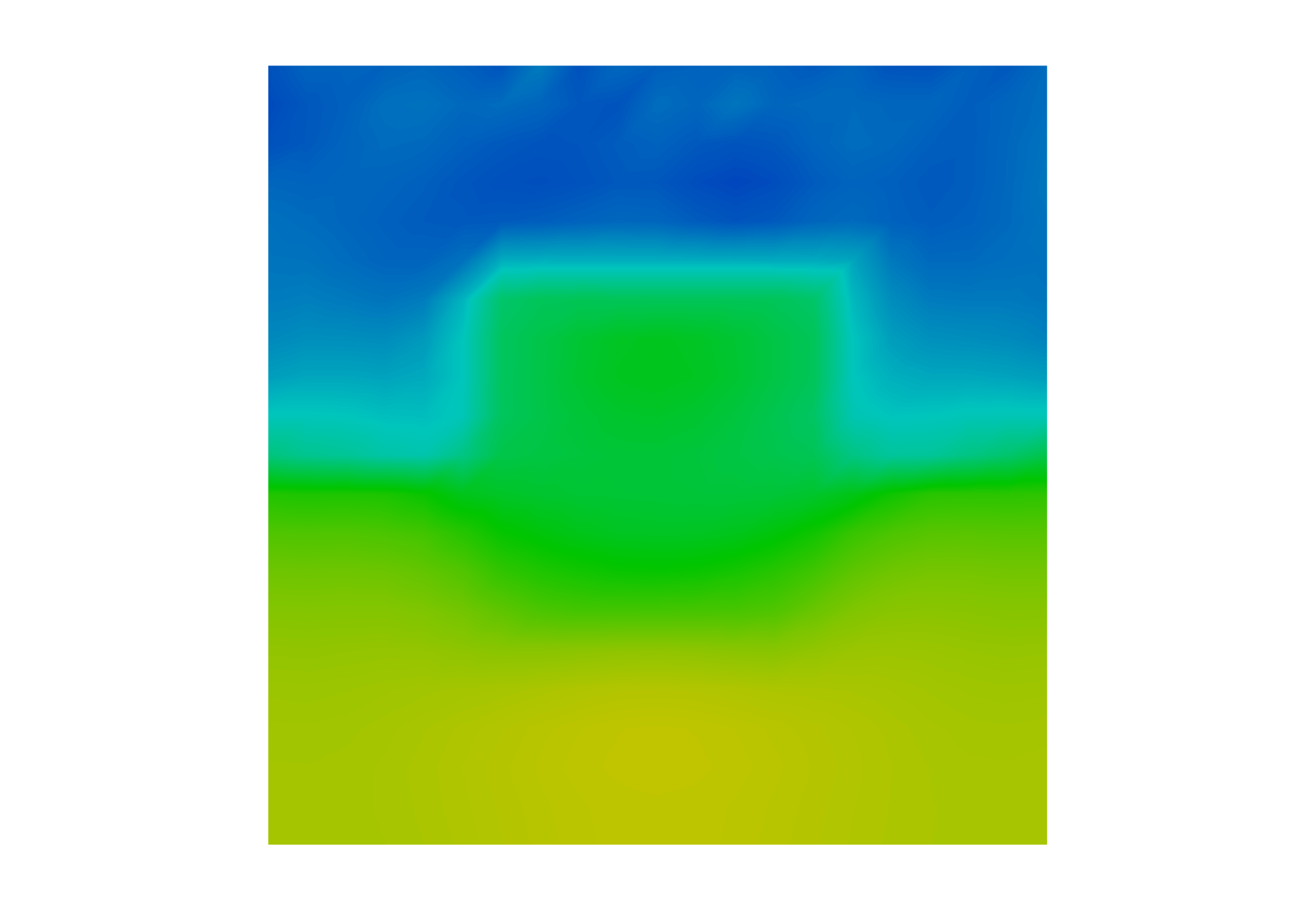} 
& 
\includegraphics[width=\wca, trim=290 70 290 70, clip=true]
{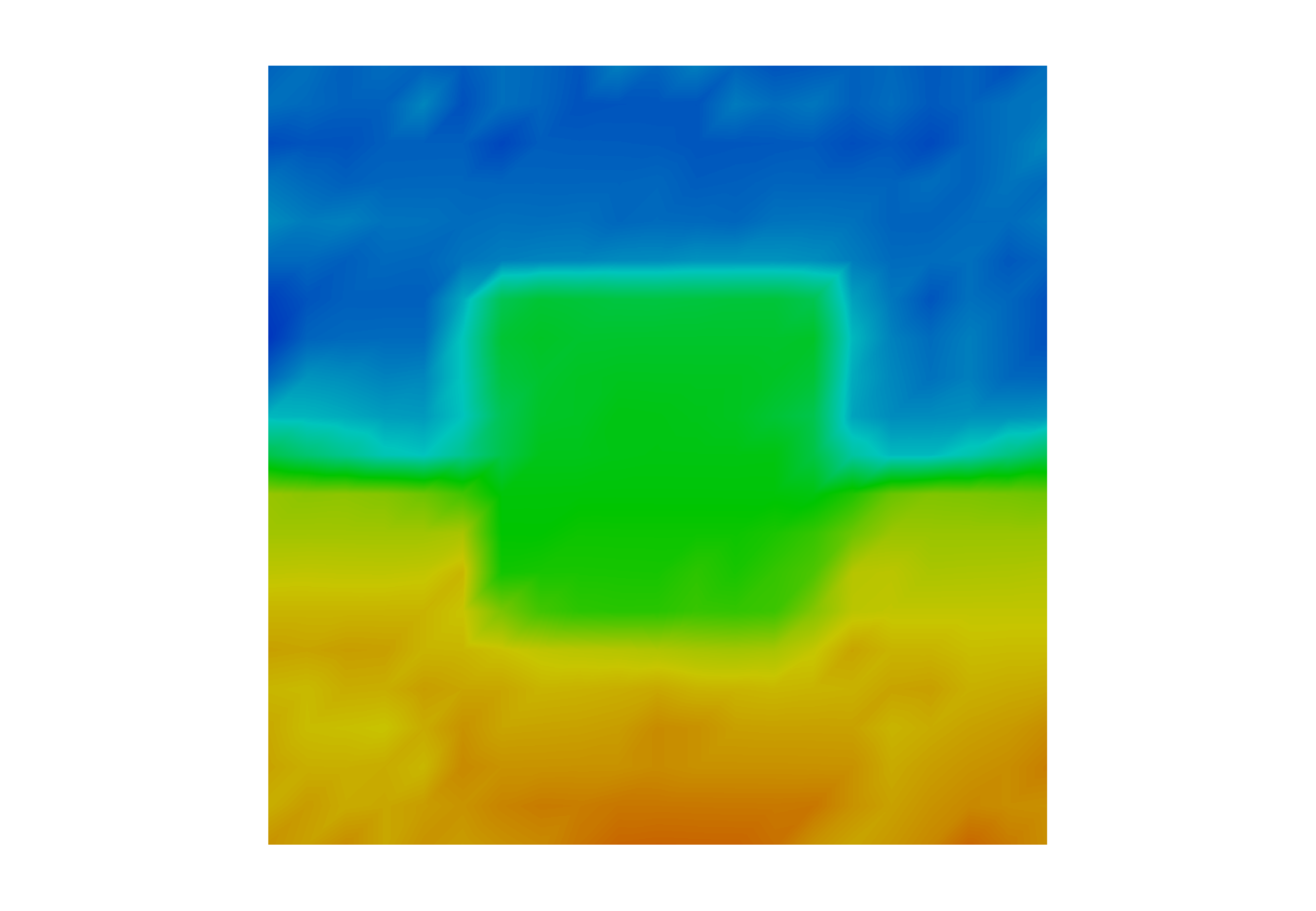} 
&
\includegraphics[width=\wca, trim=240 40 240 40, clip=true]
{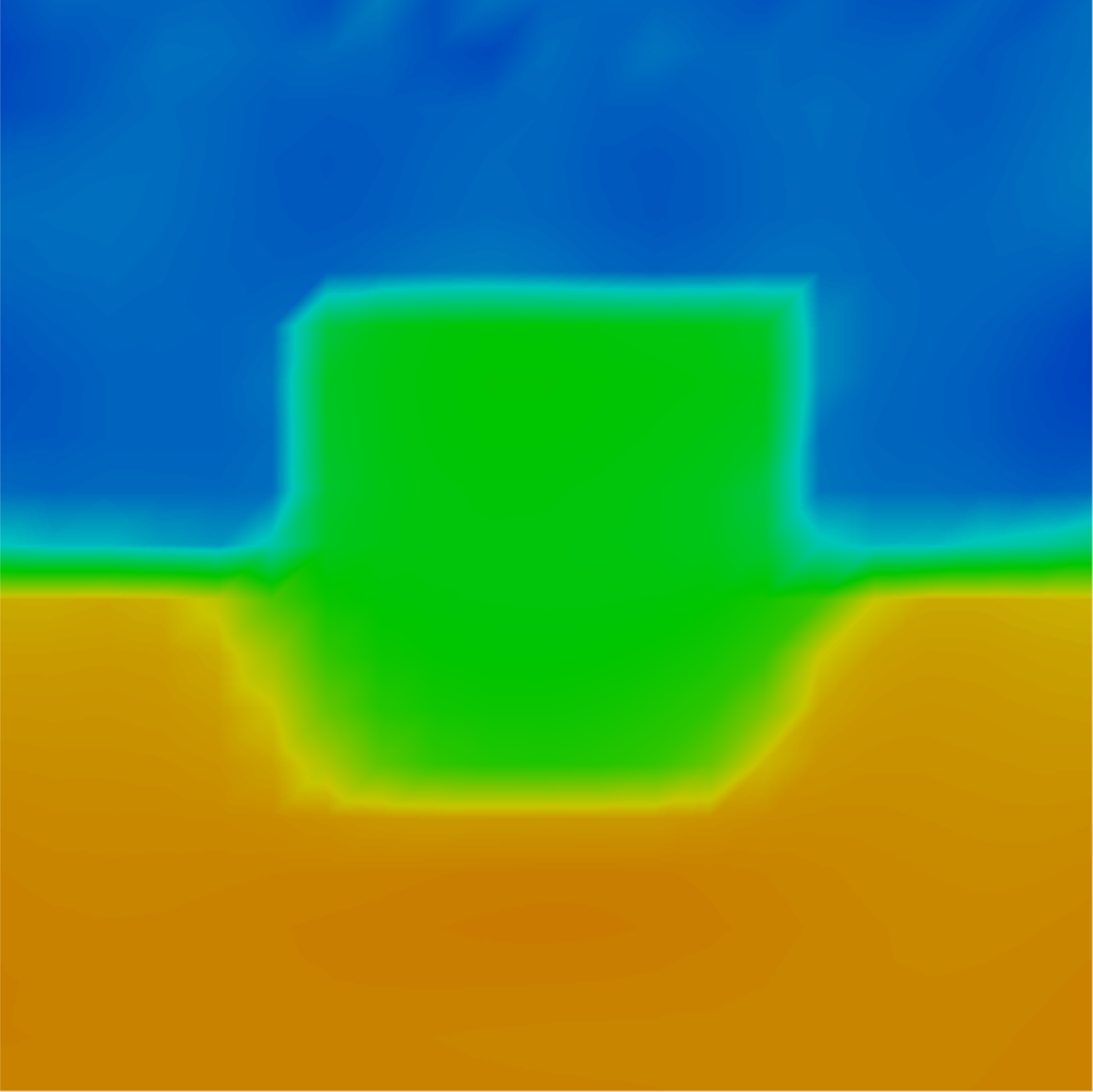}
& 
\includegraphics[width=\wca, trim=240 40 240 40, clip=true]
{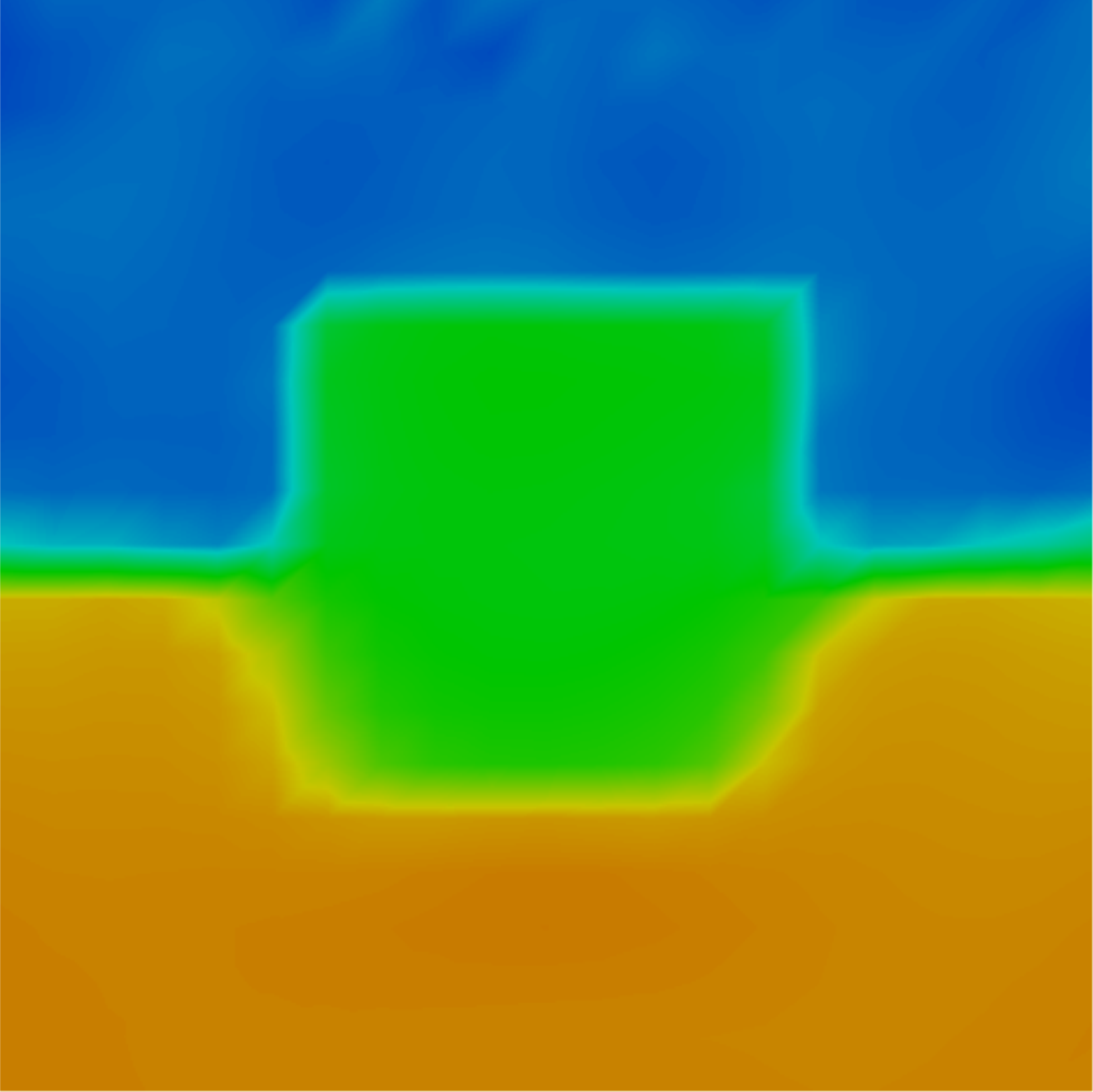}
\\
& (i) cross-gradient & (ii) norm. cross-gd &
(iii) vectorial TV & (iv) nuclear norm 
\end{tabular}
\caption{Reconstructions for the parameter fields (a) $\alpha$ and (b)
  $\beta$, obtained by solving~\eqref{eq:acousticinv}
  regularized with (i) the cross-gradient ($\gamma= 10^{-2}$) combined
  with two independent TV regularizations, (ii) the normalized
  cross-gradient ($ \gamma=9 \cdotp 10^{-6}$ and
  $\varepsilon=10^{-6}$) combined with the same independent TV
  regularizations, (iii) the VTV joint regularization ($\gamma=7
  \cdotp 10^{-6}$ and $\varepsilon=10^{-3}$), and (iv) the nuclear
  norm joint regularization ($\gamma=7 \cdotp 10^{-6}$ and
  $\varepsilon=10^{-3}$).  The parameters for the independent TV
  regularizations are the ones selected for the independent inverse
  problems (see caption in figure~\ref{fig:targetab}).  The
  legend is as in figure~\ref{fig:targetab}.}
\label{fig:acoustic-reconstruct}
\end{figure}

In figure~\ref{fig:targetab}iii, we show the reconstructions
of parameters~$\alpha$ and $\beta$ obtained by solving~\eqref{eq:acousticinv}
with independent TV regularizations.  
Whereas parameter~$\alpha$ is well reconstructed, the reconstruction
for $\beta$ is rather poor.
We next solve~\eqref{eq:acousticinv} with the proposed joint
regularization terms.
The results are
shown in figure~\ref{fig:acoustic-reconstruct}, and the corresponding values of
the relative medium misfit are given in table~\ref{tab:all-med}.

The different reconstructions for~$\alpha$
(figure~\ref{fig:acoustic-reconstruct}a) do not differ significantly from each
other.
However, the use of joint regularization improves the quality of the
reconstruction for the parameter~$\beta$
(figure~\ref{fig:acoustic-reconstruct}b). Whereas the use of the
cross-gradient only results in marginal improvement compared to the
reconstruction in figure~\ref{fig:targetab}d, the use of the
normalized cross-gradient allows recovery of the interfaces more
clearly. The best reconstructions are obtained using the VTV
or the nuclear norm joint regularizations.

\subsection{Joint inverse problem with different physics}
\label{sec:poissonacoustic}

\begin{figure}%
\centering
\begin{tabular}{r@{}cc}
\includegraphics[height=\wca, trim=100 150 1180 150, clip=true]
{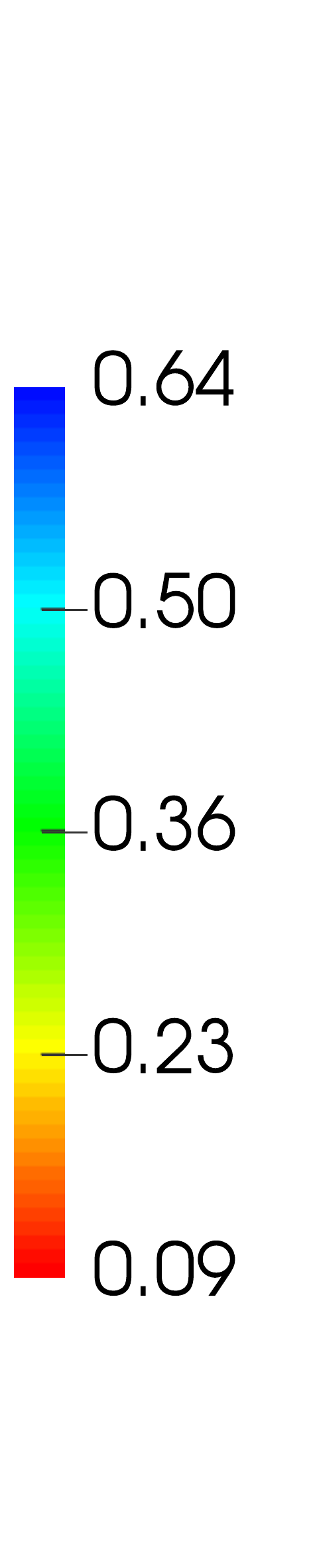} &
\tikzsetnextfilename{poissonacoustic-poisson}
\begin{tikzpicture}
\node[anchor=south west, inner sep=0] (image1) at (0,0)
{\includegraphics[width=\wca, trim=290 70 290 70, clip=true]
{fig/acousticelliptic/targetpoisson}};
\begin{scope}[x={(image1.south east)},y={(image1.north west)}]
\foreach \x in {1,2,...,20} {
\foreach \y in {1,2,...,20} {
\filldraw[fill=black,draw=white] (0.047619*\x,0.047619*\y) circle (0.002);}}
\end{scope}
\end{tikzpicture}
&
\includegraphics[width=\wca, trim=290 70 290 70, clip=true]
{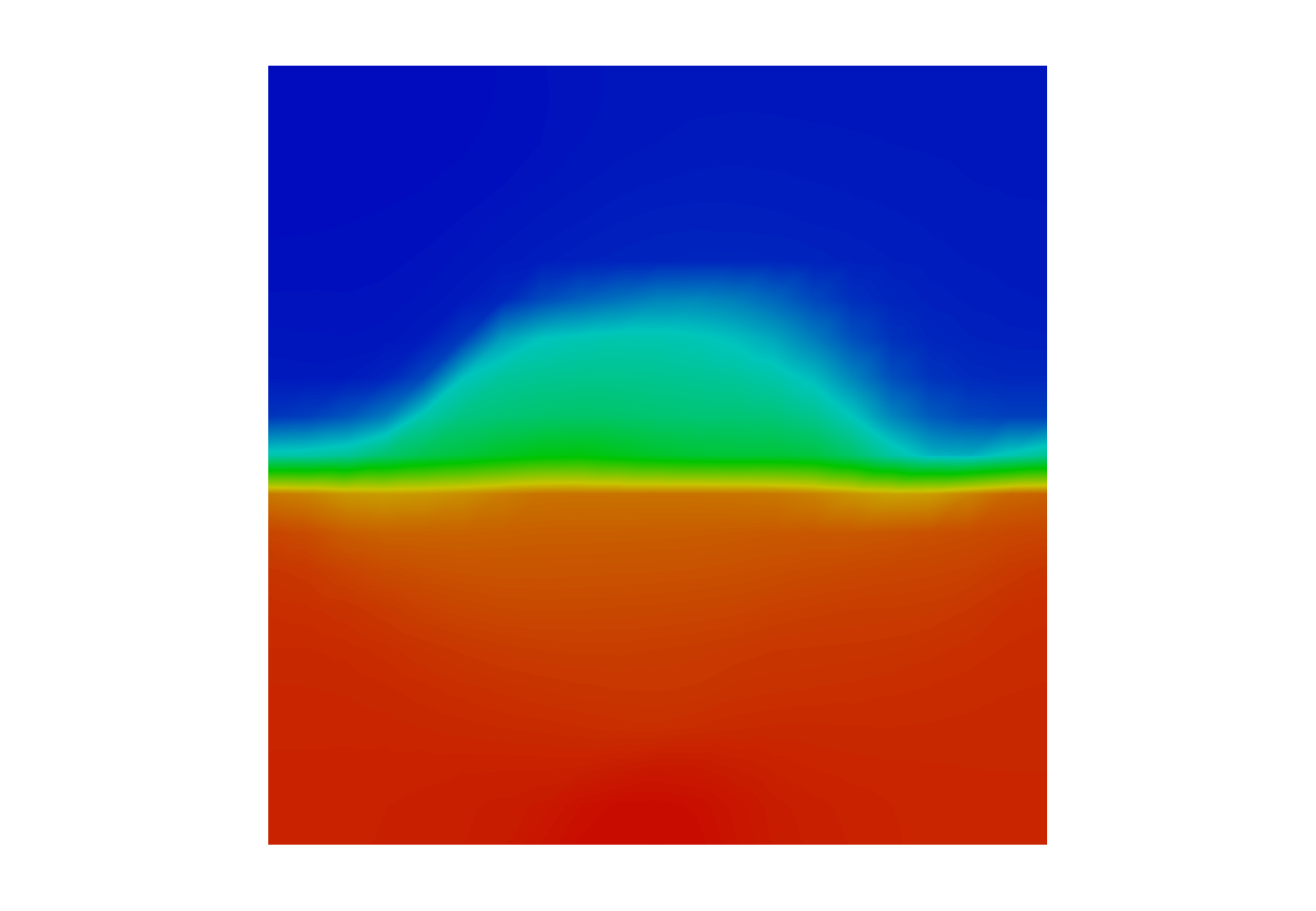} \\
 & (a) truth & (b) reconstruction
\end{tabular}
\caption{Plots of (a) truth parameter field for~$m$ in the Poisson inverse
problem~\eqref{eq:poisson2}, and (b) its reconstruction ($\gamma_m=2\cdotp
10^{-8}$ and $\varepsilon=10^{-3}$) with initial parameter field set to a
constant value of~$0.625$. The white dots in (a) indicate the location of the
pointwise observations.}
\label{fig:poisson}
\end{figure}

\begin{figure}%
\centering
\begin{tabular}{r@{}ccc}
\includegraphics[height=\wca, trim=100 150 1180 150, clip=true]
{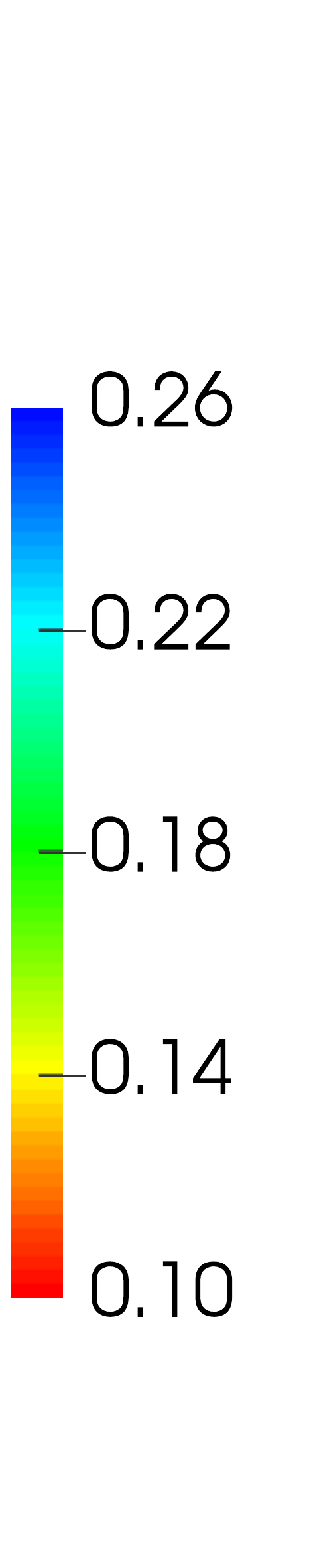} &
\tikzsetnextfilename{poissonacoustic-acoustic}
\begin{tikzpicture}
\node[anchor=south west, inner sep=0] (image1) at (0,0)
{\includegraphics[width=\wca, trim=290 70 290 70, clip=true]
{fig/acousticelliptic/targetacoustic}};
\begin{scope}[x={(image1.south east)},y={(image1.north west)}]
\foreach \x in {1,2,...,20} {
\draw[draw=black] (0.047619*\x,1.03) node[regular polygon,regular polygon sides=3,scale=0.3,fill=green,rotate=180,draw] {};}
\draw[draw=black] (0.5,0.1) node[star,fill=yellow,star points=9,scale=0.4,draw] {};
\end{scope}
\end{tikzpicture}
&
\includegraphics[width=\wca, trim=290 70 290 70, clip=true]
{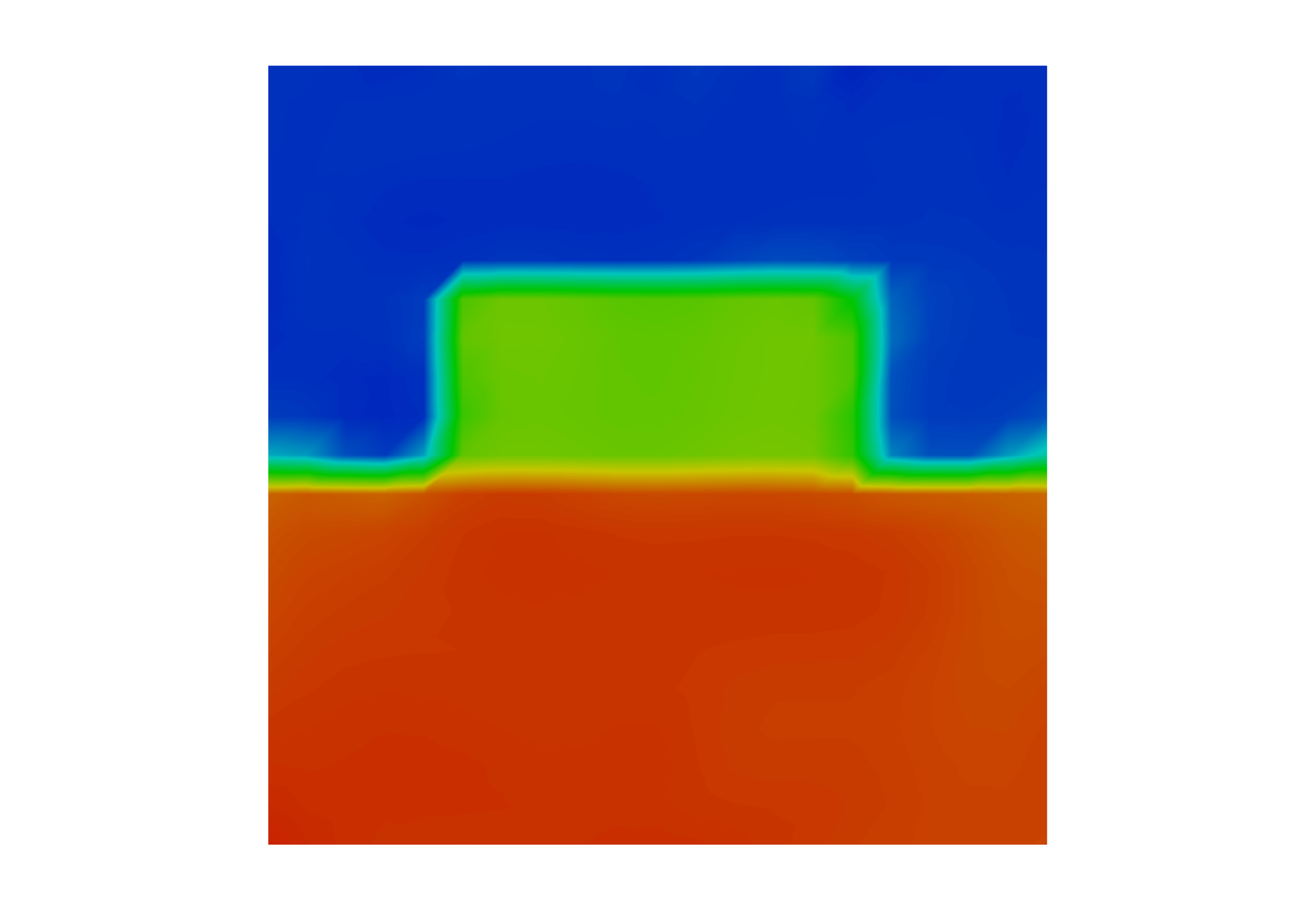} &
\includegraphics[width=\wca, trim=290 70 290 70, clip=true]
{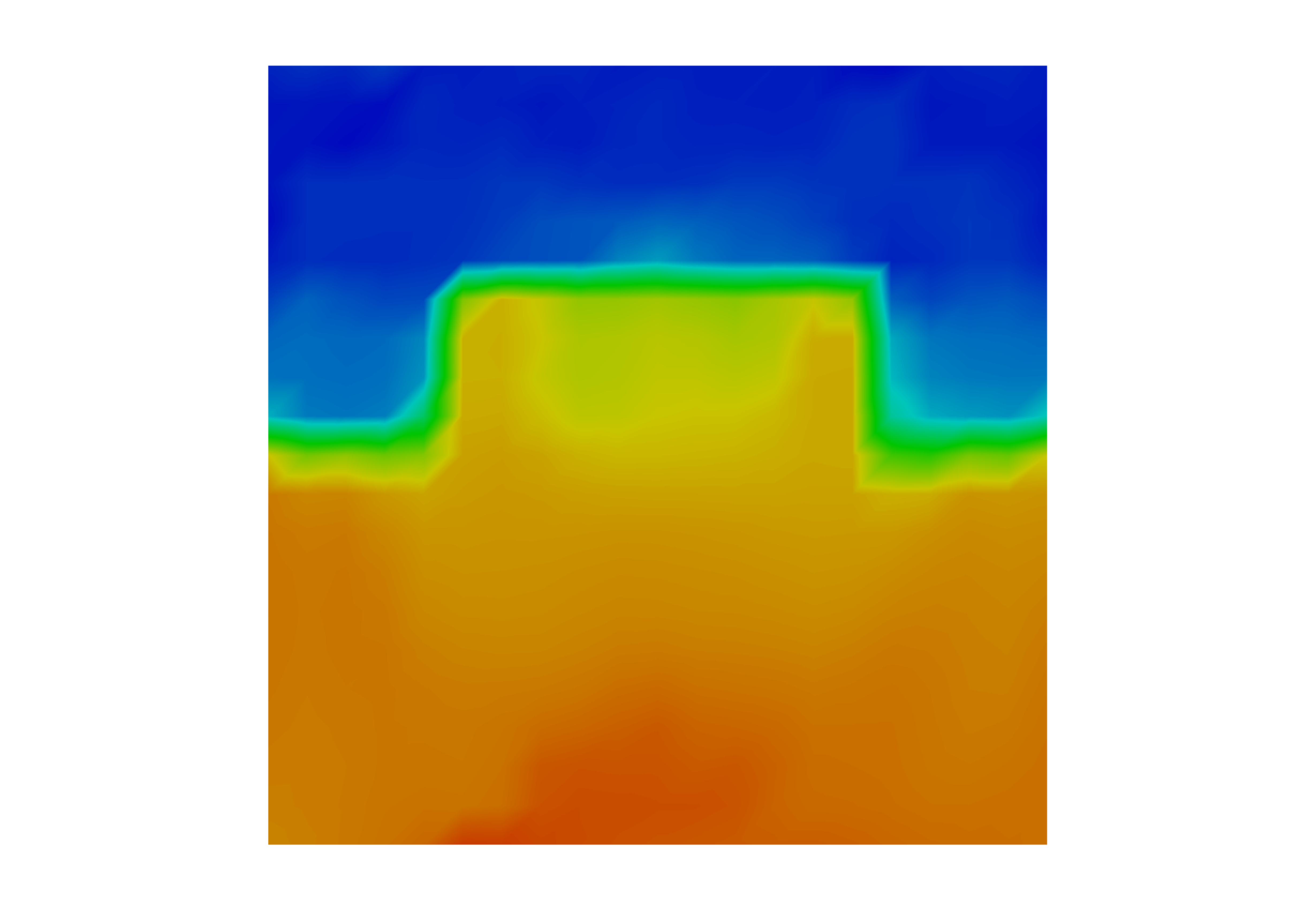} \\
 & (a) truth & (b) reconstruction 2~Hz & (c) reconstruction 4~Hz
\end{tabular}
\caption{Plots of (a) truth parameter field for~$\alpha$ in the acoustic inverse
problem~\eqref{eq:acousticinv2}, and (b,c) its reconstructions ($\gamma_\alpha=3
\cdotp 10^{-8}$ and $\varepsilon=10^{-3}$) with initial value for the parameter
field set to~$0.25$, and a source~$f_\alpha$ of frequency (b) 2~Hz, and (c) 4~Hz.
The green triangles in (a) indicate the locations of the pointwise observations,
and the yellow star in (a) indicates the location of the source.}
\label{fig:acoustic}
\end{figure}

As a last problem, we study a joint inverse problem~\eqref{eq:joint2}
governed by two different physics models; namely, we combine a Poisson
inverse problem and an acoustic wave inverse problem (assuming the
density~$\rho$ is known).  
This inverse problem is intended as a model problem for joint
seismic-electromagnetic inversion in the electromagnetic low frequency limit.
The Poisson inverse problem is identical to
the one used in section~\ref{sec:poissonpoisson},
\begin{equation} \label{eq:poisson2} \begin{gathered}
\min_m \left\{ \frac12 | Bu - \dd |^2  + \gamma_m
\int_\Omega \sqrt{ |\nabla m|^2 + \varepsilon} \, dx
\right\}, \quad \text{where} \\
\left\{ \begin{aligned}
- \nabla \cdotp ( e^m \nabla u) & = 1 \, \text{ in } \Omega, \\
u & = 0  \, \text{ on } \partial \Omega.
\end{aligned} \right.
\end{gathered} \end{equation}
The observation operator $B$ extracts the state $u$ at $20 \times 20$
equally distributed points over the entire domain (white dots in
figure~\ref{fig:poisson}a).  The data are polluted with 1\%~Gaussian
noise.  For the acoustic wave inverse problem, we set $\beta \equiv
1$, and invert only for the parameter~$\alpha = 1/\kappa = 1/c^2$,
\begin{equation} \label{eq:acousticinv2} 
\begin{gathered}
\min_{\alpha} \left\{ \frac1{2}
\int_0^T |B u(t) - \dd(t)|^2 \, dt + \gamma_\alpha \int_\Omega
\sqrt{ |\nabla \alpha|^2 + \varepsilon} \right\} \, dx, \quad
\text{where } \\
\left\{
\begin{aligned}
\alpha \ddot{u} - \Delta u = f_\alpha,
\quad & \text{in } \Omega \times (0,T), \\ 
u(\xx,0) = \dot{u}(\xx,0)  = 0, \quad & \text{in } \Omega, \\
\nabla u \cdotp \nn = 0, \quad & \text{on } \partial \Omega_n \times
(0,T), \\
\nabla u \cdotp \nn = -\sqrt{\alpha} \dot{u}, \quad &
\text{on } \partial \Omega_a \times (0,T).
\end{aligned}  \right.
\end{gathered}
\end{equation}
We use a single source~$f_\alpha$ with frequency 2~Hz or 4~Hz, located at~$(0.5, 0.1)$
(yellow star in figure~\ref{fig:acoustic}a), and 20~pointwise observations
equally spaced along the top boundary (green triangles in
figure~\ref{fig:acoustic}a).
The boundary conditions, the noise level in the data, the mesh, and the
numerical discretization are as in section~\ref{sec:acoustic}.
The initial guess for the Poisson parameter field~$m$ (resp.~for the acoustic
parameter field~$\alpha$) is set to a constant field with value $0.625$ (resp.~$0.25$),
corresponding to the value in the upper layer of the truth parameter
field, in
blue in figure~\ref{fig:poisson}a (resp.~figure~\ref{fig:acoustic}a).

\begin{figure}%
\centering
\begin{tabular}{l@{\hspace{.08in}}c@{\hspace{.01in}}c@{\hspace{.01in}}c@{\hspace{.01in}}c}
\rotatebox[origin=l]{90}{\hspace{0.3in}(a) $m$}
&
\includegraphics[height=\wca, trim=290 70 290 70 , clip=true]
{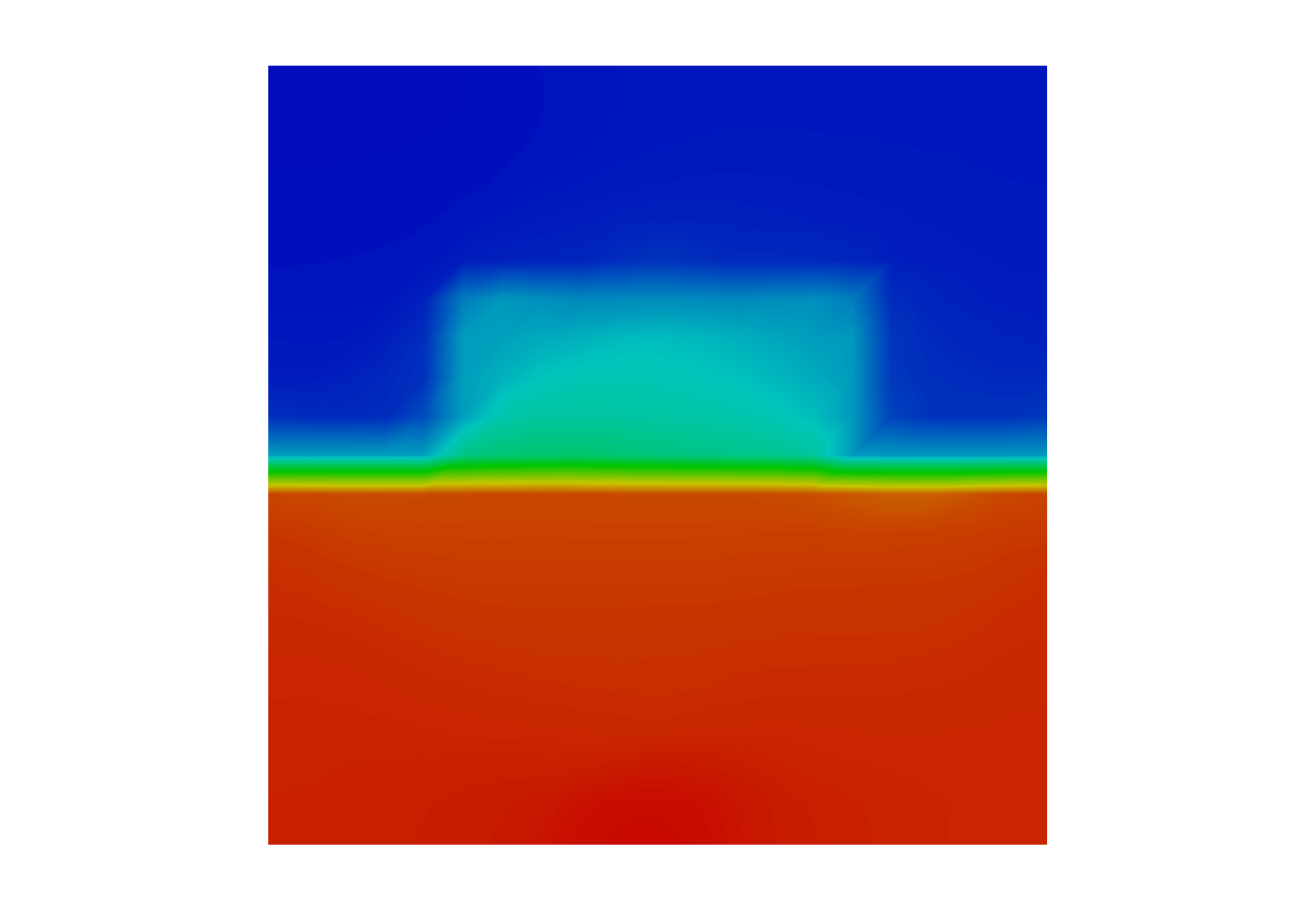}
& 
\includegraphics[height=\wca, trim=310 70 310 70, clip=true]
{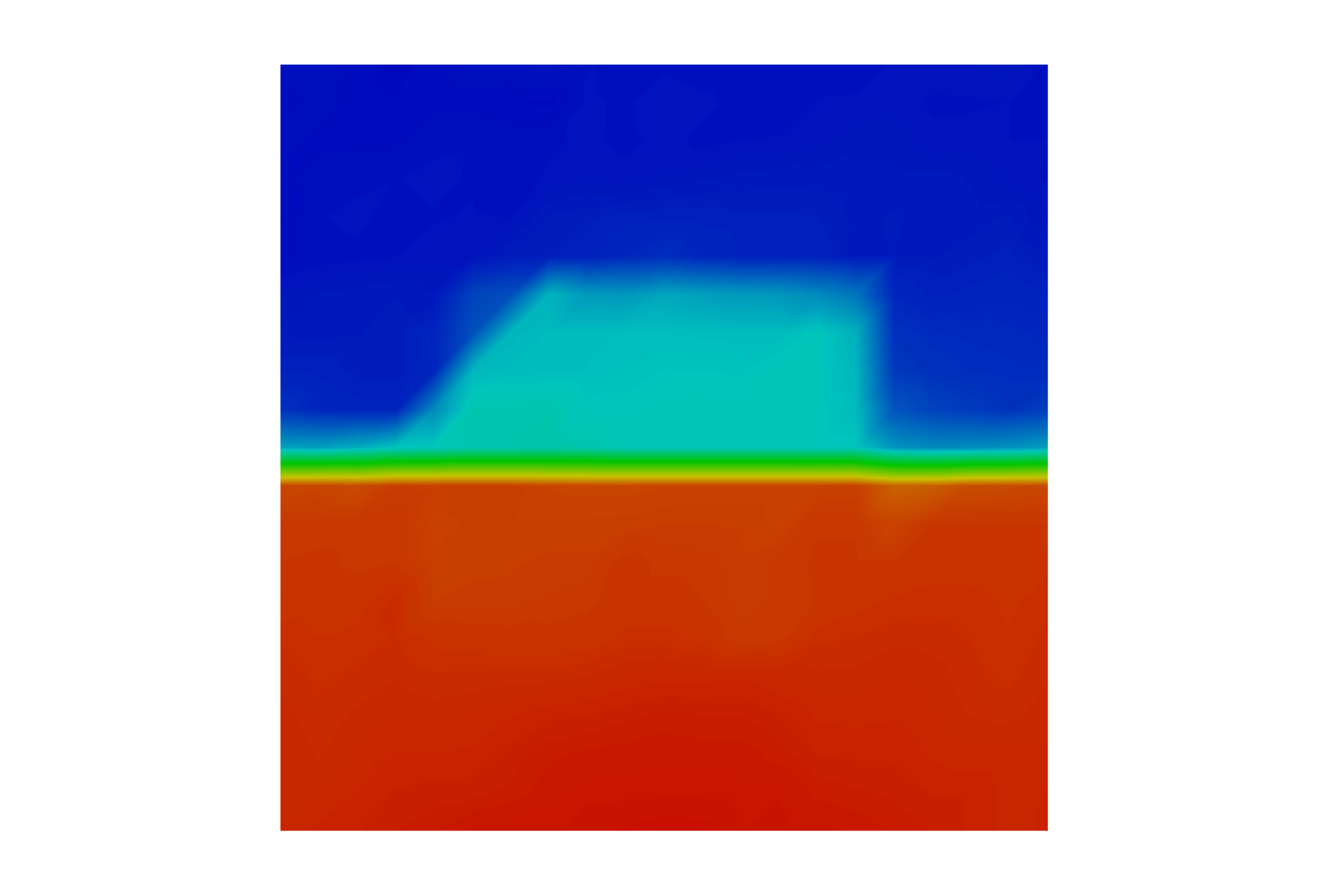}
& 
\includegraphics[height=\wca, trim=290 70 290 70, clip=true]
{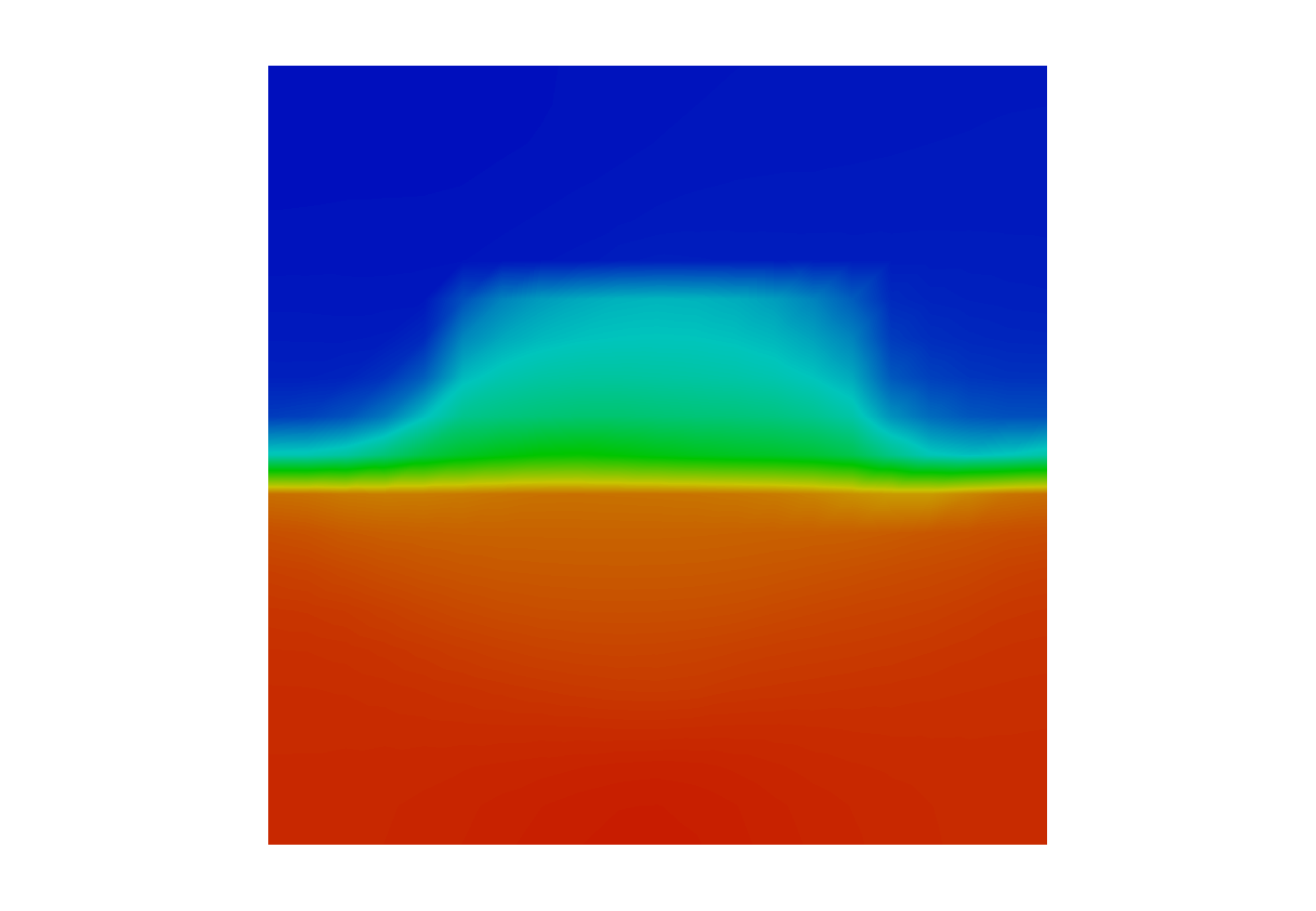}
& 
\includegraphics[height=\wca, trim=290 70 290 70, clip=true]
{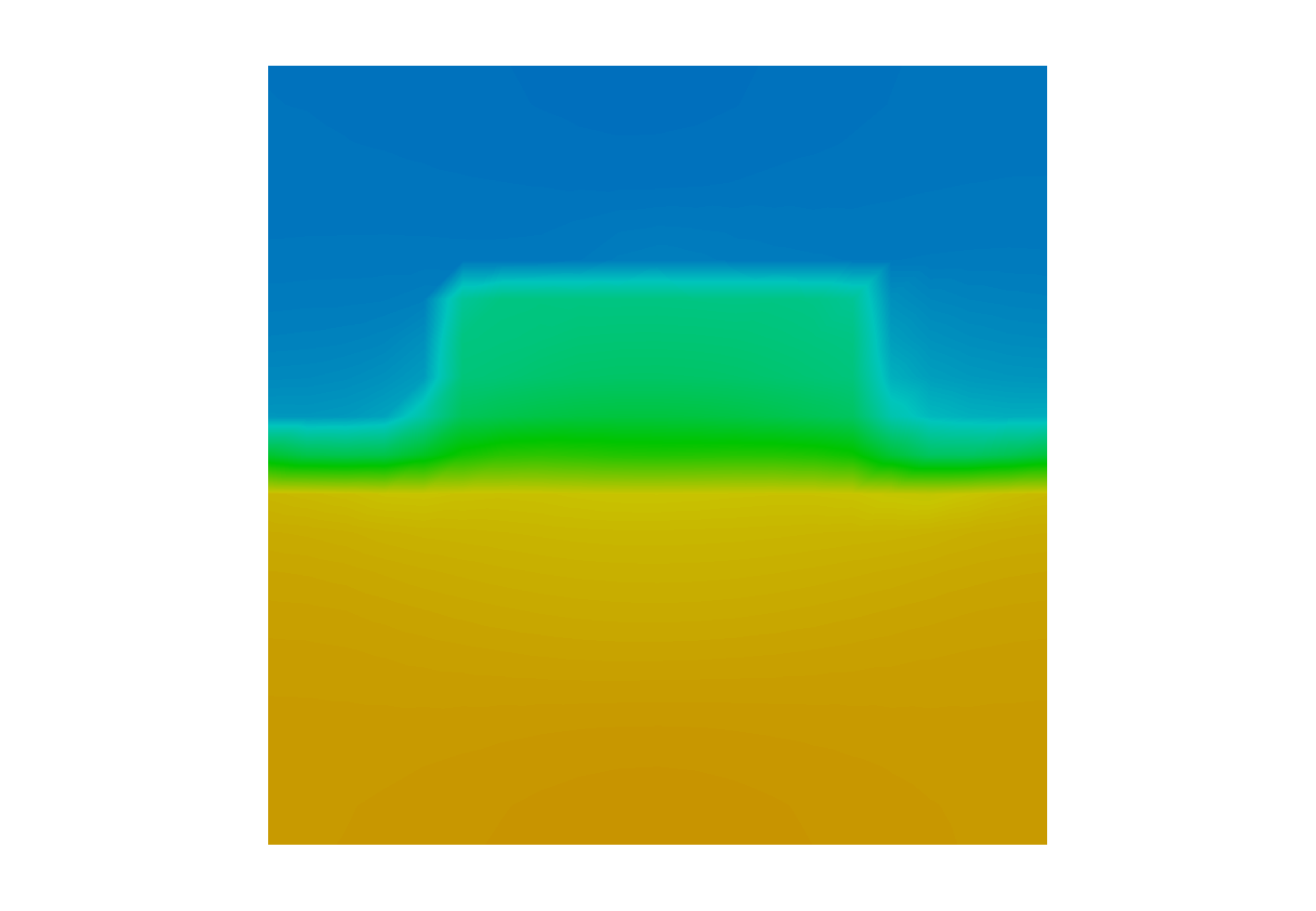}
\\
\rotatebox[origin=l]{90}{\hspace{0.35in}(b) $\alpha$}
&
\includegraphics[height=\wca, trim=290 70 290 70, clip=true]
{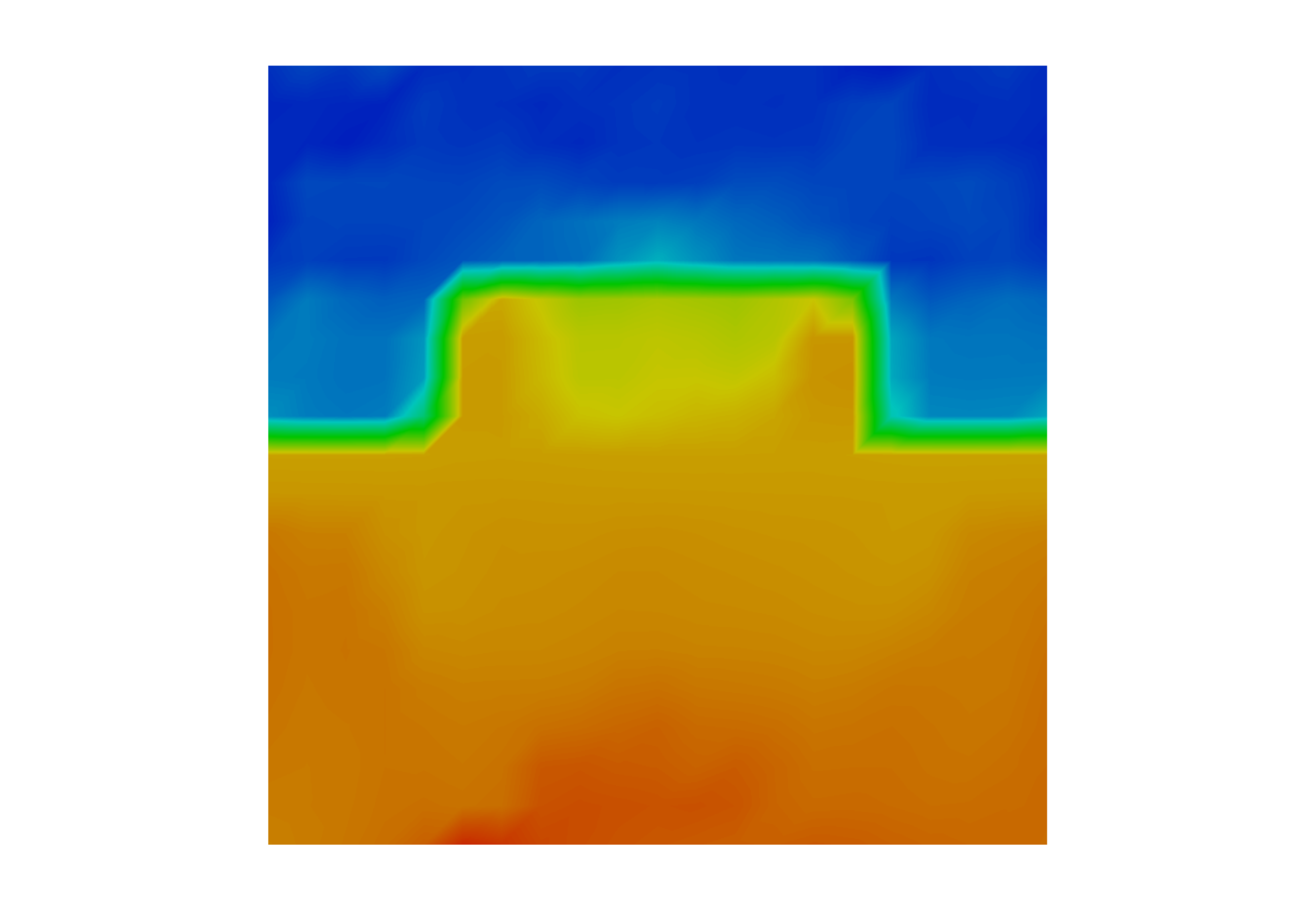} 
& 
\includegraphics[height=\wca, trim=310 70 310 70, clip=true]
{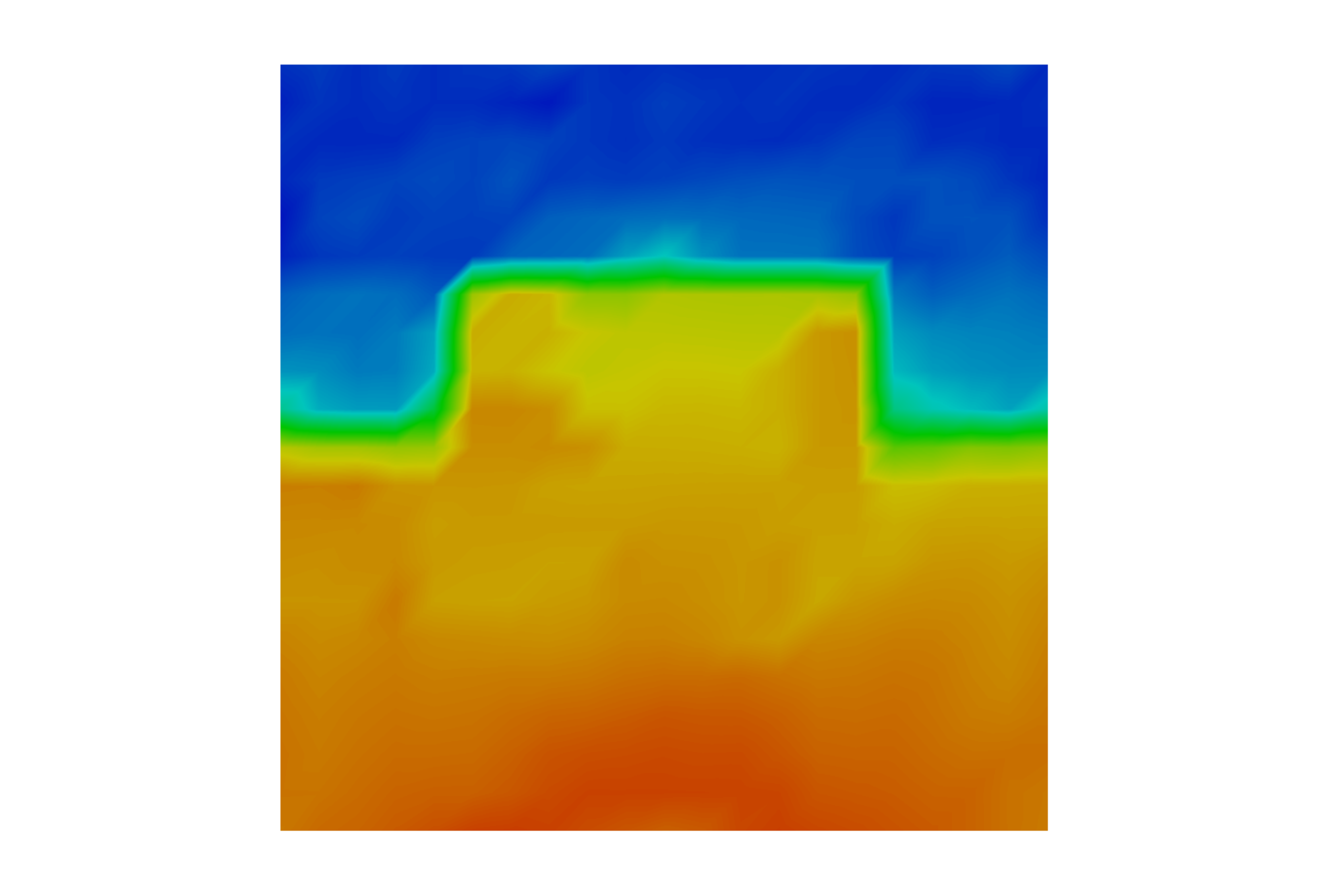} 
&
\includegraphics[height=\wca, trim=290 70 290 70, clip=true]
{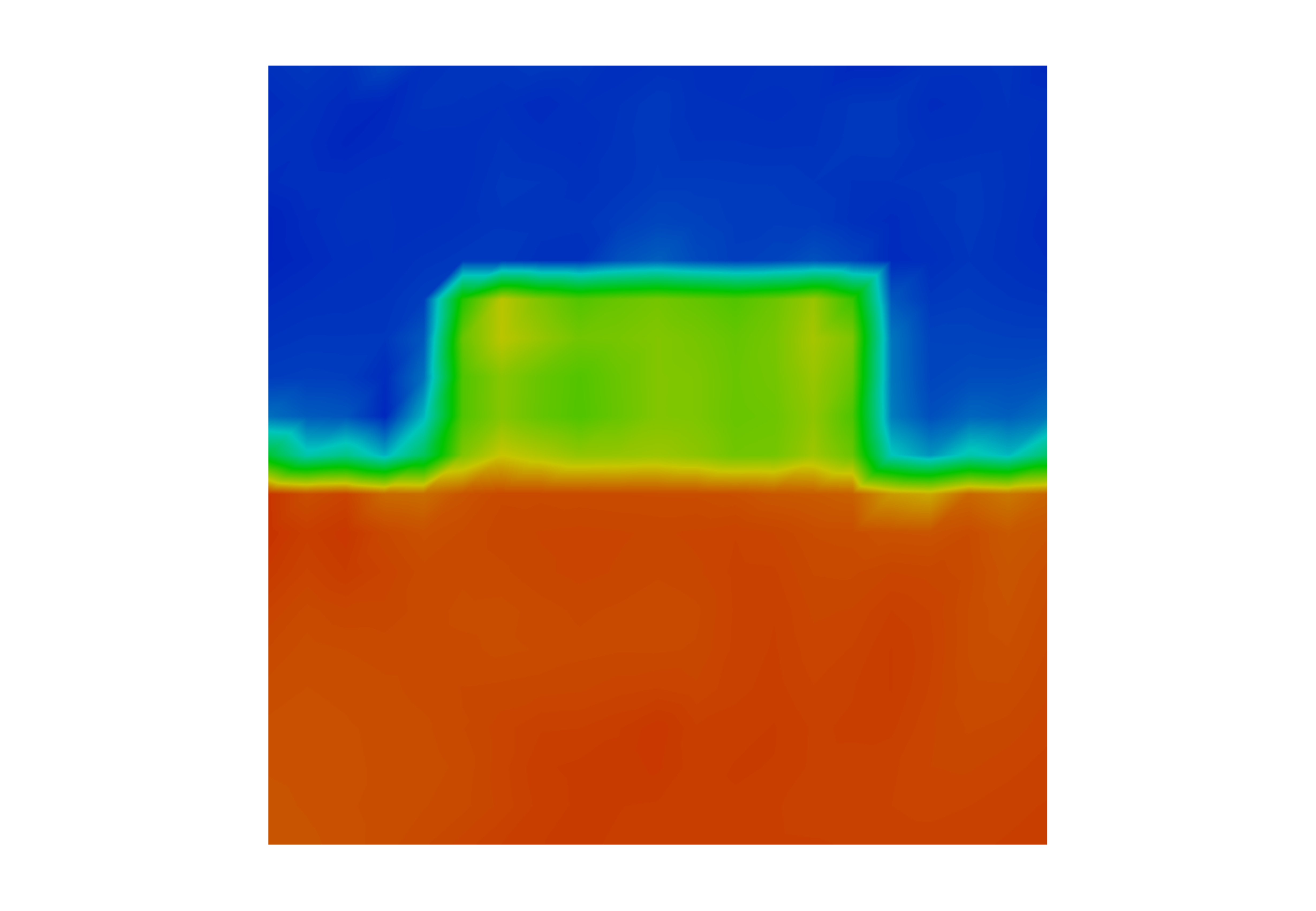}
& 
\includegraphics[height=\wca, trim=290 70 290 70, clip=true]
{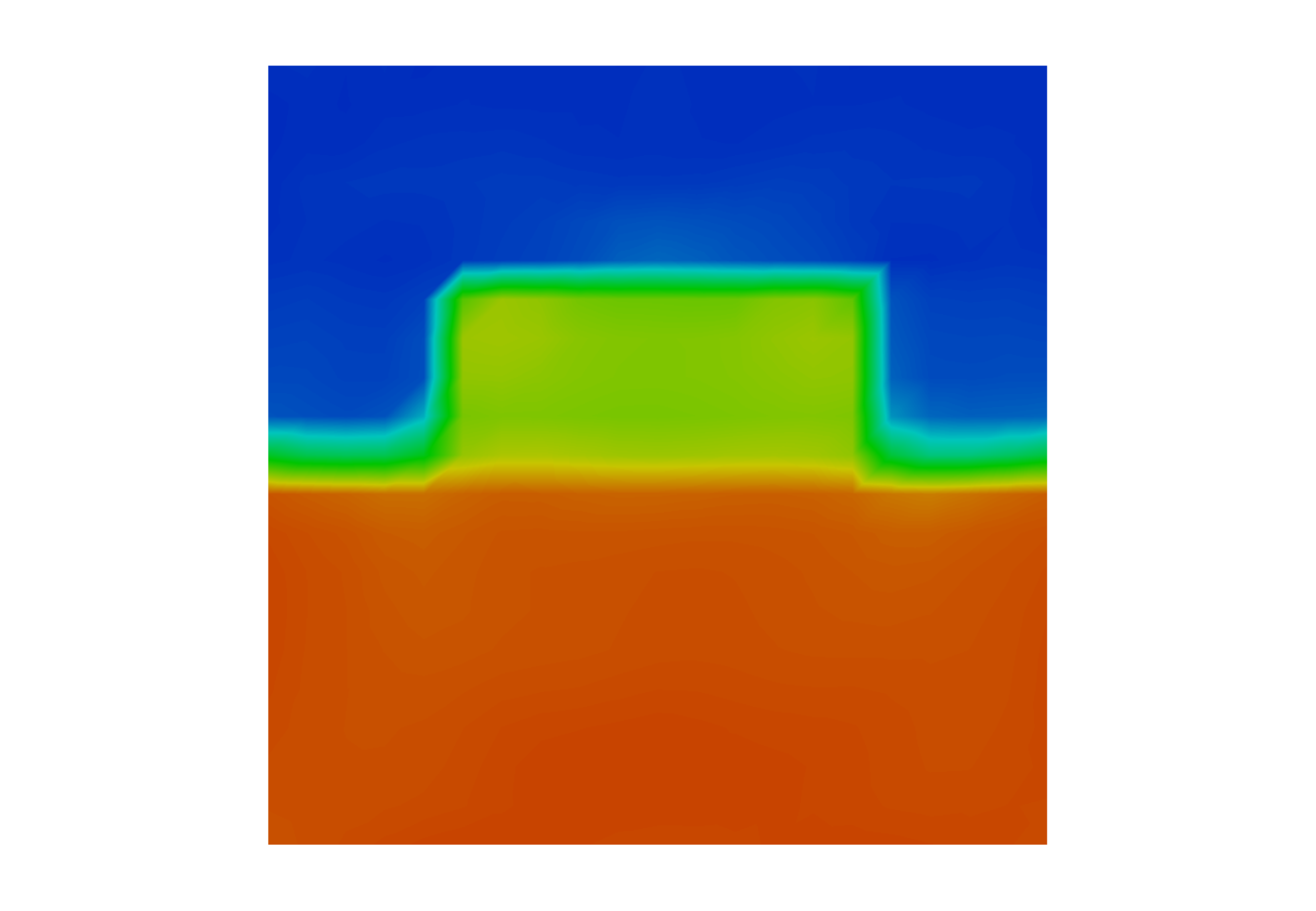}
\\
& (i) cross-gradient & (ii) norm. cross-gd &
(iii) vectorial TV & (iv) nuclear norm 
\end{tabular}
\caption{Reconstructions for the parameter fields (a) $m$
  in~\eqref{eq:poisson2} and (b) $\alpha$ in~\eqref{eq:acousticinv2},
  obtained by solving a joint inverse problem with seismic
  source~$f_\alpha$ of frequency 4~Hz, and regularized with (i) the
  cross-gradient ($\gamma=8 \cdotp 10^{-7}$) combined with two TV
  regularizations, (ii) the normalized cross-gradient ($ \gamma=8
  \cdotp 10^{-8}$ and $\varepsilon=10^{-5}$) combined with the same TV
  regularizations, (iii) the VTV joint regularization ($\gamma=4
  \cdotp 10^{-8}$ and $\varepsilon=10^{-3}$), and (iv) the nuclear
  norm joint regularization ($\gamma=5 \cdotp 10^{-7}$ and
  $\varepsilon=10^{-3}$).  The parameters for the independent TV
  regularizations are as for the independent inverse problems (see
  captions of figures~\ref{fig:poisson} and~\ref{fig:acoustic}).
  Legend is the same as in figures~\ref{fig:poisson} and
  \ref{fig:acoustic}.}
\label{fig:poissonacoustic-reconstruct}
\end{figure}

As reference, we first solve the inverse problem for the
parameters~$m$ and~$\alpha$ when \eqref{eq:poisson2}
and~\eqref{eq:acousticinv2} are solved independently.
The results for the Poisson inverse problem~\eqref{eq:acousticinv2}
are shown in figure~\ref{fig:poisson}b, where it can be seen that
the horizontal interface is well reconstructed, but the shape of the
rectangular perturbation is smeared out.
For the acoustic wave inverse problem~\eqref{eq:acousticinv2}, we show two
reconstructions in figure~\ref{fig:acoustic}, one with a source~$f_\alpha$ of
frequency 2~Hz (figure~\ref{fig:acoustic}b), and one with a source~$f_\alpha$ of
frequency 4~Hz (figure~\ref{fig:acoustic}c).
While the reconstruction at 2~Hz is excellent, the reconstruction at 4~Hz
lacks sufficient low-frequency information and appears to converge toward a
local
minimum, missing the horizontal discontinuity present in the truth parameter
field (figure~\ref{fig:acoustic}a).
The reconstructions for all four joint inverse problems, with a seismic source
$f_\alpha$ of frequency 4~Hz, are shown in
figure~\ref{fig:poissonacoustic-reconstruct}, and the corresponding values of
the relative medium misfit are given in table~\ref{tab:all-med}.

The use of the cross-gradient or its normalized
variant improves the reconstruction for the Poisson
parameter~$m$ (figures~\ref{fig:poissonacoustic-reconstruct}a, (i) and
(ii)), compared to the reconstruction from the Poisson inverse
problem~\eqref{eq:poisson2} alone (figure~\ref{fig:poisson}b).
However, neither of the cross-gradient terms brings any improvement to
the reconstruction of the acoustic wave velocity
(figures~\ref{fig:poissonacoustic-reconstruct}b, (i) and (ii)); in
particular, the reconstructions do not show the horizontal
discontinuity that was missing in the reconstruction of the acoustic
wave velocity alone (figure~\ref{fig:acoustic}c).
On the other hand, the use of either the VTV joint regularization, or the
nuclear norm joint regularization, leads to significant improvements in the
reconstruction of the acoustic wave velocity
(figures~\ref{fig:poissonacoustic-reconstruct}b, (iii) and (iv)). Both
reconstructions contain all features of the truth parameter field
(figure~\ref{fig:acoustic}a); most noticeably, the horizontal discontinuity that
was missing in the independent reconstruction (figure~\ref{fig:acoustic}c) is
now fully reconstructed.
The use of the VTV joint regularization provides only marginal improvement to
the reconstruction of the Poisson parameter~$m$, in terms of relative medium
misfit (table~\ref{tab:all-med}); however, 
the shape of the rectangular perturbation, which was smeared out in the
reconstruction from the Poisson inverse problem alone
(figure~\ref{fig:poisson}b), is clearer in
figure~\ref{fig:poissonacoustic-reconstruct}a(iii).
The reconstruction of the Poisson parameter obtained with the nuclear
norm joint regularization indicates that the optimization converged to
a local minimum. Despite all discontinuities present in the truth
parameter field (figure~\ref{fig:poisson}a) being clearly
reconstructed in figure~\ref{fig:poissonacoustic-reconstruct}a(iv),
the values of the parameters are significantly different.  Similar, or
worse, performance was observed when setting~$H_0$ to be a multiple of
the identity matrix in the BFGS solver~\cite{NocedalWright06}.
Moreover, almost identical results were obtained when solving the
Poisson-acoustic joint inverse problem, regularized by VTV, using the
BFGS method described in \ref{sec:solvernn}.  We therefore
conjecture that the poor performance of the nuclear norm joint
regularization, in the case of a multi-physics joint inverse problem,
can be attributed to the use of a gradient-based method for the
solution of the joint inverse problem.  The significant difference in
the structure of the gradients, coming from the Poisson and acoustic
wave inverse problems, dictate the use of a Newton method, which is
affine-invariant, in
order to balance the individual search directions.
This conjecture is supported by previous results found in the literature. For
instance, in the context of a joint full waveform inversion for the conductivity
and permittivity of a medium, the authors
in~\cite{LavoueBrossierMetissierEtAl14} found the reconstructions obtained using
the L-BFGS method to be highly sensitive to the scaling of the parameter fields.
The authors of~\cite{PetersHerrmann14} report similar difficulties when
employing a quasi-Newton method on a cross-well example, inverting for
compressibility and anisotropy, and study alternative formulations to remedy
this problem.

\section{Conclusion}
\label{sec:ccl}

We conducted a systematic review of regularization terms for joint
inverse problems governed by PDEs with infinite-dimensional parameter
fields.  We considered two types of joint inverse problems: (1) those
coupling several uncoupled physics forward problems via joint
regularization terms, and (2) those in which all inversion parameters
depend on the same physics. Based on a review of the literature, we
identified three joint regularization terms for this study that are
tractable for large-scale PDE constrained joint inverse problems.  The
cross-gradient is a popular choice in geophysical applications and
seeks to align level sets of the parameter fields. The normalized
cross-gradient was designed to overcome some of the potential
weaknesses of the cross-gradient term.  The vectorial total variation
is an extension of total variation regularization to joint inverse
problems, and originated from the imaging community.  In addition, we
introduced a fourth novel joint regularization term based on the
nuclear norm of a gradient matrix.
The comparison of these joint regularization terms was carried out for
three problems: (1) a joint Poisson inverse problem for which the
truth parameter fields are known to share a similar structure, (2) an
acoustic wave inverse problem in which we invert for the bulk modulus and
the density, and (3) a joint Poisson--acoustic wave inverse problem,
providing an example of  multiple physics joint inversion. 

Based on this study, we recommend use of the vectorial total variation
joint regularization. It leads to superior reconstructions in all our
examples. Moreover, we have available a scalable, efficient
primal-dual nonlinear optimization solver and Hessian preconditioner
for joint inverse problems regularized with this term
~\cite{Crestel17}.
The nuclear norm joint regularization showed encouraging results, even
leading to slightly better reconstructions than the vectorial total
variation for some examples. However, its numerical realization is
challenging since it is not twice differentiable as required by
Newton's method. 
For piecewise-homogeneous parameter fields, the cross-gradient similarity term does not improve
significantly over independent reconstructions.  In particular, it can
fail to reconstruct some edges entirely, since the cross-gradient term
vanishes at points where one parameter field is constant.
The normalized cross-gradient similarity term leads to a joint inverse
problem that is challenging to solve numerically.  Even though it
improves on the cross-gradient, the improvement is generally minimal,
and the reconstructions do not compare favorably with the ones
obtained with vectorial total variation.
Compared to the cross-gradient approaches, an additional advantage of
the VTV and nuclear norm functionals
is that they also act as
regularizations, making individual regularization functionals
unnecessary. This reduces the number of hyperparameters or regularization
weights that must be chosen (see
Table~\ref{tab:hyper}), thereby simplifying the inverse problem.

\appendix

\section{Summary of numerical optimization techniques for the solution
  of regularized inverse problems}
\label{sec:numopt}

In this section, we describe the large-scale numerical optimization
methods used for our numerical examples.
As already discussed in the introduction, the solution of
PDE-constrained optimization problems typically requires iterative
methods. These methods require first (and ideally, also second)
derivatives of the objective function with respect to the parameter
fields~\cite{NocedalWright06,Vogel02}.
These derivatives can  be computed efficiently using adjoint
methods~\cite{Reyes15, HinzePinnauUlbrichEtAl09,Troltzsch10}. In
particular, the computation of a gradient requires one solve of the
governing state equation, the solution of an adjoint equation and the
evaluation of the expression for the gradient. Moreover, adjoint
methods allow the computation of directional second derivatives by solving
two linear PDEs, one a linearization of the state equation, and the
other one a linearization of the adjoint.  Since these PDE solves
usually dominate all other
required operations, one often measures the complexity of
PDE-constrained optimization algorithms by the number of required PDE solves.
Line search and trust-region
methods are employed to globalize local
optimization methods \cite{NocedalWright06}. We restrict our
description here to the former, since we use line search methods in this paper.
In the remainder of this section, we denote the parameter
field pair by $m=(m_1,m_2)$, the objective function
by~$\mathcal J(m)$ and use upper indices to denote iteration numbers.

\subsection{Line-search Newton-CG for cross-gradient and VTV regularizations}
In the $k$-{th} iteration, we update the medium
parameters~$m^{(k)}$ along a search direction~$p^{(k)}$ by computing
$m^{(k+1)} = m^{(k)} + \alpha^{(k)} p^{(k)}$ with an appropriate step
length $\alpha^{(k)}>0$. To ensure convergence, the
search direction must be a descent direction, i.e., it must satisfy
$\langle g^{(k)}, p^{(k)} \rangle < 0$, where $g^{(k)}$ is the
gradient of $\mathcal J$ with respect to $m$ evaluated at $m^{(k)}$,
and $\langle\cdot\,,\cdot\rangle$ is an appropriate inner product.  The
step length~$\alpha^{(k)}$ could be chosen to minimize the objective
functional along this search direction~$p^{(k)}$.  However, solving this
minimization problem exactly is too expensive for large-scale
applications, since a single evaluation of the objective functional
requires the solution of the state PDE, potentially multiple times
(e.g., $N_s$ times in the example of section~\ref{sec:acoustic} which
has multiple sources).  Instead, we seek an approximate
minimizer that satisfies the following Armijo condition to ensure
sufficient descent,
\begin{align} 
& \mathcal J(m^{(k)} + \alpha^{(k)} p^{(k)}) \leq \mathcal J(m^{(k)}) + c_1 \alpha^{(k)}
  \langle g^{(k)},p^{(k)}\rangle,  \label{eq:Wolfe1} %
\end{align} 
with $0 < c_1 < 1$.
To ensure sufficiently large step lengths, we use
backtracking line
search~\cite{NocedalWright06} to
find a step length that satisfies \eqref{eq:Wolfe1}. That is, the step
length is computed by starting from an initial
guess~$\alpha_0^{(k)}>0$, and is reduced until the sufficient descent
condition~\eqref{eq:Wolfe1} is satisfied. When computing the search
direction for a Newton-type method (see next paragraph), we use
$\alpha_0^{(k)}=1$, since this is guaranteed to be a successful step
length in a neighborhood of a minimizer~\cite{NocedalWright06}.

The choice of good search directions is crucial in PDE-constrained
optimization. In the steepest descent method, one chooses the search
direction as the negative gradient, i.e., $p^{(k)}
=-g^{(k)}$. Unfortunately, the resulting algorithm usually converges
slowly in the presence of stretched contour lines of the objective $\mathcal
J$, a consequence of the typical ill-posedness of inverse problems.
The Newton direction is given by the solution of the linear system $H(m^{(k)})
p^{(k)} = -g^{(k)}$, where $H(m^{(k)})$ is the Hessian, i.e., the
second derivative of $\mathcal J$, evaluated at $m^{(k)}$. The
direction $p^{(k)}$ arising as solution of this equation is a
descent direction only if the Hessian is positive definite, which may not
be the case, in particular far away from the minimizer.
When the Hessian is indefinite, one
solution is to replace the Hessian with a positive definite
approximation, a common choice being the Gauss-Newton
Hessian~\cite{NocedalWright06}. This approximation is obtained
by setting the adjoint variables to zero in the
computation of the Hessian.  Another option is to retain the full
Hessian but solve the Newton system approximately, in a way that
guarantees the computed solution to be a descent direction.  Since for
large-scale problems exactly constructing the Hessian is infeasible,
we solve the Newton system using the conjugate gradient~(CG)
method. This requires only Hessian-vector products as provided by the
adjoint method. CG is a Krylov subspace
iterative method that solves linear systems of equations for symmetric
positive definite matrices. When applied to a Newton system
with an indefinite Hessian it will eventually encounter a negative
curvature direction and fail.
However, one can use the previous iterate before the iteration failed
as the search direction, which is guaranteed to be a descent
direction~\cite{NocedalWright06}.  Additionally, regardless of the
definiteness of the Hessian, we terminate the CG iterations prematurely
to avoid over-solving, that is, we solve the Newton system with a
coarse termination tolerance, thus applying just a few iterations of
the CG method~\cite{DemboEisenstatSteihaug82}.  As the optimization
iteration converge, the tolerance is gradually decreased to allow
increasingly accurate computation of the Newton search direction in
order 
to benefit from the fast local convergence properties of the Newton
method. For our experiments with cross-gradient regularization, its
normalized version and the vector-TV regularization, we use the
Newton-CG method with backtracking line search described above. For
the nuclear norm regularization, we do not use directional second
derivatives, but instead approximate the Newton direction using a BFGS
method, as summarized next.

\subsection{BFGS method for nuclear norm regularization}
\label{sec:solvernn}

To solve joint inverse problems regularized with the nuclear norm
joint regularization (section~\ref{sec:nn}), we use a BFGS
quasi-Newton method with damped update~\cite{NocedalWright06}.
That is,
we find the search direction~$p^{(k)}$ by computing $p^{(k)} = -
B^{(k)} g^{(k)}$, where $g^{(k)}$ is again the gradient of the objective
function
and $B^{(k)}$ is a positive definite
approximation of the inverse of the Hessian.
This approximation is updated at each iteration with the rank-2 update
\begin{equation} \label{eq:bfgsd}
B^{(k+1)} = (I-\rho_k r^{(k)} (y^{(k)})^T) B^{(k)} (I-\rho_k y^{(k)} (r^{(k)})^T) + \rho_k r^{(k)} (r^{(k)})^T,
\end{equation}
where $y^{(k)}$ is the difference between the gradient at steps~$k+1$
and~$k$, $\rho_k \coloneqq 1/(y^{(k)})^T r^{(k)}$, and $r^{(k)}$ is the damped form
of~$s^{(k)}$, the difference between the parameter at steps~$k+1$ and~$k$,
and is defined as $r^{(k)} \coloneqq \theta_k s^{(k)} + (1-\theta_k)B^{(k)}
y^{(k)}$, with
\[ \theta_k \coloneqq \left\{ \begin{aligned}
 1 , \quad & \text{ if } (s^{(k)})^T y^{(k)} \geq \alpha (y^{(k)})^T B^{(k)} y^{(k)}, \\
 \frac{ (1-\alpha) (y^{(k)})^T B^{(k)} y^{(k)}}{(y^{(k)})^T B^{(k)} y^{(k)} - (s^{(k)})^T y^{(k)}} , 
\quad & \text{
otherwise.}
\end{aligned} \right. \]
The classical BFGS method requires the curvature condition $(s^{(k)})^T
y^{(k)} > 0$ to be satisfied at all steps.  This condition is
necessary to maintain positive definiteness of $B^{(k)}$ for all $k$.
However, the curvature condition can be guaranteed to be satisfied
only when the
objective function is strictly convex, which is typically not the case
for
nonlinear inverse problems.  Using a damped update allows us to apply
a backtracking line search, while avoiding skipping some updates of
$B^{(k)}$ entirely.  In our numerical experiments, we found that
$\alpha = 0.2$ worked well.  The BFGS formula~\eqref{eq:bfgsd}
requires the initialization $B^{(0)}$.  BFGS-type methods perform
well when the difference between the initial Hessian approximation and
the true Hessian is a compact
operator~\cite{Griewank87}. 
Thus, we take $B^{(0)}$ as the inverse of the
Hessian of the regularization.  This quantity is not
available for the nuclear norm joint regularization.  However, VTV and
the nuclear norm joint regularization come from the same family of
joint regularizations, differing only by the matrix norm
employed~\cite{Holt14}. Since matrix norms are equivalent in
finite dimensions, we set $B^{(0)}$ to the inverse of the Hessian of the
VTV joint regularization at the parameter $m^{(k)}$. 

\section{Number of hyperparameters for each joint regularization}

\begin{table}[h!]
\centering
\caption{Number of hyperparameters for a joint inverse problem with 2
  parameter fields.}
\begin{tabular}{ll|cccc}
 & & \multicolumn{4}{c}{joint regularization} \\ \hline
 & & cross-grad& n-cross-grad & vectorial TV & nuclear norm \\ 
\multicolumn{2}{l|}{$\gamma$} & \multicolumn{4}{c}{} \\
 & TV & 2 & 2 & -- & -- \\
 & joint & 1 & 1 & 1 & 1 \\
\multicolumn{2}{l|}{$\varepsilon$} & \multicolumn{4}{c}{} \\
 & TV & 1 & 1 & -- & -- \\
 & joint & 1 & 1 & 1 & 1 \\ \hline
\multicolumn{2}{l|}{total} & 5 & 5 & 2 & 2
\end{tabular}
\label{tab:hyper}
\end{table}

\section{Table of relative medium misfits for examples}
In table~\ref{tab:all-med}, the relative misfits for the examples
presented in section~\ref{sec:numerics} are summarized.
\begin{table}[h!]
\centering
\caption{Relative medium misfits (in $L^2$-norm) for the
examples in section~\ref{sec:numerics}.}
\begin{tabular}{lrr|rr|rr|rr}\hline
& \multicolumn{2}{c|}{sec.~\ref{sec:coincide}}
& \multicolumn{2}{c|}{sec.~\ref{sec:coincide2}} 
& \multicolumn{2}{c|}{sec.~\ref{sec:acoustic}} 
& \multicolumn{2}{c}{sec.~\ref{sec:poissonacoustic}} 
\\
& $m_1$ & $m_2$ & $m_1$ & $m_2$ & $\alpha$ & $\beta$ 
& $m$ & $\alpha$ \\ \hline
independent & 23.2\% & 5.1\% & 46.9\% & 5.1\% 
& 2.8\% & 0.8\% 
& 9.0\% & 9.9\% \\
& & & & & & \\
cross-grad & 22.3\% & 5.2\%            & 46.1\% & 5.6\% 
& 3.1\% & 0.7\% 
& 4.9\% & 11.0\% \\
n-cross-grad & 21.2\% & 5.0\% & 46.7\% & 5.0\% 
& 2.5\% & 0.4\% 
& 4.9\% & 10.7\% \\
vectorial TV & 20.2\% & 5.1\% & 41.1\% & 5.2\% 
& 2.4\% & 0.2\% 
& 8.9\% & 3.3\% \\
nuclear norm & 20.2\% & 4.8\%              & 40.8\%  & 5.0\%
& 2.4\% & 0.2\%
& 20.6\% & 4.5\% \\ \hline
\end{tabular}
\label{tab:all-med}
\end{table}

\ack The authors would like to thank David Keyes (KAUST) and George
Turkiyyah (AUB) for very helpful discussions that inspired this
work. They would also like to thank Sergey Fomel (UT-Austin) for
referring them to~\cite{ManukyanMaurerNuber16,LiLiangAbubakarEtAl13},
Jan Modersitzki (L\"ubeck) for directing their attention
to~\cite{HaberModersitzki06b}, Nick Alger (UT-Austin) for help with
figure~\ref{fig:levelsets}, and two anonymous referees whose comments
helped significantly to improve the manuscript.  This work was
partially supported by AFOSR grant FA9550-17-1-0190, DOE grant
DE-SC0009286, KAUST award OSR-2016-CCF-2596, and NSF grants
ACI-1550593, DMS-1723211, CBET-1507009, and CBET-1508713.

\section*{References}
\bibliographystyle{unsrt}
\bibliography{ccgo}

\end{document}